\documentclass[preprint,11pt,oneside]{elsarticle}

\usepackage{amsfonts}
\usepackage{amsthm}
\usepackage{amsmath}
\usepackage{amssymb}
\usepackage{mathtools}

\usepackage{url}
\usepackage{multicol}
\usepackage{epsfig}
\usepackage[percent]{overpic}
\usepackage{graphics}
\usepackage{psfrag}
\usepackage{verbatim}
\usepackage{xcolor}
\usepackage{tabularx}
\usepackage{booktabs}
\usepackage[all]{xy}
\usepackage{makeidx}
\usepackage{enumerate}
\usepackage[letterpaper,margin=2.7cm]{geometry}
\usepackage{enumitem}
\usepackage{tikz}
\usepackage{caption}
\usepackage{subcaption}
\usepackage{float}
\newfloat{Algorithm}{tbp}{lop}
\newfloat{Listing}{tbp}{lop}
\usepackage[boxed]{algorithm2e}
\usepackage{listings}
\definecolor{mygreen}{RGB}{28,172,0} 
\definecolor{mylilas}{RGB}{170,55,241}
\SetKwComment{Comment}{\ \ \%\ }{}
\SetCommentSty{mycommfont}
\usepackage{placeins}

\theoremstyle{plain}

\theoremstyle{definition}

\newtheoremstyle{myremark}
  {3pt}
  {3pt}
  {\small \rmfamily}
  {5pt}
  {\rmfamily}
  {:}
  {.5em}
  {}

\theoremstyle{myremark}

\DeclareMathAlphabet{\pazocal}{OMS}{zplm}{m}{n}

\usepackage{mathrsfs}

\def\txtd{{\textnormal{d}}}
\def\txte{{\textnormal{e}}}

\newcommand{\eqdef}{\mathrel{\mathop:}=}

\newcommand{\be}{\begin{equation}}
\newcommand{\ee}{\end{equation}}
\newcommand{\benn}{\begin{equation*}}
\newcommand{\eenn}{\end{equation*}}
\newcommand{\bea}{\begin{eqnarray}}
\newcommand{\eea}{\end{eqnarray}}
\newcommand{\beann}{\begin{eqnarray*}}
\newcommand{\eeann}{\end{eqnarray*}}

\newcommand{\myendex}{$\blacklozenge$\end{ex}}
\newcommand{\myendexerc}{$\lozenge$\end{exerc}}
\newcommand{\myendpexerc}{$\lozenge$\end{pexerc}}

\newcommand{\commentout}[1]{}

\usepackage{graphicx}

\definecolor{blun}{RGB}{28,172,0}
\usepackage[bookmarks,colorlinks,linkcolor=blun,citecolor=blun,urlcolor=blun,pdfauthor={EsthandKuehnSoresina}]{hyperref}

\begin{document}

\title{Numerical continuation for fractional {PDEs}:\\ sharp teeth and bloated snakes}

\author[ubi]{No\'emie Ehstand}\ead{n.ehstand@ifisc.uib-csic.es}
\author[tum]{Christian Kuehn}\ead{ckuehn@ma.tum.de}
\author[graz]{Cinzia Soresina\corref{cor3}}\ead{cinzia.soresina@uni-graz.at}

\cortext[cor3]{Corresponding author}
\address[ubi]{IFISC, University of the Balearic Islands and Spanish National Research Council,\\ Campus Universitat de les Illes Balears, E-07122 Palma, Mallorca, Spain}
\address[tum]{Zentrum Mathematik, Technische Universit\"at M\"unchen,\\ Boltzmannstr. 3, 85748 Garching bei M\"unchen, Germany}
\address[graz]{Institut f\"ur Mathematik und Wissenschaftliches Rechnen, Heinrichstr. 36, 8010 Graz, Austria}
\journal{Communications in Nonlinear Science and Numerical Simulation}
\begin{abstract}
Partial differential equations (PDEs) involving fractional Laplace operators have been increasingly used to model non-local diffusion processes and are actively investigated using both analytical and numerical approaches. 
The purpose of this work is to study the effects of the spectral fractional Laplacian on the bifurcation structure of reaction--diffusion systems on bounded domains. In order to do this we use advanced numerical continuation techniques to compute the solution branches. 
Since current available continuation packages only support systems involving the standard Laplacian, we first extend the \texttt{pde2path} software to treat fractional PDEs (in the spectral definition). The new capabilities are then applied to the study of the Allen--Cahn equation, the Swift--Hohenberg equation and the Schnakenberg system (in which the standard Laplacian is replaced by the spectral fractional Laplacian). In particular, we investigate the changes in snaking bifurcation diagrams and in the spatial structure of non-trivial steady states upon variation of the order of the fractional Laplacian. Our results show that the fractional order induces significant qualitative and quantitative changes in the overall bifurcation structures, of which some are shared by the three systems. This contributes to a better understanding of the effects of fractional diffusion in generic reaction--diffusion systems.
\end{abstract}

\begin{keyword}
spectral fractional Laplacian \sep \texttt{pde2path} \sep Balakrishnan formula \sep numerical continuation \sep homoclinic snaking
\end{keyword}

\maketitle


\section{Introduction} \label{sec:intro}
In recent years increasing research attention has been devoted to the study of partial differential equations (PDEs) involving so called fractional Laplace operators, which generalize the notion of derivative to non-integer orders. Their study is not only interesting from a purely mathematical point of view, but has proven extremely useful for modeling super-diffusion processes, which naturally appear in many applications in physics, probability, biology, ecology, medicine and economics~\cite{metzler_random_2000}. Examples include models for complex phase transitions~\cite{sire_fractional_2008}, for chemical transport in heterogeneous aquifers~\cite{Adams1992}, for pattern formation in coral reefs~\cite{Somathilake2018}, for larval growth and recruitment in a turbulent environment~\cite{burrow2008levy}, for foraging of marine predators~\cite{viswanathan2010ecology}, for cardiac electrical propagation~\cite{bueno2014fractional} and for mechanisms in options pricing~\cite{Levendorski2004}. In particular, super-diffusion phenomena are associated to L{\'e}vy flights, which, in contrast to Brownian motion, are characterized by the presence of long jumps and non-local interactions \cite{garbaczewski2019fractional}.

From a mathematical point of view, the fractional Laplacian does not have a unique definition. While the different definitions are equivalent on unbounded domains (e.g. Fourier, singular integral, Balakrishnan formula or harmonic extension)~\cite{Kwanicki2017}, this is usually not the case on bounded domains. In fact, on bounded domains, boundary conditions must be incorporated in these representations in mathematically distinct ways, leading to non-equivalent definitions (e.g. spectral and Riesz definitions).  For a detailed review we refer to~\cite{di2012hitchhiker, Lischke_what_2018}.  As a consequence the choice of the most suitable definition for a specific application is a crucial point, and there is still no agreement in the literature. The study of fractional PDEs, and in general of nonlocal phenomena, is indeed an actual and active research topic.  Nevertheless, it is widely accepted that the fractional Laplacian constitutes the counterpart in the theory of non-local operators of the standard Laplacian in the classical theory of partial differential equations. 
Despite the difference of the two operators, most of the relevant questions associated to the Laplacian (e.g.~comparison principles, linear and nonlinear boundary value problems, and regularity results) have an equivalent for the fractional Laplacian, and many techniques for solving partial differential equations also apply to fractional PDEs. \\

From a dynamical systems perspective, one of the most widely studied class of PDEs are standard reaction--diffusion problems
$$\partial_t u=D\Delta u+F(u,\mu),$$
where~$t\in \mathbb{R}^+$, $x\in \Omega\subseteq\mathbb{R}^n$ are the time and space variables,~$u(x,t)\in \mathbb{R}^N$ is the vector of unknowns at $(x,t)$, the matrix~$D\in \mathbb{R}^{N\times N}$ contains the diffusion coefficients,~$\Delta$ is the classical Laplacian operator, $F:\mathbb{R}^N\times \mathbb{R} \rightarrow \mathbb{R}^N$ corresponds to non-linear reaction terms, and~$\mu\in \mathbb{R}$ is a parameter. In many cases a natural starting point in the study of such equations is to consider the non-linear elliptic problem
$$0=D\Delta u+F(u,\mu),$$
which yields steady states. Of particular interest in applications is the question of how the steady states change under parameter variation. Several methods (Lyapunov--Schmidt reduction, center manifolds, amplitude equations, Turing instability) exist to address this question analytically near trivial branches of homogeneous steady states~\cite{KuehnBookPDE}. In addition to theoretical methods and criteria, the overall bifurcation structure of steady states reveals the behavior of solutions far from homogeneous solutions or the critical points. For classical reaction--diffusion systems on bounded domains (and~$n=1,2,3$), the overall structure can be numerically computed with, for instance, \texttt{pde2path}~\cite{uecker_pde2path_2014}. This package is an advanced continuation/bifurcation software based on the FEM discretization of the stationary elliptic problem exploiting the package \texttt{OOPDE}~\cite{OOPDE} for the FEM discretization. Since \texttt{pde2path} allows branch point continuation and Hopf point continuation, continuation of relative equilibria (e.g., traveling waves and rotating waves), branch switching from periodic orbits (pitchfork bifurcation of periodic orbits and period doubling), it is quite flexible. Recently, it has also been used to treat cross-diffusion systems~\cite{CKCS, MBCKCS}. Yet, non-standard diffusion problems are not within the classical realm of \texttt{pde2path} applications.\\

In this framework the aim of the paper is two-fold. First, we want to extend the continuation software \texttt{pde2path} which has been extensively used for classical non-linear reaction--diffusion equations, to treat fractional equations, i.e. equations involving the fractional Laplacian. Then, we test the new capabilities of the software onto three important benchmark problems: namely the Allen--Cahn equation~\cite{allen1979microscopic}, the Swift--Hohenberg equation~\cite{swift1977hydrodynamic, CrossHohenberg} and the Schnakenberg system~\cite{schnakenberg1979simple} on bounded domains in a non-local setting. All three equations have theoretical interest, concrete applications, and are currently widely used in various communities working on PDEs. The purpose is to qualitatively and quantitatively explore the features of these models in the fractional setting and to better understand the effect of fractional diffusion on the steady state bifurcation structure of generic fractional reaction--diffusion systems.\\

Let~$\Omega\subset \mathbb{R}^n$ ($n=1,2,3$) be a bounded domain with boundary~$\partial \Omega$ and~$s\in(0,1)$ be the fractional order. A generic fractional reaction--diffusion system can be formulated as
\begin{equation}
    \partial_t u= -D(-\Delta)^s {u} + F({u};\mu)=:D\Delta^s {u} + F({u};\mu),
    \label{eq:reaction-diffusion-ss-frac}
\end{equation}
endowed with suitable boundary conditions and where we have introduced the shorthand notation
\begin{equation*}
\Delta^s:=-(-\Delta)^s
\end{equation*}
to be used throughout this work. In this paper we consider the spectral definition of the fractional Laplacian~$\Delta^s$~\cite{capella2010regularity}, that is, for any~$\omega \in C^{\infty}(\Omega)$,
\begin{equation}
\label{eq:spectral_definition}
\Delta^s \omega(x) = -\sum_{j=1}^{\infty} (-\lambda_j)^s \omega_j \phi_j(x), \qquad \omega_j:=\langle \omega,\phi_j \rangle_{L^2_\Omega}=\int_{\Omega}\omega \phi_j~\txtd x
\end{equation}
where~$\lambda_j$ and~$\phi_j$ the eigenvalues and the eigenfunctions, respectively, of the standard Laplacian~$\Delta$ on~$\Omega$, i.e.,~the solutions of the eigenvalue problem~$\Delta \phi_j=\lambda_j\phi_j$ with specific boundary conditions.

In order to numerically treat fractional PDEs and to embed it naturally into the numerical continuation approaches based upon FEM, we provide a FEM discretization of the fractional operator. The approach we employ to represent~\eqref{eq:spectral_definition} without a spectral Galerkin decomposition is to use the Balakrishnan representation formula \cite{balakrishnan1960}, namely for $s\in(0, 1)$
\begin{equation}\label{eq:Balakrishnan_formula-bis}
\Delta^s u(x)=-\frac{\sin(s\pi)}{\pi}\int_0^\infty \Delta(\xi I-\Delta)^{-1} u(x)\xi^{s-1}\txtd\xi.
\end{equation}
The numerical approximation of~\eqref{eq:Balakrishnan_formula-bis} is based on a sinc quadrature approximation of the involved integral with respect to $\xi$ and a discretization of the operator $\xi I-\Delta$ using the finite element method in space \cite{bonito2018numerical}. This strategy has been successfully applied recently in other contexts such as optimal control problems~\cite{dohr2019fem}. Here we integrate this discretization within the numerical continuation.

Then, using the new capabilities of the continuation software, we study the effect of the fractional derivative order on the steady state bifurcation structure of three well known PDEs in which the standard Laplacian is replaced by a spectral fractional Laplacian:
\begin{itemize}[leftmargin=*]
\item[-] The fractional Allen--Cahn equation with cubic--quintic nonlinearities
\begin{equation}\label{eq:AC-frac}
\partial_t u = \Delta^s u   + \mu u + u^3 - \gamma u^5, \quad \quad u\in \mathbb{R},\ \ \mu, \gamma \in \mathbb{R},
\end{equation}
will be considered on a bounded interval $\Omega \subset \mathbb{R}$ taken with homogenous Dirichlet boundary conditions. The model with the standard Laplacian, proposed in \cite{allen1979microscopic}, is a well-studied equation for various polynomial nonlinearities~\cite{alfaro2008singular, alama1997stationary, rabinowitz2003mixed, ward1996metastable}. In particular, the classical Allen--Cahn PDE ($s=1$) with a cubic nonlinearity is also known as the (real) Ginzburg--Landau PDE \cite{aranson2002world}, an amplitude equation or normal form for bifurcations from homogeneous states \cite{mielke2002ginzburg,schneider1996validity}. Moreover, it is also common to consider an additional quintic term added to the Ginzburg--Landau PDE \cite{kapitula1998instability,kuehn2015numerical}. Different versions of the fractional Allen--Cahn have been recently investigated, see for instance~\cite{achleitner2020metastable, akagi2016fractional, bueno2014fourier}.
\item[-] The fractional Swift--Hohenberg equation with competing cubic--quintic nonlinearities
\begin{equation}\label{eq:SH-frac}
\partial_t u =  -(1+\Delta^s)^2 u +\mu u+ \nu u^3 - u^5, \quad \quad u\in \mathbb{R},\ \ \mu, \nu \in \mathbb{R},\ \ \nu>0,
\end{equation}
will also be studied on a bounded interval $\Omega \subset \mathbb{R}$, again taken with homogeneous Dirichlet boundary conditions. The classical Swift--Hohenberg equation ($s=1$) \cite{swift1977hydrodynamic} is a widely studied model in pattern dynamics~\cite{CrossHohenberg,Schneider5,Thieleetal}. It is also a standard test case for deriving reduced amplitude/modulation equations~\cite{KuehnBookPDE,SchneiderUecker,KirrmannSchneiderMielke,ColletEckmann1,vanHarten}. Furthermore, it was found that overall bifurcation diagrams of the Swift--Hohenberg equation exhibit a process called (homoclinic) snaking~\cite{BURKE2007681,avitabile2010snake,houghton2011swift}. Recently, also the Swift--Hohenberg equation with nonlocal reaction terms has gained particular attention~\cite{KuehnThrom,morgan2014SH} but results for bifurcations in the space-fractional Swift--Hohenberg equation do not seem to be available up to now.
\item[-] The fractional Schnakenberg system is given by
\begin{equation}\label{eq:Schnak-frac}
\partial_t\begin{pmatrix}  u_1 \\ u_2 \end{pmatrix} = \begin{pmatrix} 1 & 0 \\ 0 & d \end{pmatrix} \begin{pmatrix} \Delta^s u_1\\ \Delta^s u_2 \end{pmatrix} + F(u_1,u_2;\mu), \quad \quad  u_1,u_2\in \mathbb{R}, \ \ \mu, d\in \mathbb{R}, \ \ d>1,
\end{equation}
where the reaction part is
\begin{equation}\label{eq:F_m}
F(u_1,u_2;\mu) =  \begin{pmatrix} -u_1 + u_1^2u_2 \\ \mu - u_1^2u_2  \end{pmatrix}
+ \sigma\left(u_1-\frac{1}{u_2}\right)^2\begin{pmatrix} 1 \\ -1\end{pmatrix},
\quad \quad \sigma\in \mathbb{R}.
\end{equation}
As before, we consider a bounded interval $\Omega \in \mathbb{R}$ taken with homogeneous Neumann boundary conditions. Note that with $\sigma=0$ the reaction part~\eqref{eq:F_m} reduces to the classical formulation of the Schnakenberg system~\cite{schnakenberg1979simple}, which is one of the prototype reaction--diffusion systems exhibiting Turing patterns~\cite{murrayIIspatial}. A nonzero $\sigma$ modifies the primary bifurcation from super- to subcritical, leading to to the appearance of snaking branches between periodic branches~\cite{Uecker2014,uecker_pde2path_2014}.
\end{itemize}

As indicated above, the classical Swift--Hohenberg equation and the Schnakenberg system both exhibit snaking behavior~\cite{BURKE2007681,avitabile2010snake,houghton2011swift,Becketal,KnoblochLloydSandstedeWagenknecht,McCallaSandstede,MercaderBatisteAlonsoKnobloch,Uecker2014}: upon variation of the main bifurcation parameter, a non-trivial branch of solutions emerges from a branch point and later on undergoes several fold bifurcations leading to a snaking-type structure of the bifurcation diagram and to multi-stability. Since solutions on these branches are often important localized patterns, a detailed analysis is required. In particular, it is important to understand, how the snaking scenario changes under nonlocal influence as already carried out for nonlocal reaction terms in~\cite{morgan2014SH}. Here we show that, in the Swift--Hohenberg equation, the snaking structure gets ``bloated'' upon decreasing the fractional power $s\in(0,1)$, i.e. its width increases. In the Schnakenberg system, we observe a similar effect as well as further classes of deformations of the branches. Further, in all three studied systems, the solutions tend to sharpen in the transition layers, a ``sharpening of the teeth'' of spatially oscillating profiles.

\medskip
The paper is organized as follows. The software \texttt{pde2path} as well as the discretization and the implementation of the fractional operator are detailed in Section~\ref{sec:numerics}, followed by some preliminary results in order to validate the discretization method. In Sections \ref{sec:AC-theory}, \ref{sec:SHe-theory} and \ref{schnak} we briefly recall the standard version of each benchmark problem (the Allen--Cahn equation, the Swift--Hohenberg equation and the Schnakenberg system respectively), including the standard steady state bifurcation diagrams as well as typical solutions. Sections \ref{chapter:AC-results}, \ref{chapter:SHe-results} and \ref{chapter:Schnakenberg-system-results} are devoted to the presentation of the results obtained by exploiting the new continuation capabilities, highlighting the effects of fractional diffusion on the bifurcation structure of the considered problems. Finally, in Section \ref{sec:conclusions} some concluding remarks can be found. The Matlab code is reported in the Appendix, together with some explanations on the Matlab implementation. Some supplementary videos are available at~\cite{VideoFolder}.

\section{Numerical continuation of fractional PDEs in \texttt{pde2path}}\label{sec:numerics}

As mentioned in the introduction, we want to adapt numerical continuation in the context of the Matlab bifurcation package \texttt{pde2path} to treat fractional reaction--diffusion equations. In this section, we briefly present the software in its standard setting, that is for PDEs with standard Laplacian. Then we explain the FEM discretization of the fractional Laplacian and modifications made to the FEM discretization in \texttt{pde2path} in order to study fractional systems. We finally present numerical tests in order to validate our method. The purpose of this section is not to give a complete overview of the software for beginner users, hence we do not explain in detail the basic setup; see~\cite{dohnal2014pde2path, rademacher2018oopde, uecker_pde2path_2014, uecker2019hopf, uecker2019pattern} for complete guides and for the notation adopted in the following.

\subsection{Standard setting}
The continuation software \texttt{pde2path} can deal with the general class of PDEs of the form
\begin{equation*}
\partial_t u = -G(u,\mu) \eqdef \nabla \cdot (c\otimes \nabla u) - au + b\otimes\nabla u+f,
\end{equation*}
where $u=u(x,t)\in\mathbb{R}^N$, $x\in \Omega \subset \mathbb{R}^n$, for $n = {1,2,3}$ and $\Omega$ a bounded domain, $\mu \in \mathbb{R}^p$ is a parameter vector, and the diffusion, advection and linear tensors, denoted by $c\in \mathbb{R}^{N\times N\times n \times n},\, b\in \mathbb{R}^{N\times N\times n}$ and $\, a\in \mathbb{R}^{N\times N}$ respectively, as well as the nonlinearity $f$ can depend on $x$, $u$, $\nabla u$ and parameters. We refer to \cite{uecker_pde2path_2014} for a detailed explanation. The problem is endowed with generalized Neumann boundary conditions, also known as mixed conditions,
\begin{equation} \label{eq:BC_for_RD_pde2path}
\mathbf{n}\cdot (c \otimes \nabla u) + qu = g,
\end{equation}
where $\mathbf{n}$ is the outer normal to the boundary, $q\in\mathbb{R}^{N\times N}$ and  $g\in\mathbb{R}^{N}$. This formulation covers Neumann and Robin boundary conditions, as well as Dirichlet boundary conditions via large stiff spring pre-factors in $q$ and $g$. The software spatially discretizes the problem via the finite element method (FEM) exploiting the OOPDE toolbox \cite{OOPDE}, leading to a high-dimensional ODE problem. Note that a code extension is also available for periodic boundary conditions \cite{dohnal2014pde2path} and cross-diffusion terms \cite{CKCS}.

In order to study the stationary problem related to equation~\eqref{eq:reaction-diffusion-ss-frac}, we restrict this introduction to the case $a$ and $b$ null tensors, $c$ such that $c_{ijkh}=d_i\delta_{ij}\delta_{kh},\, i,j=1,\dots,N,\,k,h=1,\dots, n$, and $\mu\in\mathbb{R}$, leading to
\begin{equation}\label{eq:general_RD_pde2path}
G(u;\mu) = -D\Delta u-f = 0,
\end{equation}
where $D\in \mathbb{R}^{N\times N}$ is the diagonal matrix of the diffusion coefficients $d_i$.
In \texttt{pde2path} the problem~\eqref{eq:general_RD_pde2path}, together with boundary condition~\eqref{eq:BC_for_RD_pde2path}, is converted into the algebraic system
\begin{equation}\label{eq:algebraic_system_for_G}
    \mathtt{G}(\mathtt{u},\mu) = K_{tot}\mathtt{u} - F_{tot}(\mathtt{u}) = 0,
\end{equation}
where $\mathtt{u}\in \mathbb{R}^{n_pN}$ contains the nodal values of $u$ with $n_p$ mesh points. The matrix $K_{tot}\in \mathbb{R}^{n_pN\times n_pN}$ corresponds to the finite element discretization of the $D$ and $q$ terms from equations~\eqref{eq:general_RD_pde2path} and~\eqref{eq:BC_for_RD_pde2path}, while $F_{tot}\in \mathbb{R}^{n_pN}$ corresponds to the discretization of $f$ and $g$ in~\eqref{eq:general_RD_pde2path} and~\eqref{eq:BC_for_RD_pde2path}. Precisely,
\begin{equation}\label{eq:matrixtot}
K_{tot}=K+s_{BC}Q_{BC},\qquad F_{tot}=F-s_{BC}G_{BC},
\end{equation}
where $K$ is the stiffness matrix and corresponds to $D$, $s_{BC}$ is the stiff spring pre-factor used to approximate Dirichelet boundary conditions (BC) and $F,\,Q_{BC},\,G_{BC}$ corresponds to $f,\,q,\,g$ respectively. Details on the discretization for different types of boundary conditions can be found in \cite{rademacher2018oopde}.

\subsection{Fractional PDEs in \texttt{pde2path}}

In order to adapt the continuation software \texttt{pde2path} to treat fractional reaction--diffusion equations~\eqref{eq:reaction-diffusion-ss-frac} endowed with homogeneous Neumann or Dirichlet boundary conditions, we need to discretize the spectral fractional Laplacian~\eqref{eq:spectral_definition} keeping the structure of equation~\eqref{eq:algebraic_system_for_G}. We exploit the Balakrishnan formula representation \cite{balakrishnan1960, Kwanicki2017} in the version appearing in \cite[formula (4), p.260]{yosida1968functional}, reported in equation~\eqref{eq:Balakrishnan_formula-bis}.
Note that ``Balakrishnan formula" refers to both the direct fractional Laplacian representation, $\Delta^su$ as above, and the inverse fractional Laplacian representation, $\Delta^{-s}u$ see \cite{dohr2019fem, bonito2018numerical}. In~\cite{dohr2019fem, bonito2018numerical} the inverse formula has been discretized via a sinc quadrature for the integral as well as a finite element discretization in space. Here we adapt this method to the direct formula~\eqref{eq:Balakrishnan_formula-bis}.\\
First, the substitution $\xi = \txte^\eta$ ($d\xi = \txte^\eta d\eta$) in~\eqref{eq:Balakrishnan_formula-bis} leads to
\begin{equation*} 
\Delta^su(x)=\frac{\sin(s\pi)}{\pi}\int_{-\infty}^\infty (\txte^\eta I-\Delta)^{-1}\Delta u(x)\txte^{\eta s} ~\txtd\eta.
\end{equation*}
Using a truncated equally spaced quadrature (sinc quadrature \cite{lund1992sinc}) one can then approximate the integral as
\begin{equation*}
\Delta^s u(x) \approx \frac{\sin(s\pi)}{\pi} \kappa \sum_{l=-n_{-}}^{n^{+}}\txte^{\kappa ls}(\txte^{\kappa l}I-\Delta)^{-1}\Delta u(x),
\end{equation*}
where $\kappa>0$ is the quadrature step size and $n_+$, $n_-$ are positive integers suitably chosen (see \ref{app:code}). A fully discrete version of equation \eqref{eq:Balakrishnan_formula-bis} is obtained as
\begin{equation}\label{eq:quadrature_approx}
\Delta^s u(x) \approx -\frac{\sin(s\pi)}{\pi} \kappa \sum_{l=-n_{-}}^{n_{+}}\txte^{\kappa ls}v_l(x),
\end{equation}
where $v_l \in \mathbb{U}(\mathcal{\tau})$ (being $\mathbb{U}(\mathcal{\tau})$ the finite element space satisfying the boundary conditions \cite{bonito2018numerical}) is the solution of the Galerkin variational problem
\begin{equation*}
\txte^{\kappa l}\int_\Omega v_lw-\int_\Omega \Delta v_lw=-\int_\Omega \Delta u w \qquad \forall w\in \mathbb{U}(\mathcal{\tau}),
\end{equation*}
approximated via FEM as the solution of the linear system
\begin{equation}\label{eq:linsys}
(\txte^{\kappa l}M+K)\mathtt{v_l}=K\mathtt{u}.
\end{equation}
where $M$ and $K$ are the classical FEM mass and stiffness matrices respectively, and $\mathtt{v_l}$ denotes the node vector for $v_l$.
Hence, we obtained an approximation of $\Delta^s u$, which requires to solve~$n_-+n_++1$ independent linear systems.

However, to fit the \texttt{pde2path} setup (which necessitates the gradient expression of \texttt{G} in~\eqref{eq:algebraic_system_for_G}) given by $K_{tot} - \nabla F_{tot}$), we need a matrix approximation of the operator $\Delta^s$ itself. Such matrix, called $K_s$, can be built column by column by \emph{testing} the proposed approximation method on the standard basis of $\mathbb{R}^{n_p}$ of vectors $\hat{e}_i,\, 1 \leq i\leq n_p$, such that $(\hat{e}_i)_j = \delta_{ij}$ for all $1\leq j \leq n_p$. In detail, the $i$-th column of $K_s$ is obtained computing $\Delta^s \hat{e}_i$ with the proposed scheme.
Once the matrix $K_s$ has been built, we end up with an algebraic system of the form~\eqref{eq:algebraic_system_for_G} by replacing the standard stiffness matrix $K$ with $-MK_s$ in the expression of $K_{tot}$ in formula~\eqref{eq:matrixtot}.

\medskip
The complete method is presented in Algorithm \ref{alg:FractionalLaplacian}. A description of its implementation in \texttt{pde2path} can be found in \ref{app:code}. In particular, the Matlab code is given in Listing~\ref{lst:FractionalLaplacian} and the choice of the quantities $k,\,n_+,\,n_-$ as well as further implementation issues are discussed.

\begin{Algorithm}
\LinesNumbered
\begin{algorithm}[H]
\label{alg:FractionalLaplacian}
choose $\kappa$ \\
set $n_{+}$ and $n_{-}$\\
\For {$i\ =\ 1\to n_p$}{
    $t\ =\ \hat{e}_i$ \Comment{test vector}
    $z\ =\ Kt$ \\
    $col=0$\\
        \For {$l\ =\ n_{-} \to n_{+}$}{
            solve $(\txte^{\kappa l} M + K)*v\ =\ z$ for $v$\\
            $col = col + \txte^{ \kappa ls} \cdot v$\\
        }
    $col\ =\ -\dfrac{\kappa \cdot \sin(s\pi)}{\pi} \cdot col$\\
    $i^{th}$ column of $K_s$ $\leftarrow$  $col$
}
\end{algorithm}
\caption{Algorithm that allows to build the matrix $K_s$, which corresponds to a FEM discretization of the fractional Laplacian. A detailed description of the implementation of this algorithm in \texttt{pde2path} is given in \ref{app:code}.}
\end{Algorithm}


\subsection{Numerical tests} \label{sec:pre-results}
In order to validate the algorithm presented in the previous section, we first present some numerical results on the eigenvalues of the fractional Laplacian. The algorithm is then applied to a simple Poisson problem, exemplifying the effect of fractional diffusion.

\subsubsection{Spectrum of the discretized fractional Laplacian}
A first simple test consists in checking the convergence of the eigenvalues of the discretized fractional operator $K_s$ towards those of the continuous operator as the mesh-size is decreased. According to \eqref{eq:spectral_definition}, on a one-dimensional domain of length $L$, the Neumann spectral fractional Laplacian $\Delta^s$ has eigenpairs $(\phi_j, \lambda_j) = \left(\cos( {j\pi x}/{L} ), -({j \pi}/{L})^{2s}\right)$ for~$j\geq 0$.
The first $n_e$ eigenvalues $\lambda^{(h)}_j,\,j=1,\dots,n_e$ of the discretized fractional Laplacian $K_s$ can be numerically computed, and the maximum relative error can be used as a measure of the convergence
\begin{equation}
\mathrm{err_{n_e}(\lambda_j,\lambda^{(h)}_j)}:=\max_{j\leq n_e} \left( \frac{|\lambda_j-\lambda^{(h)}_j|}{|\lambda_j|} \right).
\label{eq:max-rel-er}
\end{equation}

\begin{figure}
    \centering
    \begin{overpic}[width=0.45\textwidth]{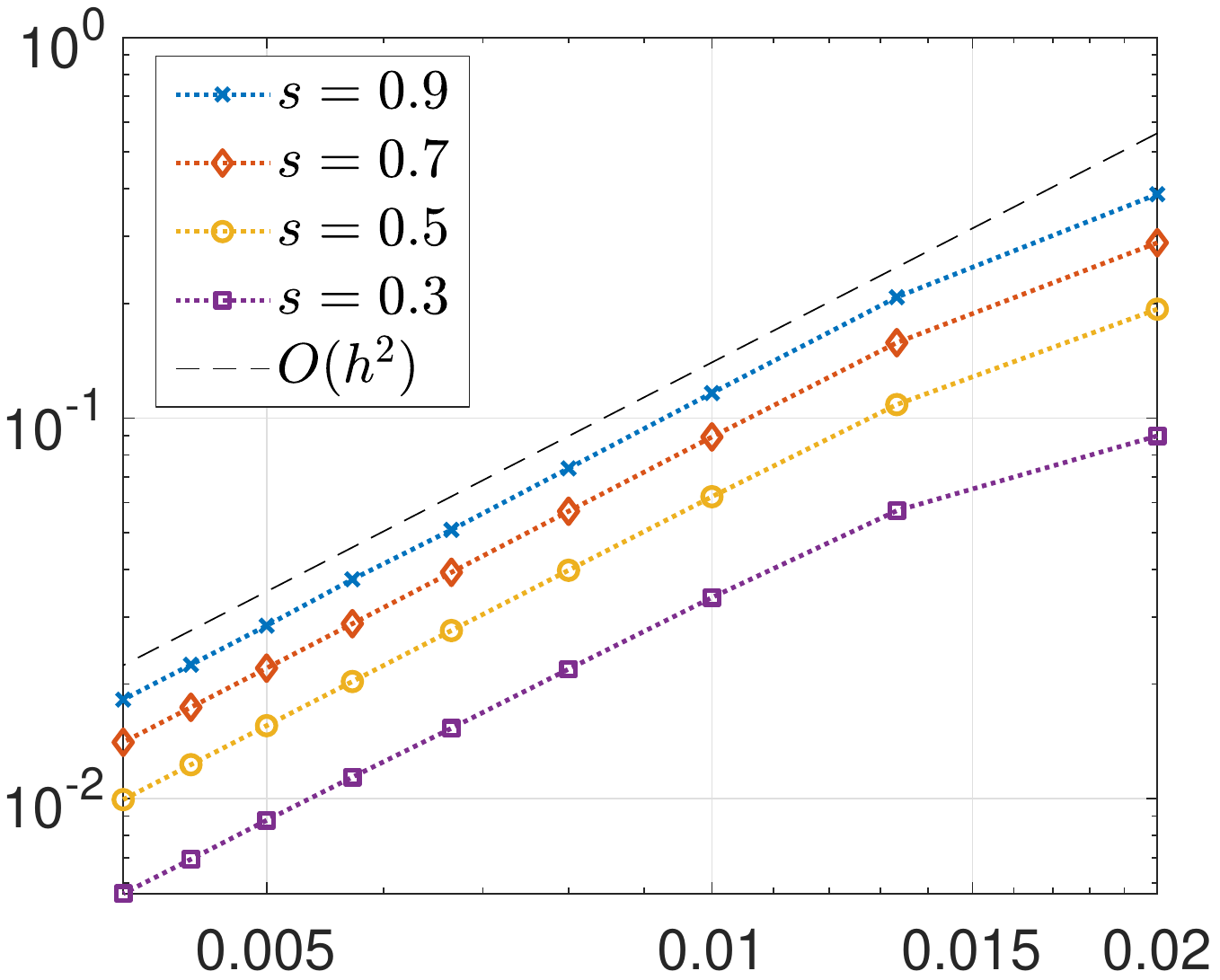}
     \put(100,5){$h$}
     \put(-10,30){\rotatebox{90}{$\mathrm{err_{n_e}(\lambda_j,\lambda^{(h)}_j)}$}}
    \end{overpic}
\caption{Convergence of the eigenvalues of the fractional Laplacian, approximated numerically via a discretization of the Balakrishnan formula based on finite elements, towards the eigenvalues of the fractional Laplacian.
On the $y$-axis the maximum relative error defined in equation \eqref{eq:max-rel-er} for decreasing values of $h$ and several fractional orders is reported, considering $n_e=40$. Convergence order 2 (black dashed line) is shown as reference, being also the order of convergence of eigenvalues for the standard Laplacian discretized via finite elements.}
    \label{fig:eig_conv}
\end{figure}

We choose a domain $\Omega=~(0,1)$ of length $L=1$ and we vary the number of discretization points $n_p$ between $50$ and $250$. For each $n_p$, we compute $K_s$ according to the method presented in Section~\ref{sec:numerics} (see Listing \ref{lst:FractionalLaplacian} for the Matlab code) and use the Matlab command \texttt{eig} to obtain the eigenvalues from which we select the first $n_e=40$ (in ascending order according to their absolute value). Figure~\ref{fig:eig_conv} shows the convergence of the maximum relative error~\eqref{eq:max-rel-er} for fractional order $s=0.9$, $s=0.5$, $s=0.3$ and $s=0.1$. Even if the a priori convergence rate of the discretization method described in the previous section has not been studied, we observe that the error decreases in $O(h^2)$ (denoted by the dashed line in Figure \ref{fig:eig_conv}). This is also the expected convergence rate of eigenvalues for the standard Laplacian discretized via finite elements method with P1 Lagrange elements.

\FloatBarrier
\subsubsection{A simple Poisson problem}
Another simple test can be done which exemplifies the effect of different fractional orders. We consider the one-dimensional Poisson problem on the interval $\Omega =~(0,1)$
\begin{equation}
     (-\Delta)^s u(x) = f(x),
    \label{eq:poisson_problem}
\end{equation}
where $f(x) = 6x+2$, the fractional operator $(-\Delta)^s$ is understood in the spectral definition~\eqref{eq:spectral_definition} and with homogeneous Dirichlet boundary conditions $u\rvert_{\partial \Omega} = 0$. We want to study the convergence of the solution computed applying the method presented in Section~\ref{sec:numerics} as the mesh size decreases. While the solution to problem~\eqref{eq:poisson_problem} for the standard Laplacian ($s=1$) reads $u_{D}(x)=-x(x-1)(x+2)$, the analytical solution for the fractional Laplacian problem is not known. Therefore, we use the numerical solution on a mesh with $500$ nodes, $\bar{u}$, as reference solution.

\begin{figure}
    \centering
    \begin{overpic}[width=0.45\textwidth]{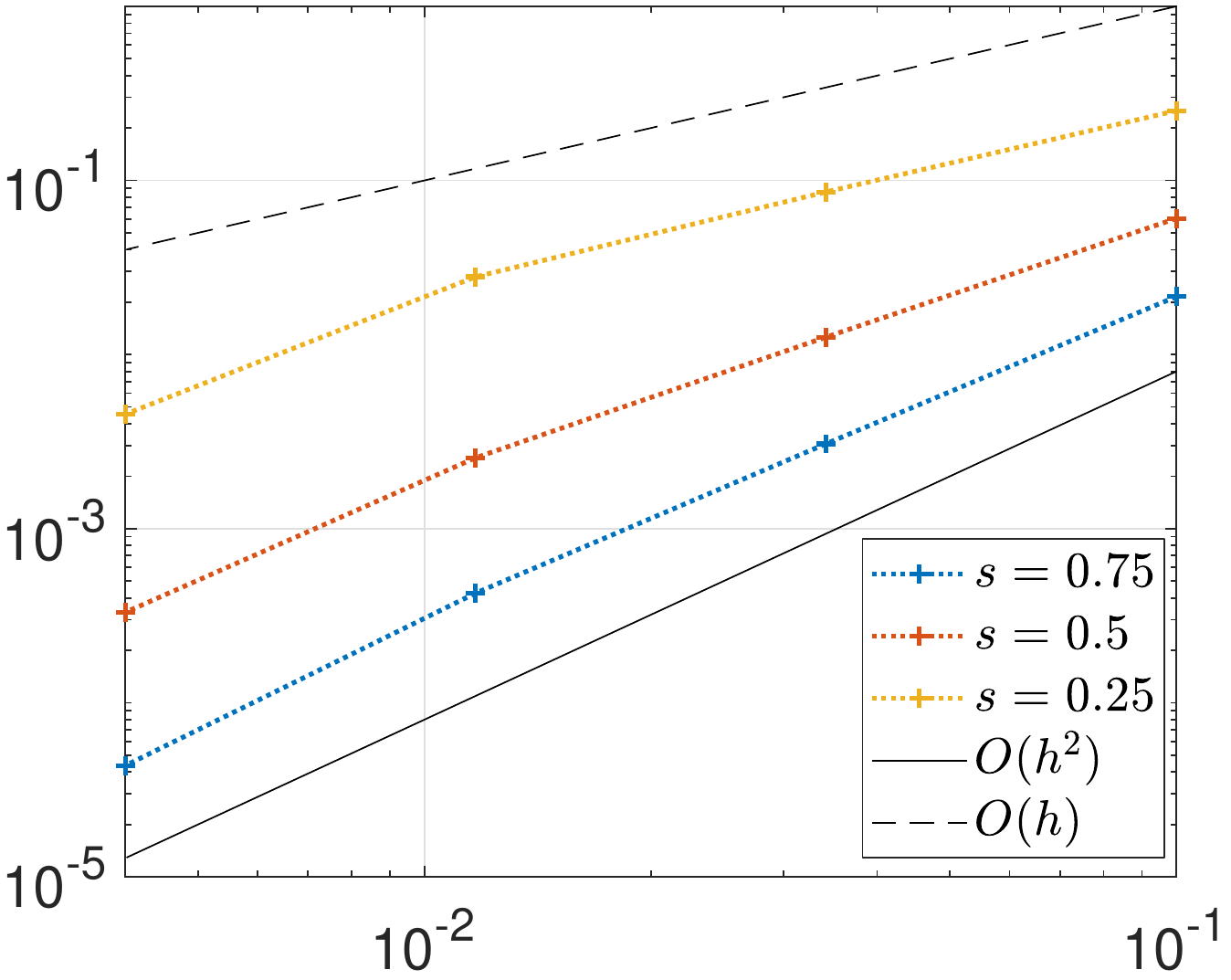}
     \put(100,5){$h$}
    \end{overpic}
    \hspace{1cm}
    \begin{overpic}[width=0.425\textwidth]{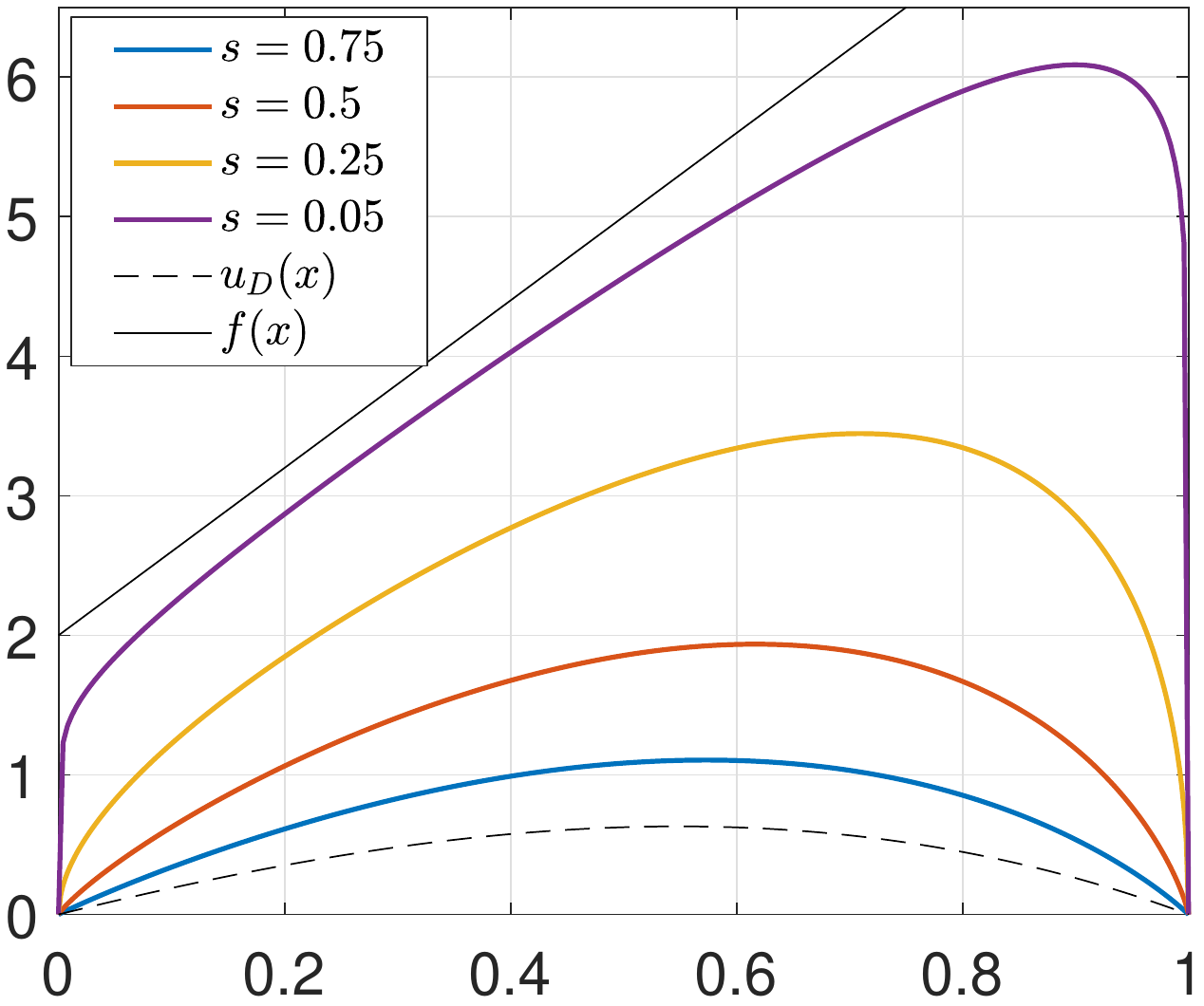}
     \put(100,5){$x$}
     \put(-5,40){\rotatebox{90}{$u(x)$}}
    \end{overpic}
\caption{Left: convergence rates of the numerical solution to~\eqref{eq:poisson_problem} with~$f(x) = 6x+2$, homogeneous Dirichlet boundary conditions on the domain $\Omega=(0,1)$, for different fractional orders $s$. Black dashed and solid lines, corresponding to order of convergence one and two respectively, are shown as reference. Right: numerical solutions for different fractional orders $s$. The dashed curve corresponds to the analytical solutions to~\eqref{eq:poisson_problem} with $s=1$ (standard Laplacian), while the black solid line indicates the reaction term.  (For interpretation of the references to color in this figure legend, the reader is referred to the web version of this article)}
\label{fig:Dirichlet-Poisson}
\end{figure}

We compute solutions for fractional orders ranging between $s=0.1$ and $s=0.9$, increasing the number of nodes from $n_p = 10$ to $n_p = 250$. The convergence towards the numerical solution as the mesh size decreases is shown in Figure~\ref{fig:Dirichlet-Poisson} (left panel) for $s = 0.75$, $s=0.5$, and $s=0.25$. We observe that the convergence rate decreases from $O(h^2)$ for $s = 0.75$ to $O(h)$ for $s=0.25$. Note that a similar effect has been observed in \cite{dohr2019fem}.

The deformation of the numerical solutions as $s$ decreases can be seen in Figure~\ref{fig:Dirichlet-Poisson} (right panel). As $s\to 0$, the solution tends to $6x+2$ which does not satisfy Dirichlet boundary conditions. Therefore this creates an abrupt change in the solution at the boundary values which could affect the convergence.

\FloatBarrier


\section{The (fractional) Allen--Cahn equation}\label{sec:ACe-all}
In the following, we first briefly recall the standard Allen--Cahn equation, its steady state bifurcation diagram and typical solutions before investigating the fractional version of the problem. The figures for the standard case, presented in Subsection \ref{sec:AC-theory} are obtained using the continuation software \texttt{pde2path}, and were already shown in \cite{uecker_pde2path_2014,rademacher2018oopde}. The results presented in Subsection \ref{chapter:AC-results} exploit the new capabilities of \texttt{pde2path} adapted to treat fractional problems, see section \ref{sec:numerics}. \\

\subsection{The Allen--Cahn equation with standard Laplacian}\label{sec:AC-theory}
The Allen--Cahn equation was proposed to model phase separation in alloy systems with multiple components. With standard diffusion and cubic-quintic non-linearity it reads
\begin{equation}\label{eq:AC35}
\partial_t u = \Delta u + \mu u + u^3 - \gamma u^5,
\end{equation}
where $u=u(x,t)\in \mathbb{R}$ is the unknown and $\mu, \gamma \in \mathbb{R}$ are parameters.

We consider the equation on a one-dimensional bounded domain $\Omega =(-L/2,L/2)$, where $L = 10$ with homogeneous Dirichlet boundary conditions, $u\rvert_{\partial \Omega} = 0$. Then $ \bar{u} = 0$ is a homogeneous stationary solution of equation~\eqref{eq:AC35} for all choices of $\mu$, $\gamma$. Further, linearizing equation~\eqref{eq:AC35} around $\bar{u}$, it can be shown that this state is stable as long as $\mu<({\pi}/{L})^2$, and first becomes unstable to modes with wavenumber $k_c={\pi}/{L}$ when the parameter reaches the critical value $\mu_c = ({\pi}/{L})^2$. Subsequent modes become unstable, at $\mu_j = \left( {j\pi}/{L} \right)^2,\; j\geq 2$.

\begin{figure}
\begin{multicols}{2}
\centering
\vspace*{\fill}
\hspace{0.3cm}
\begin{subfigure}[b]{0.5\textwidth}
    \centering
    \begin{overpic}[width=\textwidth]{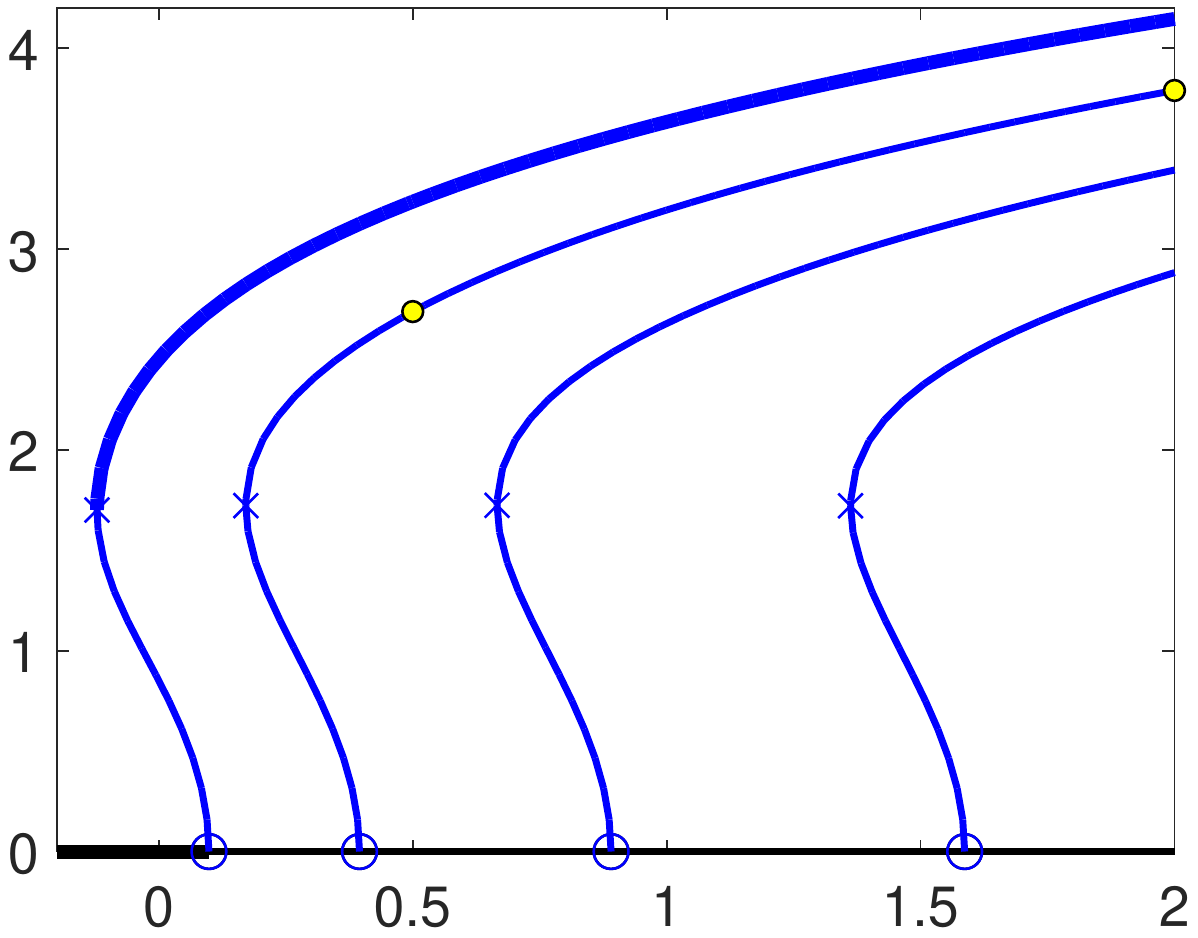}
\put(31,10){$\mathtt{B_2}$}
\put(31,45){$(a)$}
\put(92,65){$(b)$}
\put(87,0){$\mu$}
\put(-10,40){\rotatebox{90}{$\|u\|_{L^2}$}}
    \end{overpic}
\end{subfigure}\\
\vspace*{\fill}
\newpage
\vspace*{\fill}
\begin{subfigure}[b]{0.2\textwidth}
     \begin{overpic}[width=\textwidth]{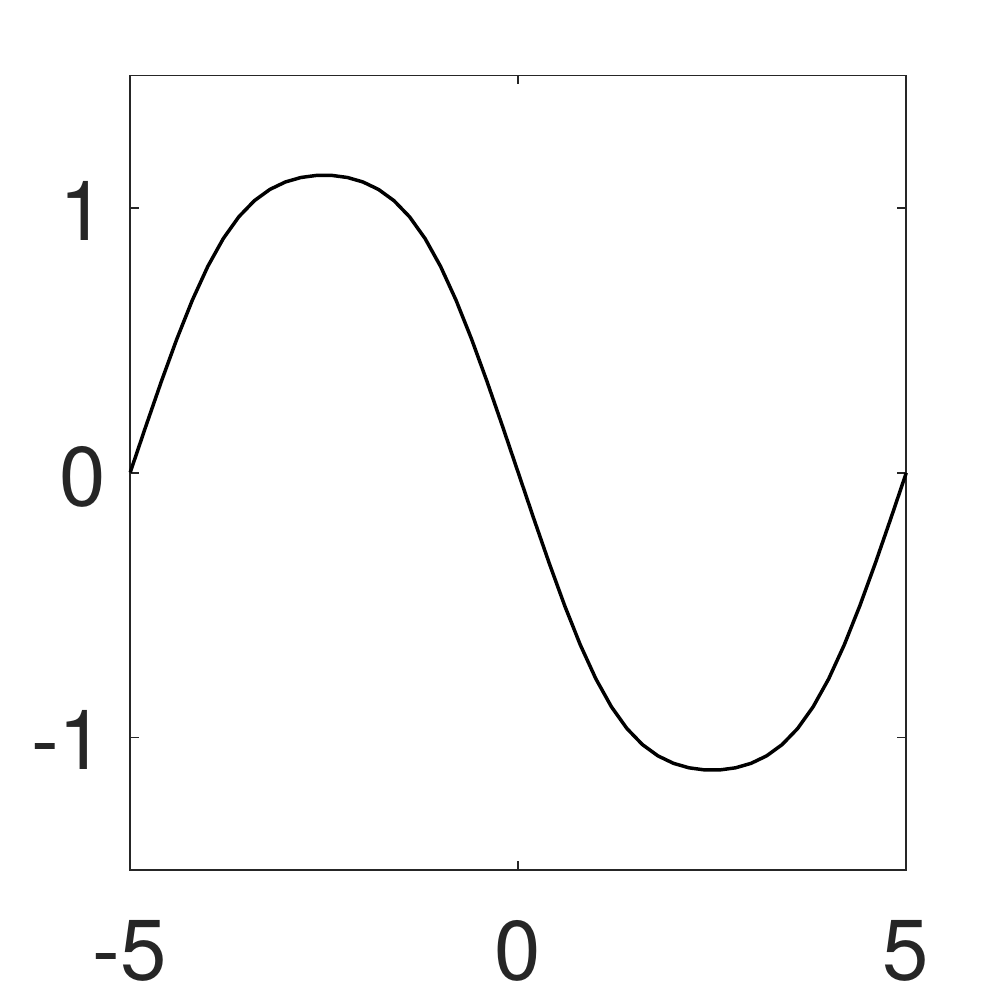}
\put(75,-2){$x$}
\put(-9,45){\rotatebox{90}{$u(x)$}}
\end{overpic}
   \caption{$\mu=0.5$}
\label{subfig:AC-b2-sol1}
\end{subfigure}
\hspace{0.5cm}
\begin{subfigure}[b]{0.2\textwidth}
     \begin{overpic}[width=\textwidth]{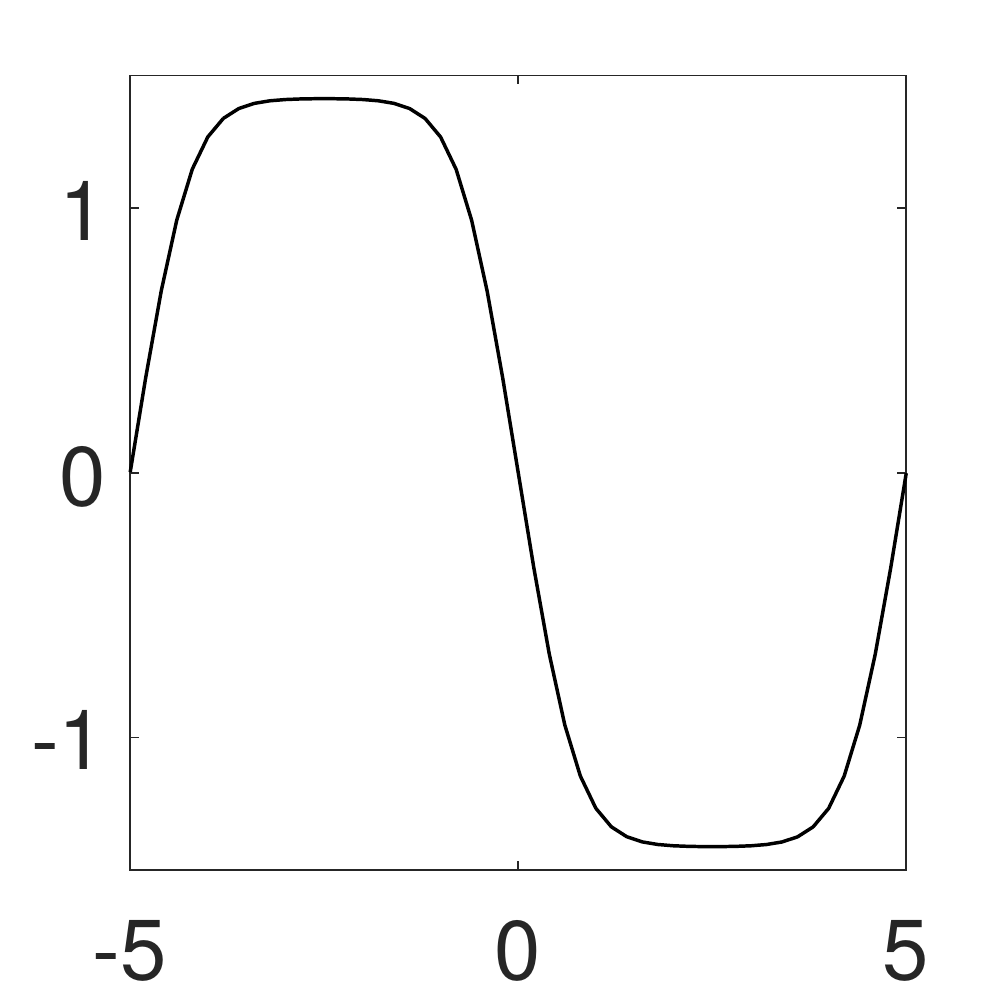}
\put(75,-2){$x$}
\put(-9,45){\rotatebox{90}{$u(x)$}}
\end{overpic}
    \caption{$\mu=2$}
    \label{subfig:AC-b2-sol4}
\end{subfigure}
\vspace*{\fill}
\end{multicols}
\vspace{-0.5cm}
\caption{Left panel: steady state bifurcation diagram for the Allen--Cahn equation~\eqref{eq:AC35} with $\gamma~=~1$, on the domain $\Omega=(-5, 5)$ with homogeneous Dirichlet boundary conditions. Right panel (\subref{subfig:AC-b2-sol1}--\subref{subfig:AC-b2-sol4}): steady state solution profiles of the Allen--Cahn equation, along the second bifurcating branch originating at the bifurcation point~$\mathtt{B_2}$. In the bifurcation diagram, thicker lines denote stable solution while thinner lines correspond to unstable ones. Circles and crosses indicate branch and fold points respectively.}
\label{fig:Dirichlet-AC}
\end{figure}
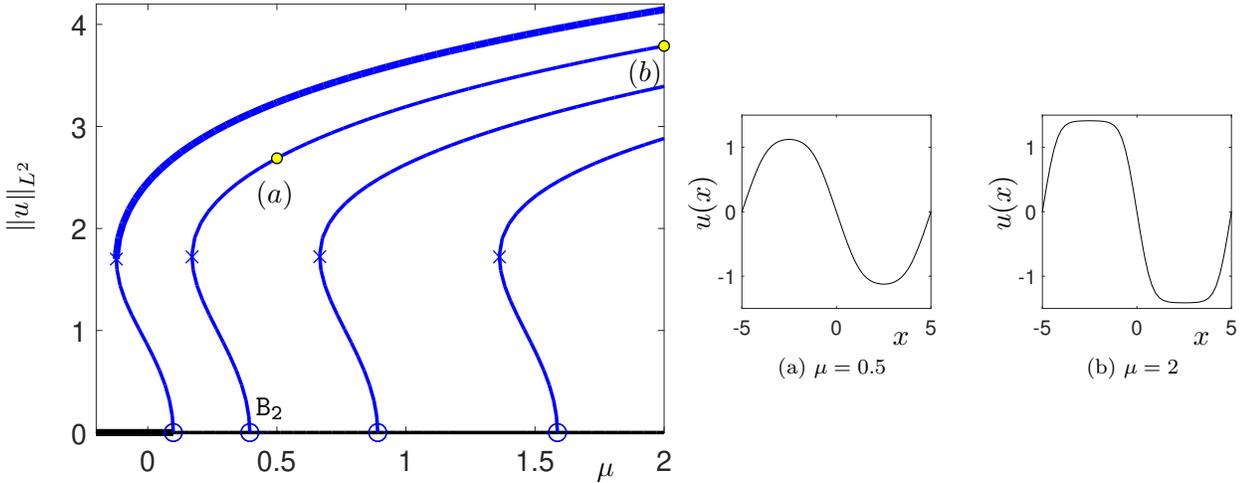

The steady state bifurcation structure of equation~\eqref{eq:AC35} can be computed numerically to obtain a detailed picture far from the homogeneous branch. The branches of non-homogeneous stationary solutions bifurcating from the homogeneous branch are illustrated in Figure~\ref{fig:Dirichlet-AC} for the case $\gamma=1$. They bifurcate via pitchfork bifurcations and experience a fold at larger amplitude. The first bifurcating branch is unstable until the fold and becomes stable after the fold. All other branches are unstable for $\mu<2$, and they stabilize for greater values of $\mu$. In addition, Figures \ref{subfig:AC-b2-sol1}--\ref{subfig:AC-b2-sol4} show solution profiles along the second branch (originating at the second bifurcation point $\mu_2 = (\pi/5)^2$, denoted $\mathtt{B_2}$ in the figure) for increasing vales of $\mu$: spatial subregions appear occupied by one of two phases between which there is an interface which sharpens as $\mu$ increases.

\FloatBarrier

\subsection{The fractional Allen--Cahn equation}\label{chapter:AC-results}
We start with the fractional Allen--Cahn equation with cubic--quintic nonlinearity~\eqref{eq:AC-frac} on a 1D domain of length $L$ with homogeneous Dirichlet boundary conditions. As for the standard Allen--Cahn equation~\eqref{eq:AC35}, the homogeneous solution $\Bar{u}=0$ is a stationary steady state of equation~\eqref{eq:AC-frac}. Furthermore, using the Fourier property of the fractional Laplacian, one can easily find the bifurcations from the homogeneous states to be located at\\[-0.2cm]
\begin{equation}\label{eq:AC-bifurcation-points-s}
\mu_j = \left( \frac{j\pi}{L} \right )^{2s}, \quad j \geq 1,
\end{equation}
from which we expect the bifurcation points to move towards $\mu=~1$ as $s$ tends to zero. To obtain a complete picture of the bifurcation diagram far from the homogeneous solution, numerical continuation techniques are required.

\begin{figure}
  \centering
  \begin{overpic}[width=0.45\textwidth]{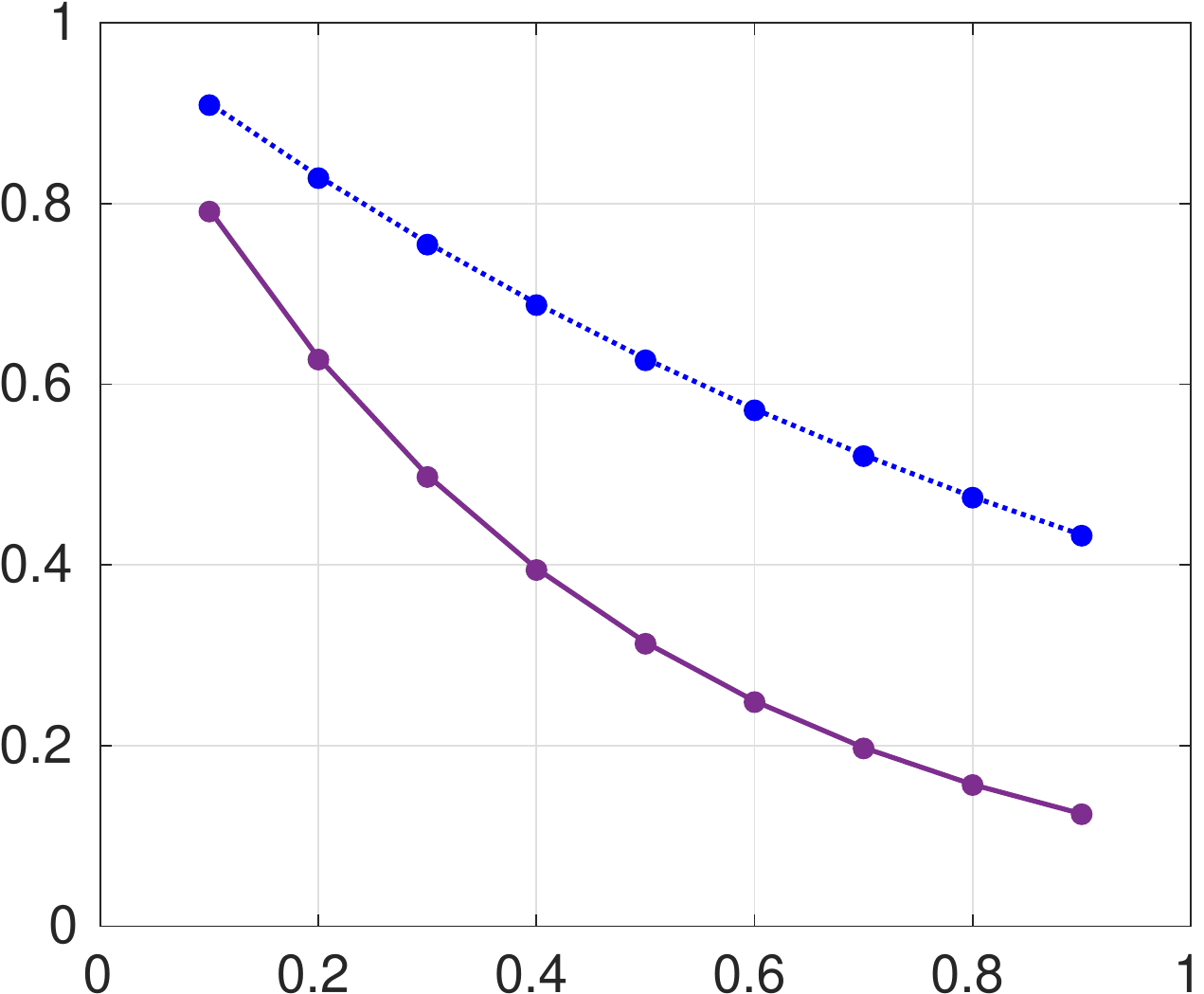}
  \put(100,5){$s$}
  \put(90,41){$\mu_2$}
  \put(90,18){$\mu_1$}
  \end{overpic}
\caption{Evolution of the location of the first (dashed blue line) and second (continuous purple line) bifurcation points on the homogeneous branch as a function of the fractional order of the Laplacian $s$ for the fractional Allen--Cahn equation~\eqref{eq:AC-frac} on a bounded domain with homogeneous Dirichlet boundary conditions. The dots indicate numerical values, while the lines correspond to values computed analytically via formula~\eqref{eq:AC-bifurcation-points-s}.  (For interpretation of the references to color in this figure legend, the reader is referred to the web version of this article)}
\label{fig:AC-bifurcation-pts-1D}
\end{figure}

The following numerical results are obtained fixing $\gamma=1$ and a domain $\Omega=(-5, 5)$ of length $L=10$. We used $n_p=301$ mesh points, with corresponding meshsize $h=0.0333$. The Dirichlet boundary conditions read $u(-5)=u(5)=0$. Figure~\ref{fig:AC-bifurcation-pts-1D} shows that, indeed, the shifting of the first and second bifurcation points as $s$ decreases agrees with the values computed using formula~\eqref{eq:AC-bifurcation-points-s}. The maximum relative error for the first bifurcation point is $4\cdot 10^{-3}$ and for the second bifurcation point it is $2\cdot 10^{-3}$.

\begin{figure}
\centering
\begin{subfigure}{.3\textwidth}
  \centering
  \begin{overpic}[width=\textwidth]{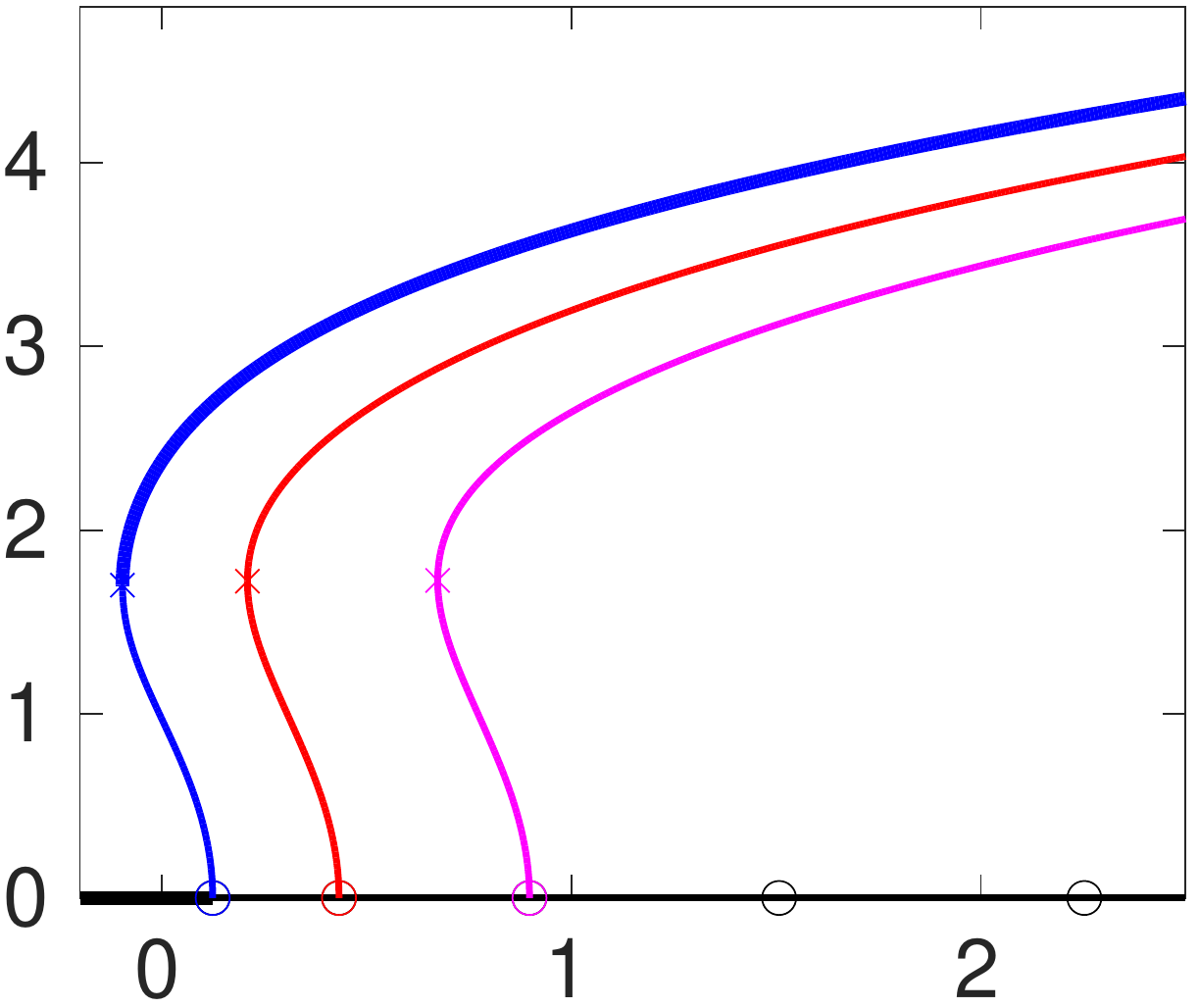}
    \put(-9,45){\rotatebox{90}{$\|u\|_{L^2}$}}
     \put(85,-5){$\mu$}
    \end{overpic}
  \caption{$s=0.9$}
  \label{fig:AC-bif-s09}
\end{subfigure}
\hspace{0.1cm}
\begin{subfigure}{.3\textwidth}
  \centering
  \begin{overpic}[width=\textwidth]{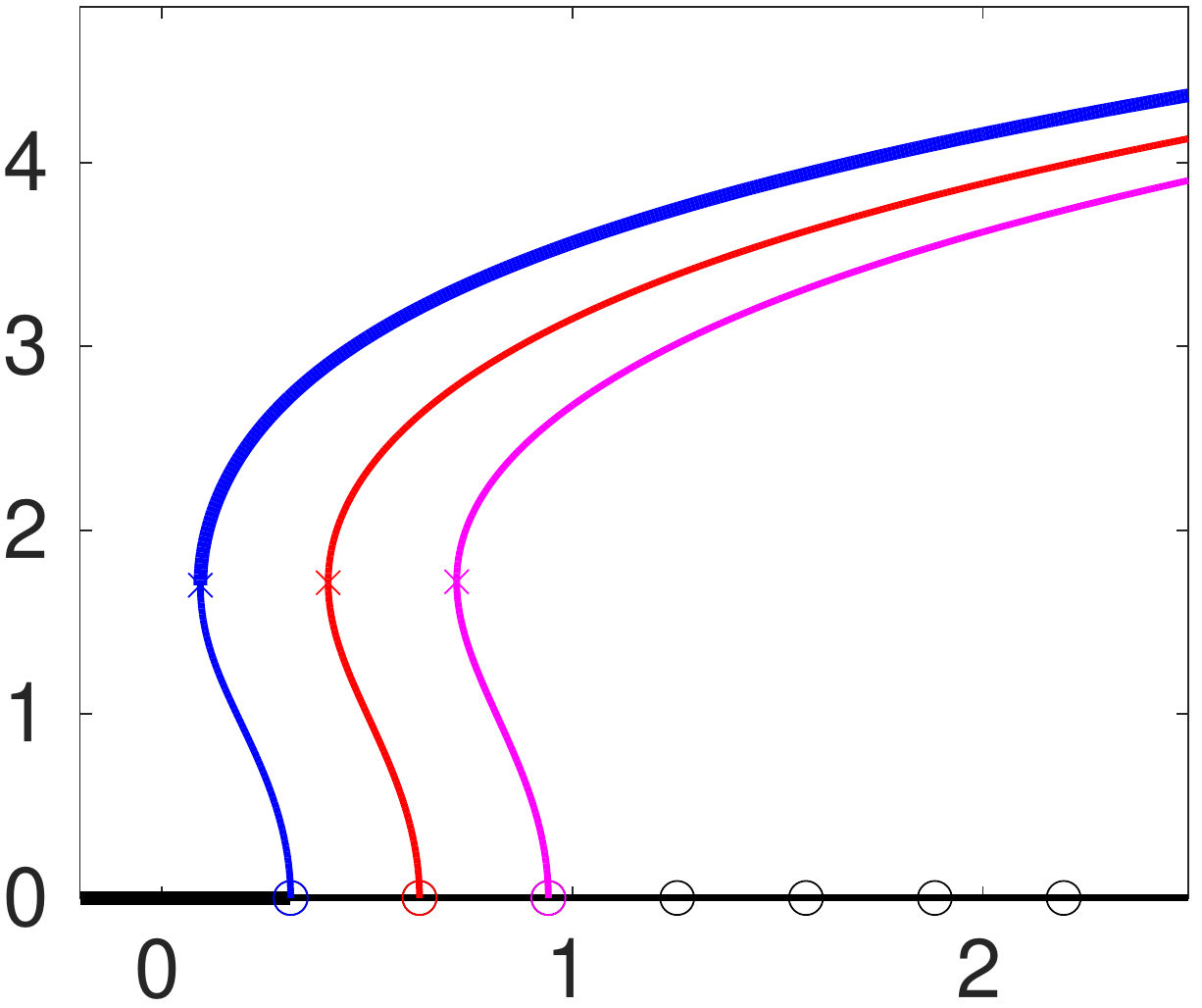}
     \put(85,-5){$\mu$}
    \end{overpic}
  \caption{$s=0.5$}
  \label{fig:AC-bif-s05}
\end{subfigure}
\hspace{0.1cm}
\begin{subfigure}{.3\textwidth}
  \centering
  \begin{overpic}[width=\textwidth]{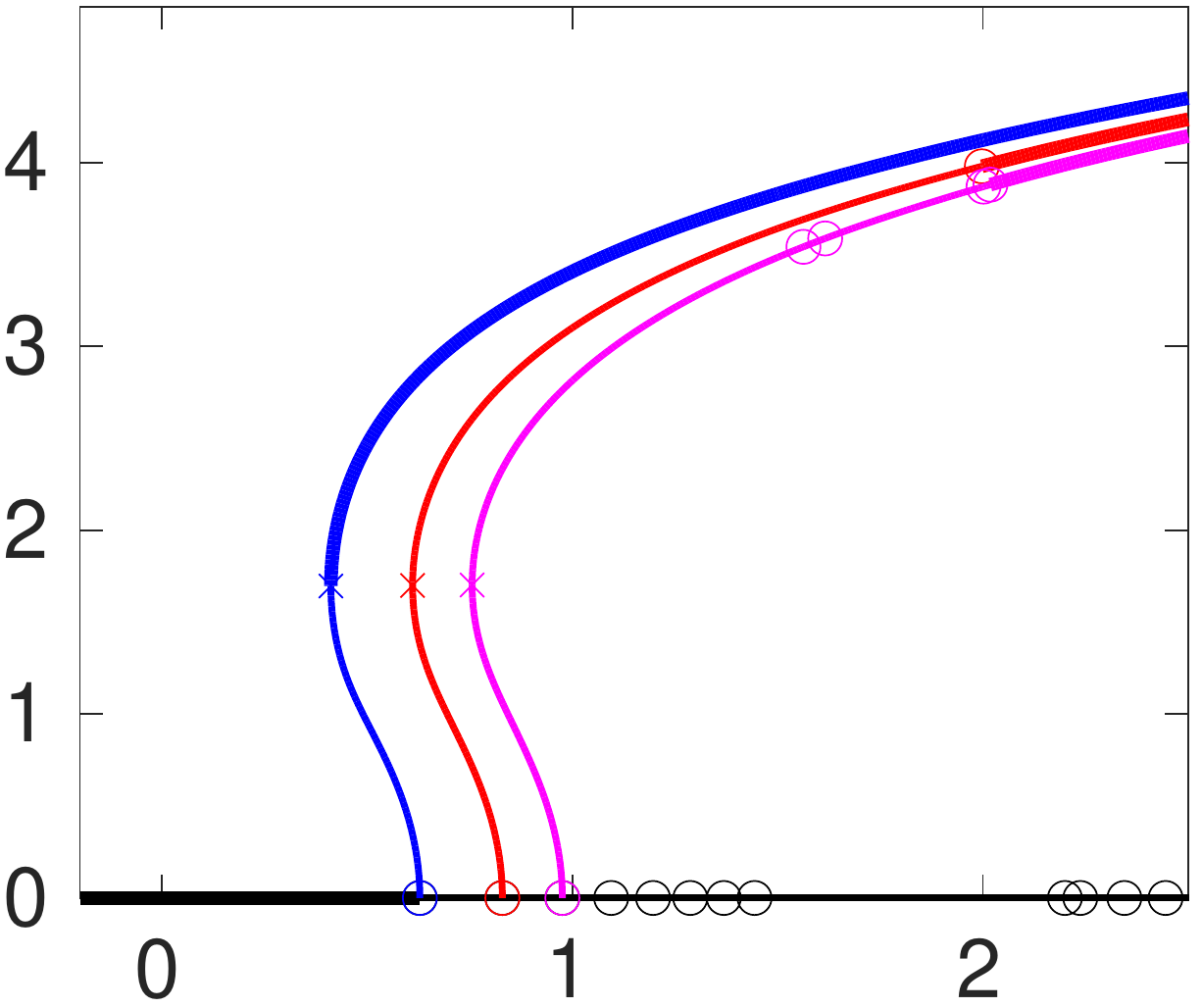}
     \put(85,-5){$\mu$}
    \end{overpic}
  \caption{$s=0.2$}
  \label{fig:AC-bif-s02}
\end{subfigure}
\caption{Bifurcation diagram of the fractional Allen--Cahn equation~\eqref{eq:AC-frac} on a bounded domain with homogeneous Dirichlet boundary conditions for decreasing values of the fractional order $s$. For sake of clarity, only the first three branches bifurcating from the homogeneous state are shown. Thick and thin lines denote stable and unstable solutions, while circles and crosses indicate branch and fold points respectively. (For interpretation of the references to color in this figure legend, the reader is referred to the web version of this article)}
\label{fig:AC-bifurcations-1D}
\vspace*{0.2cm}
\begin{multicols}{3}
\centering
\begin{subfigure}{.3\textwidth}
\centering
     \begin{overpic}[width=0.7\textwidth,trim=18 20 0 0,clip]{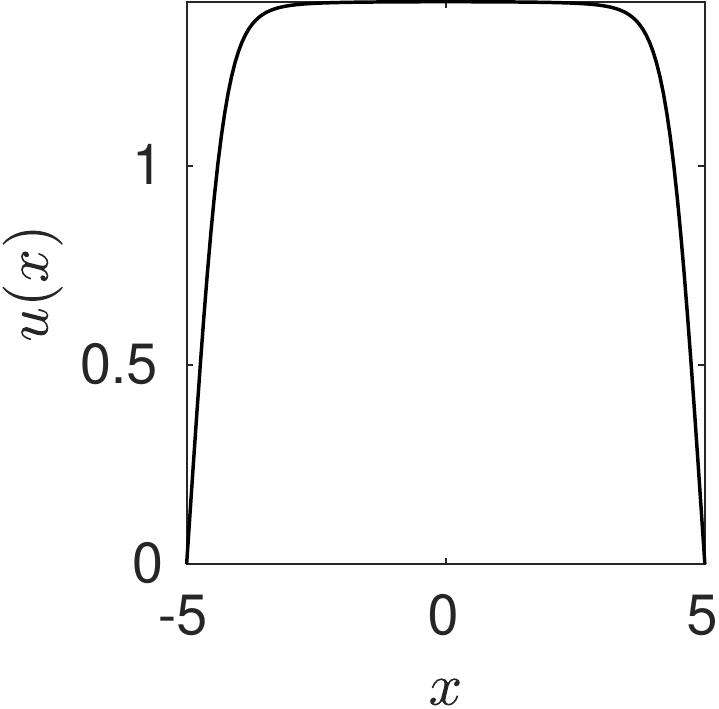}
     \put(-10,45){\rotatebox{90}{$u(x)$}}
     \put(100,10){$x$}
    \end{overpic}
  \caption{$s=0.9$, 1st branch}
  \label{fig:AC-bif-s09-sol1}
\end{subfigure}\\[0.3cm]
\begin{subfigure}{.3\textwidth}\centering
     \begin{overpic}[width=0.7\textwidth,trim=20 20 0 0,clip]{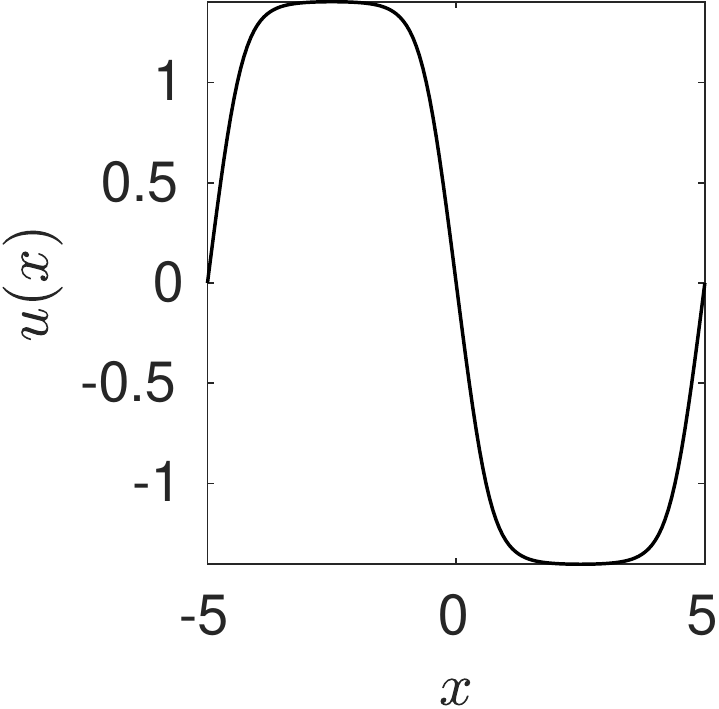}
     \put(-10,45){\rotatebox{90}{$u(x)$}}
     \put(100,10){$x$}
    \end{overpic}
  \caption{$s=0.9$, 2nd branch}
  \label{fig:AC-bif-s09-sol2}
\end{subfigure}\\[0.3cm]
\begin{subfigure}{.3\textwidth}\centering
     \begin{overpic}[width=0.7\textwidth,trim=18 20 0 0,clip]{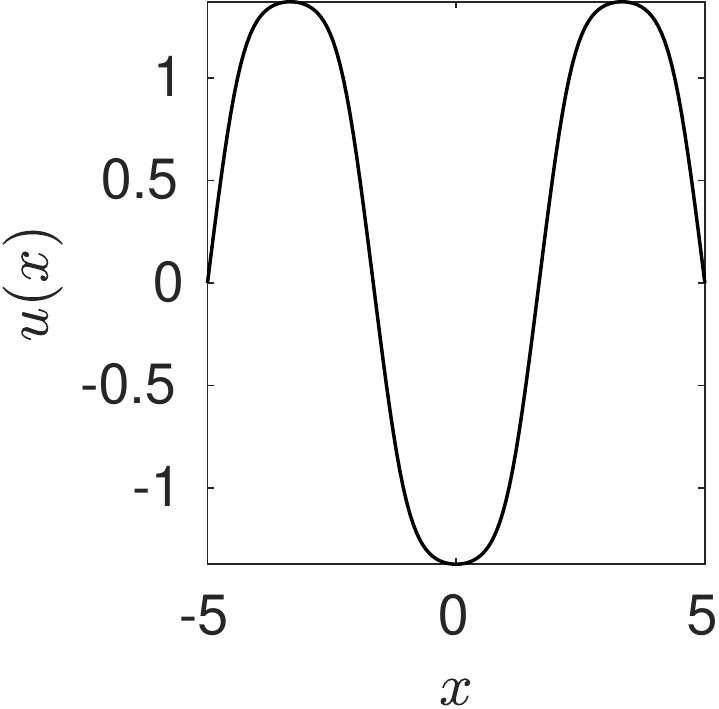}
     \put(-10,45){\rotatebox{90}{$u(x)$}}
     \put(100,10){$x$}
    \end{overpic}
  \caption{$s=0.9$, 3rd branch}
  \label{fig:AC-bif-s05-sol3}
\end{subfigure}
\newpage
\begin{subfigure}{.3\textwidth}\centering
 \centering
  \begin{overpic}[width=0.7\textwidth,trim=18 20 0 0,clip]{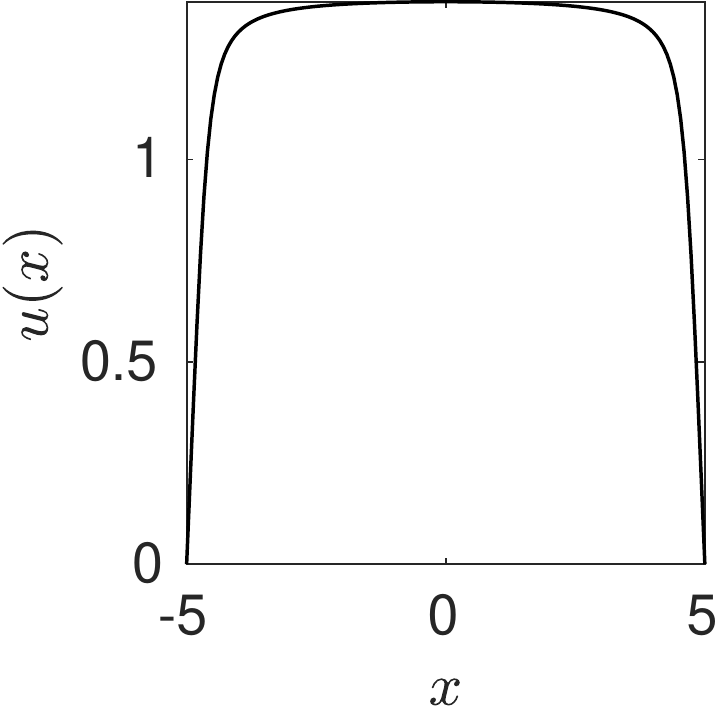}
     \put(100,10){$x$}
    \end{overpic}
  \caption{$s=0.5$, 1st branch}
  \label{fig:AC-bif-s05-1}
\end{subfigure}\\[0.3cm]
\begin{subfigure}{.3\textwidth}\centering
  \centering
  \begin{overpic}[width=0.7\textwidth,trim=18 20 0 0,clip]{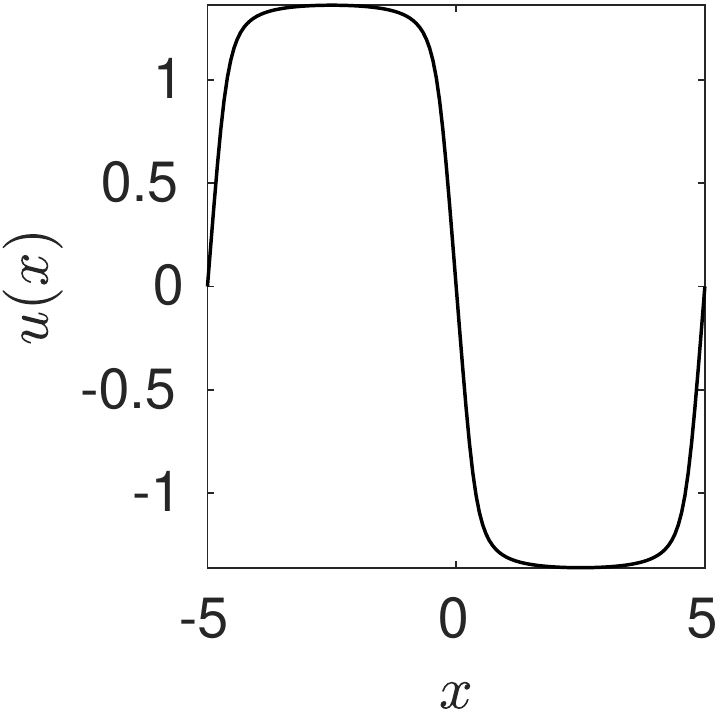}
     \put(100,10){$x$}
    \end{overpic}
 \caption{$s=0.5$, 2nd branch}
  \label{fig:AC-bif-s05-2}
\end{subfigure}\\[0.3cm]
\begin{subfigure}{.3\textwidth}\centering
  \centering
  \begin{overpic}[width=0.7\textwidth,trim=18 20 0 0,clip]{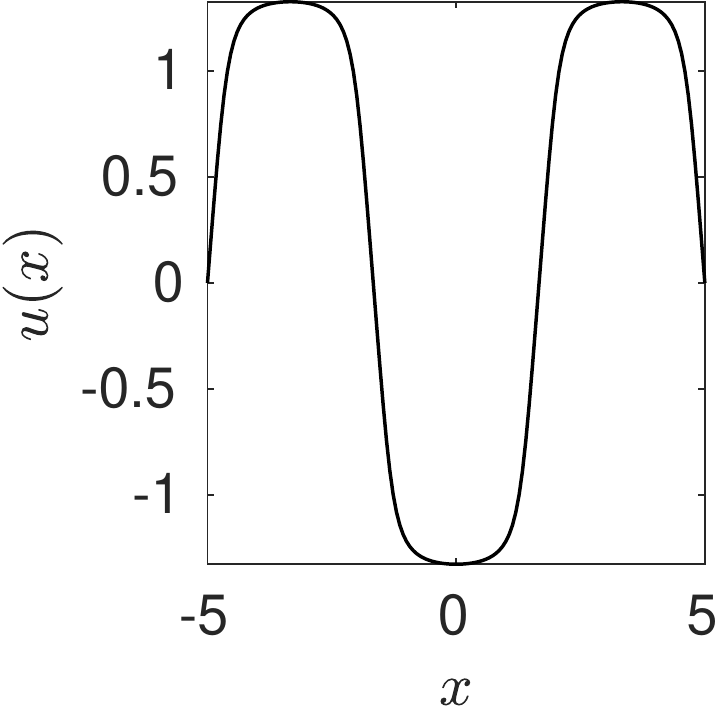}
     \put(100,10){$x$}
    \end{overpic}
  \caption{$s=0.5$, 3rd branch}
  \label{fig:AC-bif-s05-3}
\end{subfigure}
\newpage
\begin{subfigure}{.3\textwidth}\centering
  \centering
     \begin{overpic}[width=0.7\textwidth,trim=18 20 0 0,clip]{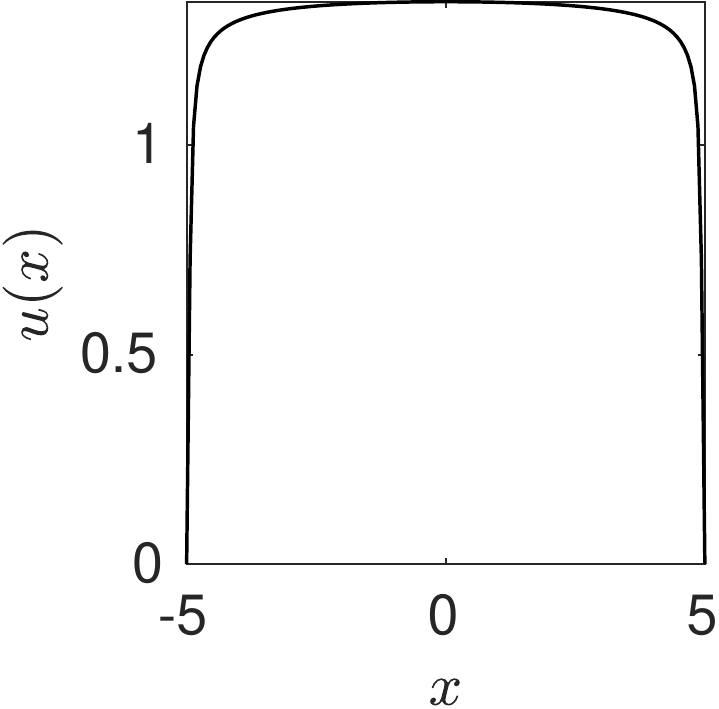}
     \put(100,10){$x$}
    \end{overpic}
  \caption{$s=0.2$, 1st branch}
  \label{fig:AC-bif-s02-1}
\end{subfigure}\\[0.3cm]
\begin{subfigure}{.3\textwidth}\centering
  \centering
     \begin{overpic}[width=0.7\textwidth,trim=18 20 0 0,clip]{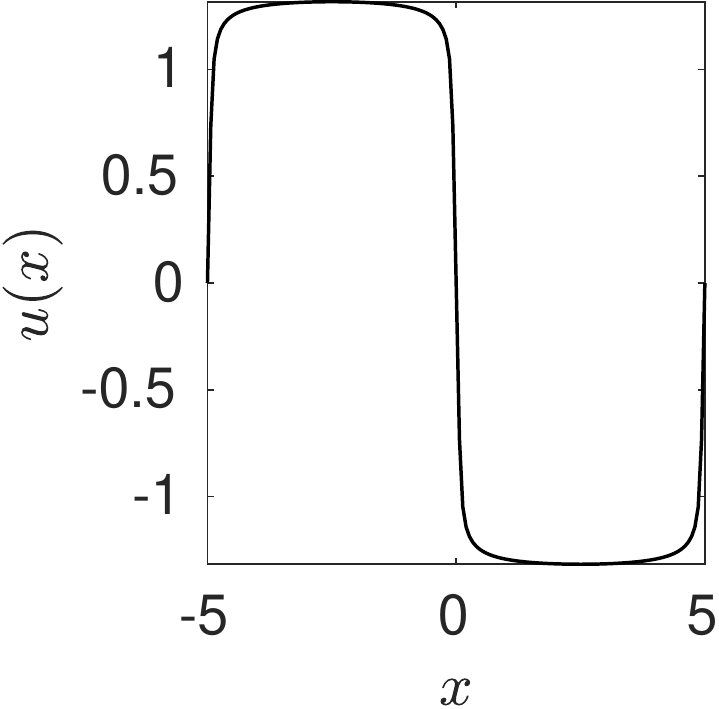}
     \put(100,10){$x$}
    \end{overpic}
  \caption{$s=0.2$, 2nd branch}
  \label{fig:AC-bif-s02-2}
\end{subfigure}\\[0.3cm]
\begin{subfigure}{.3\textwidth}\centering
  \centering
  \begin{overpic}[width=0.7\textwidth,trim=18 20 0 0,clip]{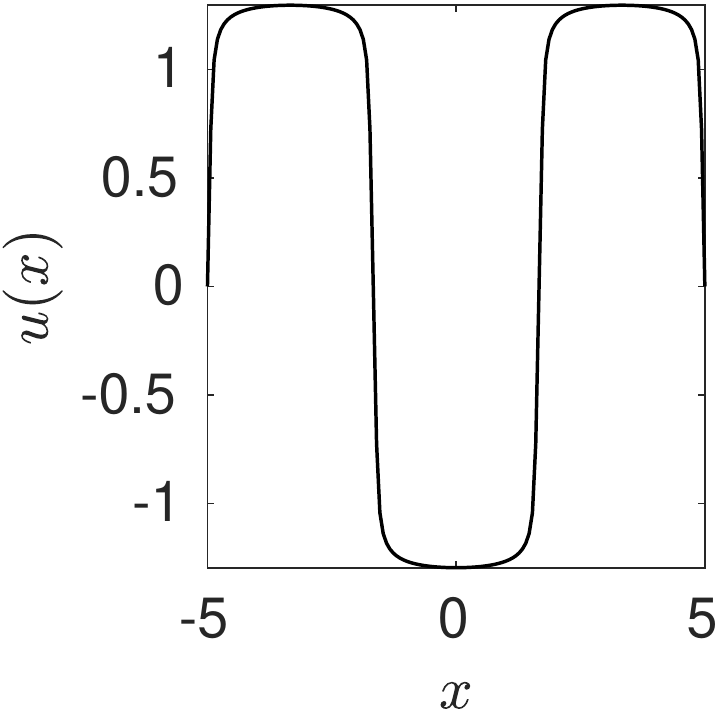}
     \put(100,10){$x$}
    \end{overpic}
  \caption{$s=0.2$, 3rd branch}
  \label{fig:AC-bif-s02-3}
\end{subfigure}
\end{multicols}
\vspace{-0.5cm}
\caption{Stationary solutions to the fractional Allen--Cahn equation~\eqref{eq:AC-frac} with homogeneous Dirichlet boundary conditions, on the first, second and third bifurcating branches of Figure~\ref{fig:AC-bifurcations-1D} at $\mu=2$, for fractional orders $s=0.9$, $s=0.5$ and $s=0.2$. Note that for $s=0.2$ all the three solutions are stable. One also observes clearly the sharpening of the layers (or ``teeth'') once the fractional order is decreased.}
\label{fig:AC-solutions-mu2-1D}
\end{figure}

Figure~\ref{fig:AC-bifurcations-1D} shows the evolution of the bifurcation diagram (with respect to the $L^2$-norm of $u$) of equation~\eqref{eq:AC-frac} as the fractional order of the Laplacian decreases. We observe the gathering of branches which is a consequence of the shifting of bifurcation points. In addition, one can see that a stable region appears on the second and third bifurcating branches at $s=0.2$ (Figure~\ref{fig:AC-bif-s02}). These stable states exist in fact for all fractional orders, however they appear at much smaller $\mu$ values when $s$ is decreased. This is better understood by looking at the solution profiles. We have seen in Section \ref{sec:AC-theory} that the steady state solutions of the Allen--Cahn equation correspond to situations where the space is divided into subregions occupied by one of two phases. We also know that the interface between these phases sharpens as $\mu$ increases. Similar observations can be made for the steady states of the fractional Allen--Cahn equation~\eqref{eq:AC-frac}. However, the ``sharpening'' of solutions occurs at much lower $\mu$, as $s$ decreases. Thus, assuming a relation between the sharpening of solutions and their stability, it seems reasonable for the stability of solutions to be shifted to lower $\mu$ values too. Figure~\ref{fig:AC-solutions-mu2-1D} illustrates the change in solution profiles as $s$ decreases, for $\mu=2$ fixed, on the first, second and third branch. Note that for each solution profile $u(x)$ in Figure~\ref{fig:AC-solutions-mu2-1D}, $-u(x)$, is also a solution to~\eqref{eq:AC-frac}. These two solutions lie on the same branch in the bifurcation diagram due to symmetry.

\FloatBarrier

\section{The (fractional) Swift--Hohenberg equation}\label{sec:SHe-all}
In this section, we first recall the standard version of the Swift--Hohenberg equation before investigating it's fractional version. As explained in the previous section, the standard steady state bifurcation diagrams and typical solutions, presented in Subsection \ref{sec:SHe-theory} are obtained using the continuation software \texttt{pde2path}, and were already shown in \cite{Uecker2014}. The results presented in Subsection \ref{chapter:SHe-results} exploit the new capabilities of \texttt{pde2path} adapted to treat fractional problems, see section \ref{sec:numerics}. \\

\subsection{The Swift--Hohenberg equation with standard Laplacian}\label{sec:SHe-theory}
The Swift--Hohenberg equation is one of the simplest phenomenological model for pattern formation. With cubic--quintic non-linearity it reads
\begin{equation}\label{eq:SHE}
\partial_t u = - (1 + \Delta)^2u + \mu u + \nu u^3 - u^5,
\end{equation}
with parameter $\mu \in \mathbb{R}$ and $\nu > 0$. This equation presents the following characteristics: it is reversible in space, i.e. it is equivalent under the transformation $x \rightarrow -x$, $u \rightarrow u$, and, due to the choice on nonlinearity, is also symmetric under the transformation $x \rightarrow x$, $u \rightarrow -u$.

We consider the one-dimensional bounded domain $\Omega$ of length $L = m\pi$, $m \in \mathbb{N}_{>0}$ with homogeneous Dirichlet boundary conditions $u\rvert_{\partial \Omega} = 0$. In this case the trivial state $\bar{u}=0$ is the only spatially homogeneous steady state of equation~\eqref{eq:SHE}. Linearizing equation~\eqref{eq:SHE} around $\bar{u}$, it can be shown that the trivial state is stable as long as $\mu < 0 $ and first becomes unstable to modes with wavenumber $k_c=1$ at $\mu_c=0$. On the bounded domain $\Omega$ further bifurcations occur at the parameter values $\mu_j=(1-(j\pi/L)^2)^2,\; j\geq 1$.

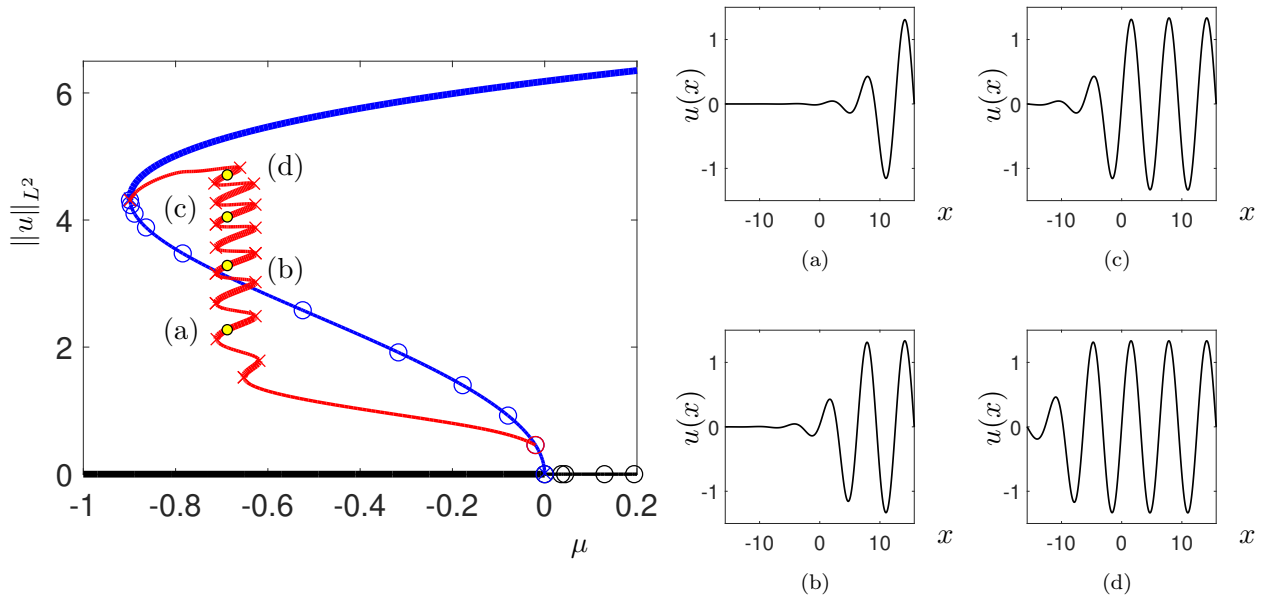
\begin{figure}
\begin{multicols}{3}
\vspace*{\fill}
\hspace{0.5cm}
\begin{subfigure}[b]{0.5\textwidth}
    \centering
    \begin{overpic}[width=\textwidth]{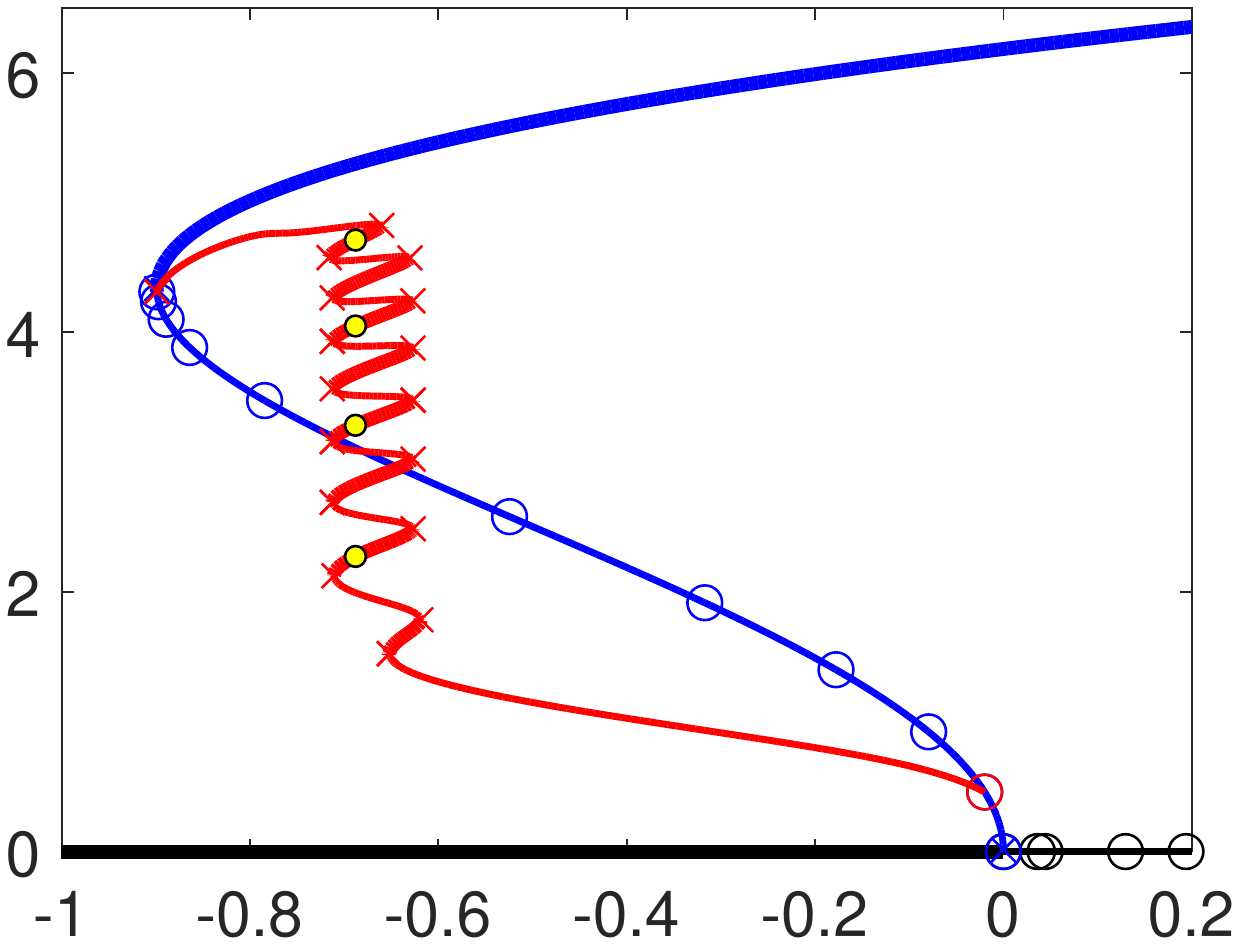}
    \put(85,-5){$\mu$}
    \put(-7,45){\rotatebox{90}{$\|u\|_{L^2}$}}
    \put(18,30){(a)}
    \put(35,40){(b)}
    \put(18,50){(c)}
    \put(35,57){(d)}
    \end{overpic}
\end{subfigure}
\vspace*{\fill}
\newpage
\hspace{3.5cm}
\begin{subfigure}[b]{0.2\textwidth}
     \begin{overpic}[width=\textwidth]{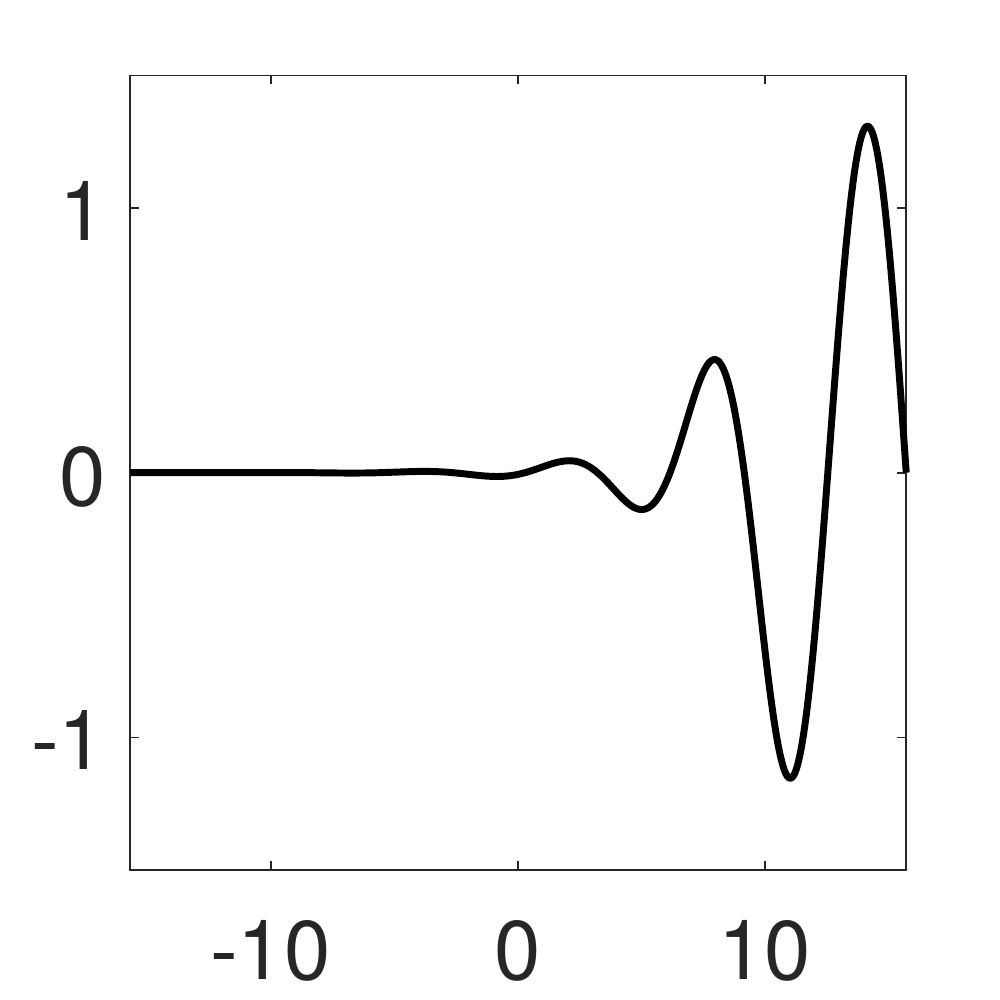}
     \put(-7,45){\rotatebox{90}{$u(x)$}}
     \put(100,5){$x$}
    \end{overpic}
     \caption{}
     \label{subfig-1:DSH}
\end{subfigure}\\

\hspace{3.5cm}
\begin{subfigure}[b]{0.2\textwidth}
     \begin{overpic}[width=\textwidth]{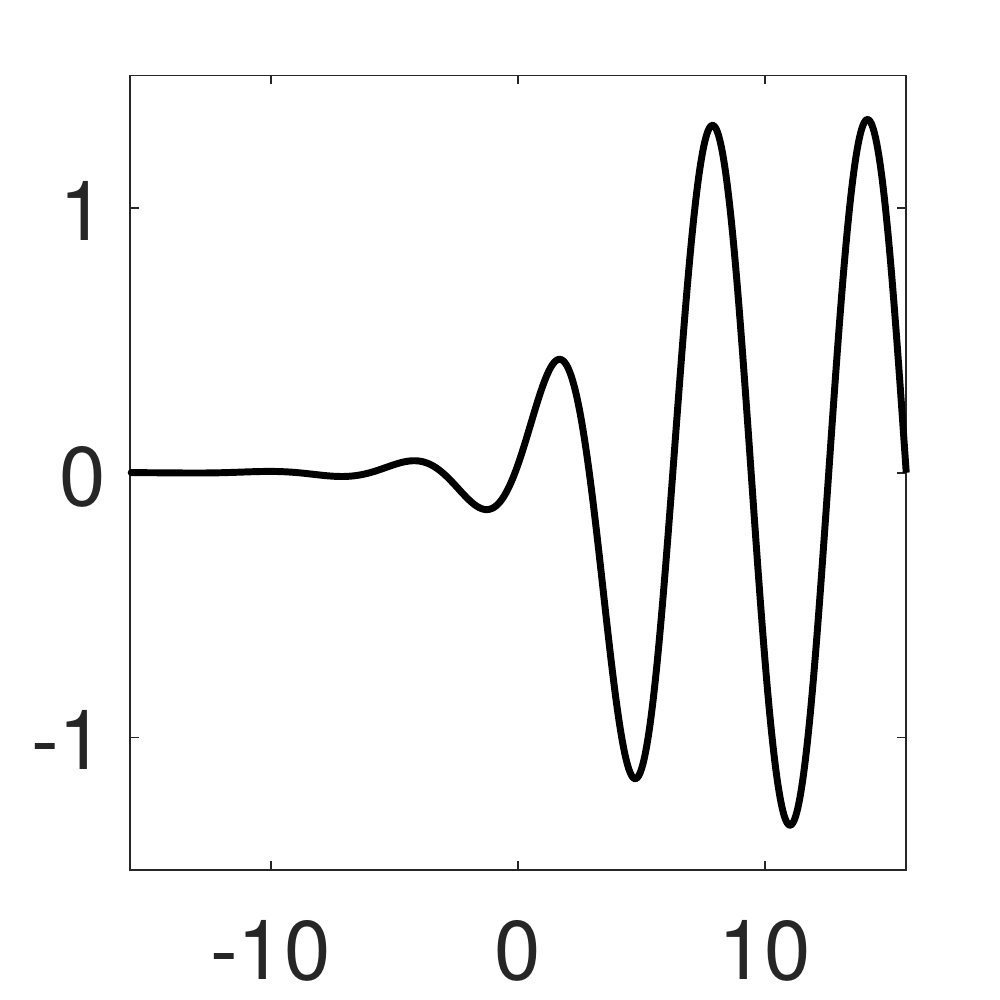}
     \put(-7,45){\rotatebox{90}{$u(x)$}}
     \put(100,5){$x$}
    \end{overpic}
     \caption{}
     \label{subfig-2:DSH}
\end{subfigure}
\newpage
\hspace{2cm}
\begin{subfigure}[b]{0.2\textwidth}
     \begin{overpic}[width=\textwidth]{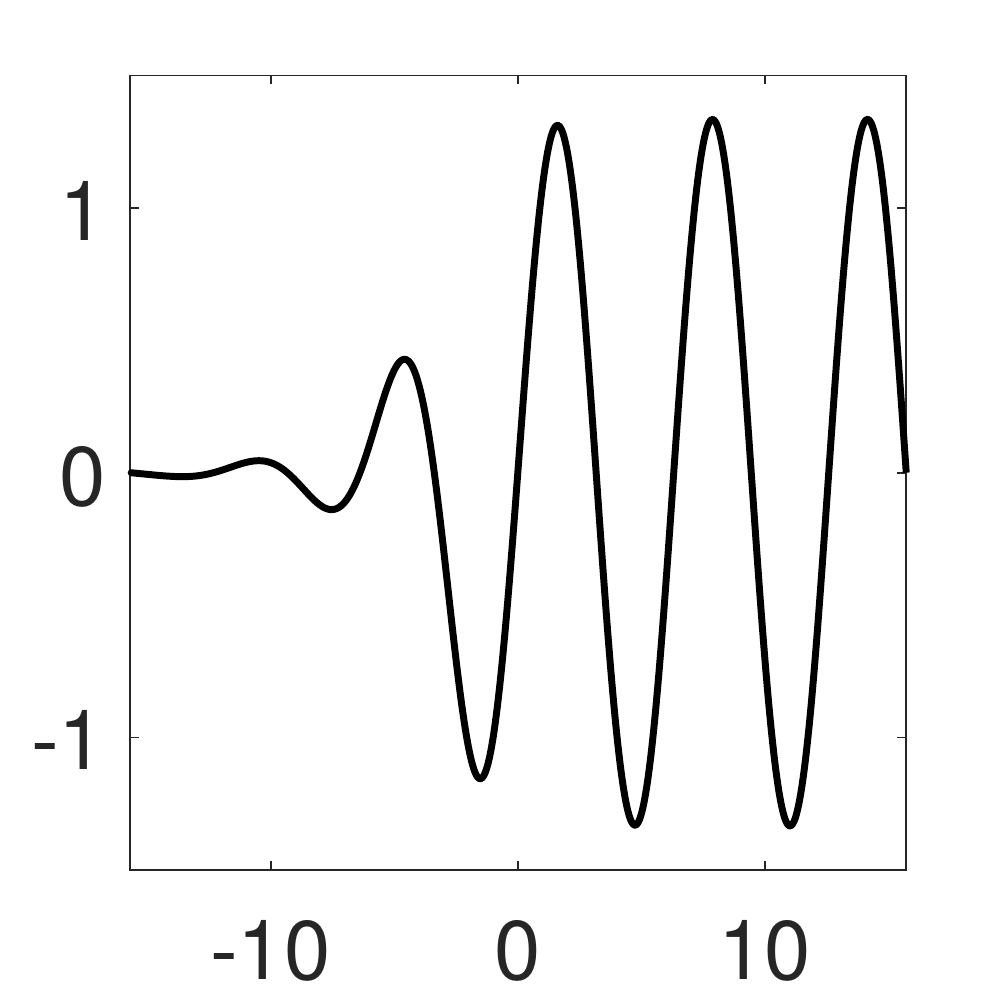}
     \put(-7,45){\rotatebox{90}{$u(x)$}}
     \put(100,5){$x$}
    \end{overpic}
     \caption{}
     \label{subfig-3:DSH}
\end{subfigure}\\

\hspace{2cm}
\begin{subfigure}[b]{0.2\textwidth}
     \begin{overpic}[width=\textwidth]{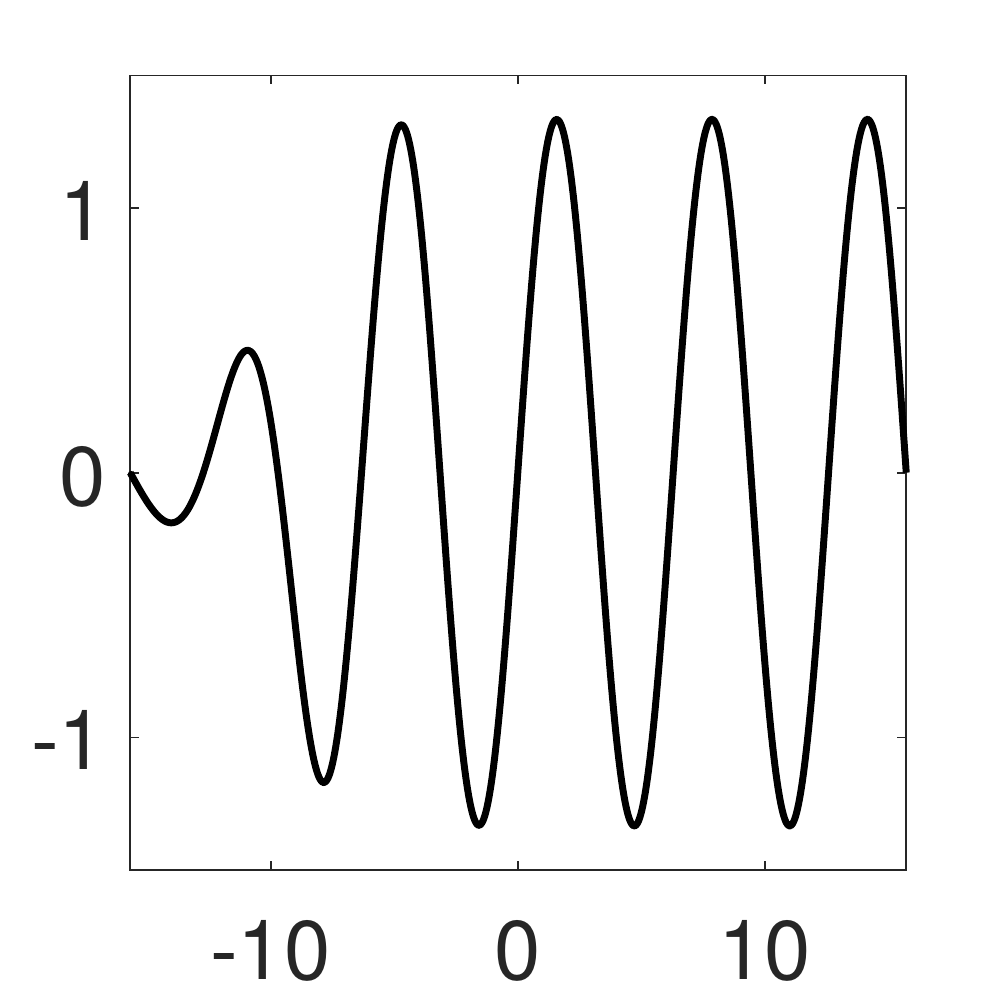}
     \put(-7,45){\rotatebox{90}{$u(x)$}}
     \put(100,5){$x$}
    \end{overpic}
     \caption{}
     \label{subfig-4:DSH}
\end{subfigure}
\end{multicols}
\vspace{-0.5cm}
\caption{Steady state bifurcation diagram of the Swift--Hohenberg equation~\eqref{eq:SHE} with $\nu=2$ on the domain $\Omega~=~(-5\pi, 5\pi)$ with homogeneous Dirichlet boundary conditions. Blue: branch of non-homogeneous stationary solutions with wavenumber $k=1$. Red: branch of front solutions illustrated in subfigures \subref{subfig-1:DSH}--\subref{subfig-4:DSH}. For sake of clarity only the first periodic and first snaking branch are shown. Thicker lines denote stable solution while thinner lines correspond to unstable ones. Circles and crosses indicate branch and fold points respectively.  (For interpretation of the references to color in this figure legend, the reader is referred to the web version of this article)}
\label{fig:SH-sn1-theory}
\end{figure}

As in the previous section, the bifurcation structure can be obtained numerically and is shown in Figure~\ref{fig:SH-sn1-theory} for $\Omega=(-5\pi,5\pi)$ and $\nu=2$. The (spatially) periodic patterns arise on the blue branch bifurcating subcritically from the homogeneous (black) branch (as shown in \cite[equation~9]{BURKE2007681} for~$\nu>0$). The subcriticality of the bifurcation is of fundamental importance for the appearance of the ``snaking'' (red) branch, which emerges on bounded domains from the branch of non-homogeneous stationary solutions close to the origin at $\mu<0$. This branch corresponds to front solutions shown in Figures \ref{subfig-1:DSH}--\ref{subfig-4:DSH}. At each turn in the snaking, one additional oscillation is added to the pattern until the domain is filled with oscillations. Then, the snaking branch reconnects close to the fold of the periodic branch.

Note that additional bifurcations further up the periodic branch lead to other snaking branches corresponding to front or localized patterns. These branches are not illustrated here for the sake of clarity.

\subsection{The fractional Swift--Hohenberg equation}\label{chapter:SHe-results}

We now study the steady state fractional Swift--Hohenberg equation with cubic--quintic non-linearity~\eqref{eq:SH-frac}, on a 1D domain of length $L$ with homogeneous Dirichlet boundary conditions. The homogeneous state $\bar{u}=0$ is a stationary solution of the fractional equation~\eqref{eq:SH-frac}. Further, according to \eqref{eq:spectral_definition} the Dirichlet spectral fractional Laplacian has eigenpairs
$(\phi_j, \lambda_j) = (\sin\left( {j\pi x}/{L} \right), -\left({j \pi}/{L} \right)^{2s})$ for~$j\geq 1$. The bifurcations from the homogeneous states occur at\\[-0.2cm]
\begin{equation}\label{eq:subsequent-bif}
\mu_j = \left(1-\left(\frac{j\pi}{L}\right)^{2s}\right)^2, \qquad j\geq 1.
\end{equation}

Thus, as in the standard problem, the homogeneous state $\bar{u}$ is stable for $\mu < 0 $ and first becomes unstable at $\mu_c=0$. In addition, note that as $s$ tends to zero, $\mu_j$ tends to $\mu_c=0$.

Substituting $(u_1,u_2) = (u, \Delta^s u)$ in the stationary equation~\eqref{eq:SH-frac}, we convert the Swift--Hohenberg equation into a system of second order PDEs
\begin{equation*}
\begin{pmatrix} 1 & 0\\ 0 & 0 \end{pmatrix} \partial_t \begin{pmatrix} u_1\\ u_2\end{pmatrix}=
    \begin{pmatrix}
    -\Delta^s u_2 -2u_2 -(1-\mu)u_1 + \nu u_1^3 - u_1^5\\
    \Delta^s u_1 - u_2
    \end{pmatrix},
\end{equation*}
which fit the \texttt{pde2path} framework based on finite elements.

To obtain the following numerical result, we fix $\nu=2$ and the domain $\Omega=(-5\pi, 5\pi)$ of length $L=10\pi$ with boundary conditions: $u(-5\pi) = u(5 \pi) = 0$ and $\Delta^su(-5 \pi) = 0 =\Delta^su(5 \pi) = 0$. We choose the meshsize $h=0.04$, corresponding to $n_p=786$ mesh points.

\begin{figure}
    \centering
    \begin{overpic}[width=0.45\textwidth]{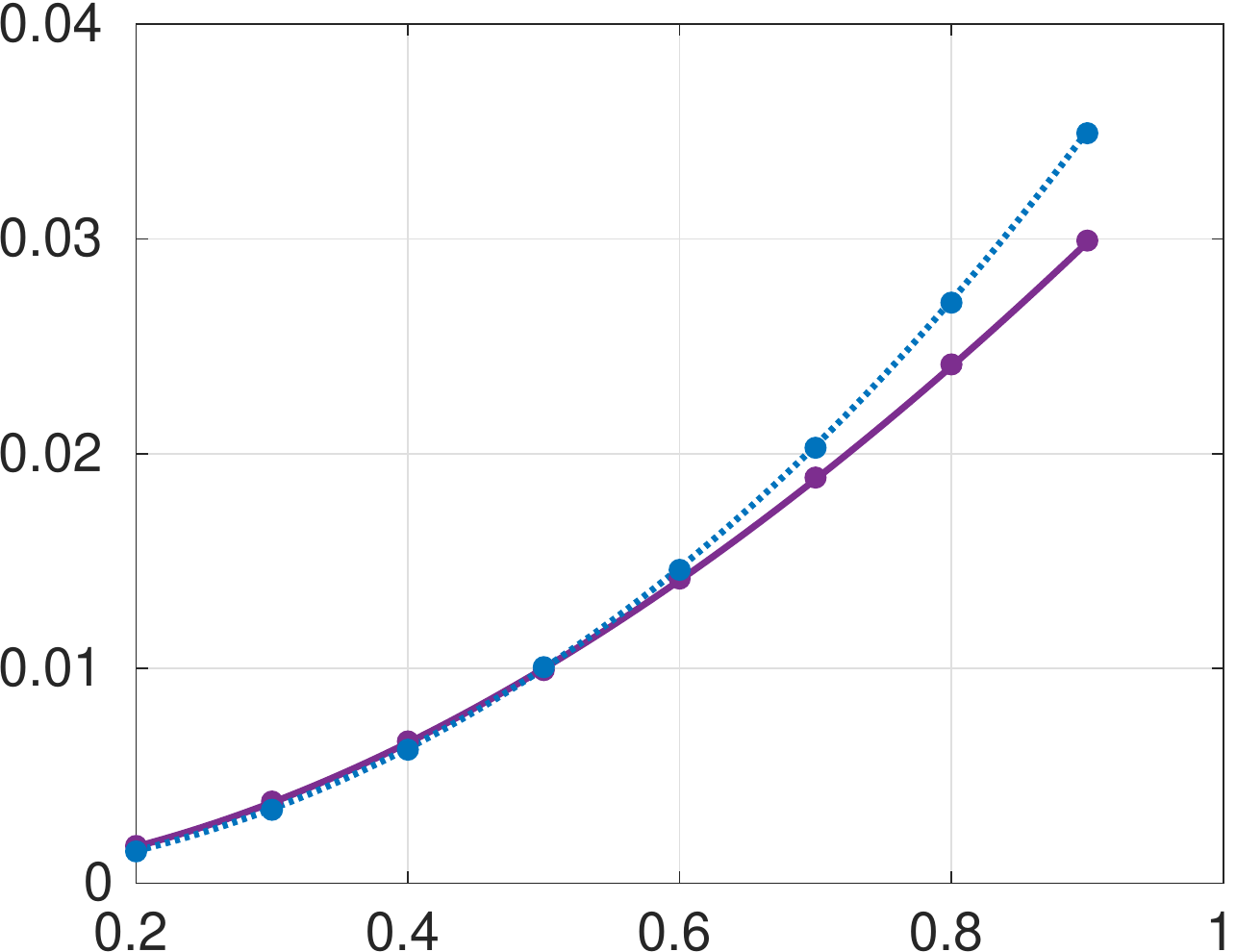}
    \put(100,5){$s$}
    \put(90,65){$\mu_3$}
    \put(90,57){$\mu_2$}
    \end{overpic}
    \caption{Evolution of the location of the second (continuous purple line) and third (dashed blue line) bifurcation points on the homogeneous branch (black in Figure~\ref{fig:sh-bif-firstImpression}) as $s$ decreases, for the fractional Swift--Hohenberg equation~\eqref{eq:SH-frac} on a bounded domain with homogeneous Dirichlet boundary conditions. Dots correspond to numerical values, while lines correspond to analytical values obtained via equation~\eqref{eq:subsequent-bif}. (For interpretation of the references to color in this figure legend, the reader is referred to the web version of this article)}
    \label{fig:sh-bif1and2-shifting}
    \end{figure}

As mentioned above, the bifurcation points on the homogeneous branch accumulate at $\mu_c=0$ as the fractional order $s$ tends to $0$. Figure~\ref{fig:sh-bif1and2-shifting} shows, as an example, the shifting of the second and third bifurcation points, corresponding to kernel vectors $(\sin\left( {9\pi x}/{L} \right), -\left({9 \pi}/{L} \right)^{2s})$ and $(\sin\left( {11\pi x}/{L} \right), -\left({11 \pi}/{L} \right)^{2s})$ respectively.
The values obtained numerically are marked with dots and the analytical values, computed using formula~\eqref{eq:subsequent-bif} with $L=10\pi$, are showed as lines on the figure. We find that the maximum relative error is on the order of $10^{-2}$ for both points and observe that the numerical values indeed coincide with the ones expected from the analysis. In addition, we observe that the two bifurcation points interchange their positions at $s=0.5$.

\begin{figure}
 \centering
 \begin{subfigure}[b]{0.45\textwidth}
         \centering
         \begin{overpic}[width=\textwidth]{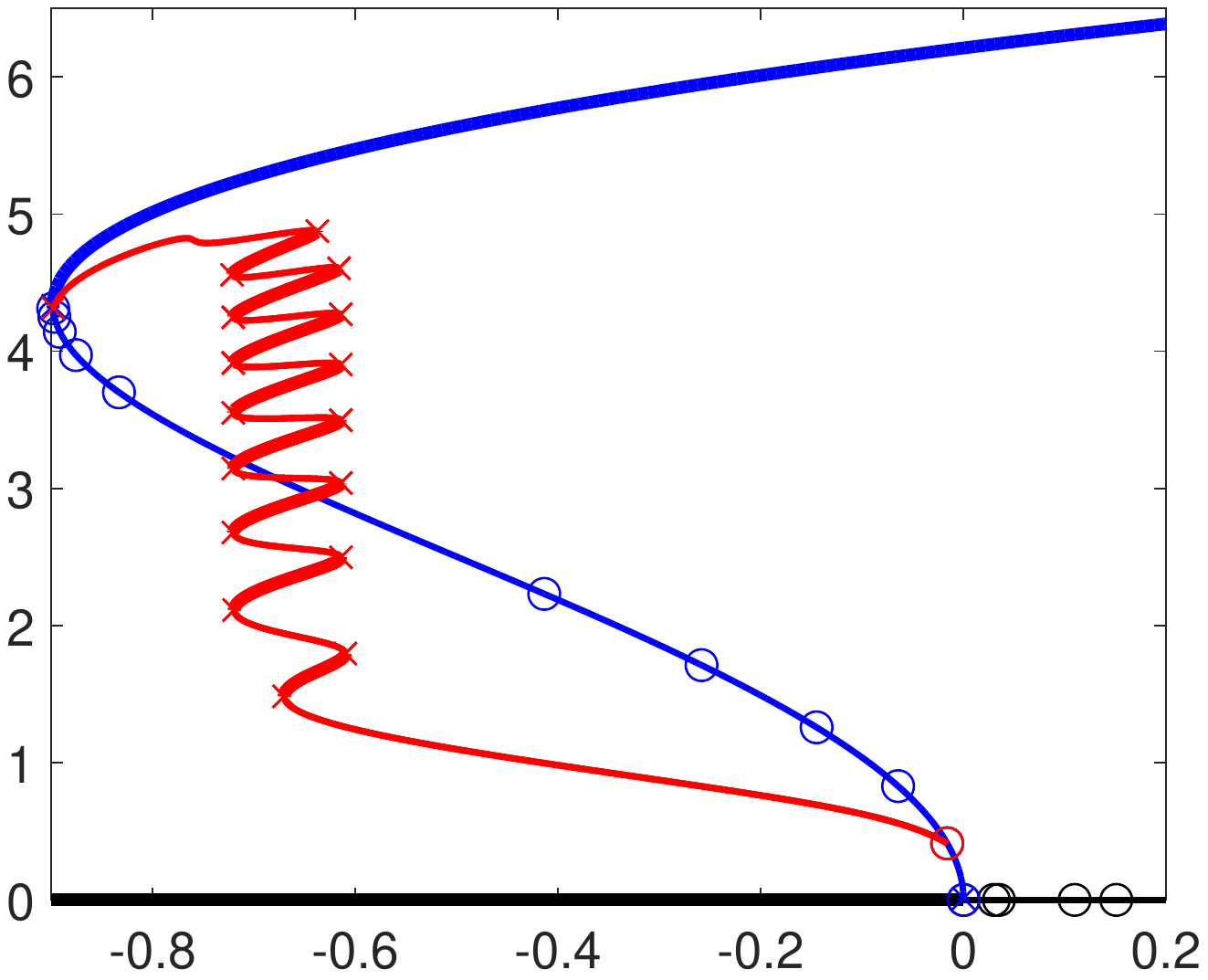}
         \put(-8,45){\rotatebox{90}{$\|u\|_{L^2}$}}
         \put(99,5){$\mu$}
         \end{overpic}
         \caption{$s=0.9$}\label{fig:sh-1Dbif-s09}
 \end{subfigure}
 \hspace{1cm}
 \begin{subfigure}[b]{0.45\textwidth}
         \centering
         \begin{overpic}[width=\textwidth]{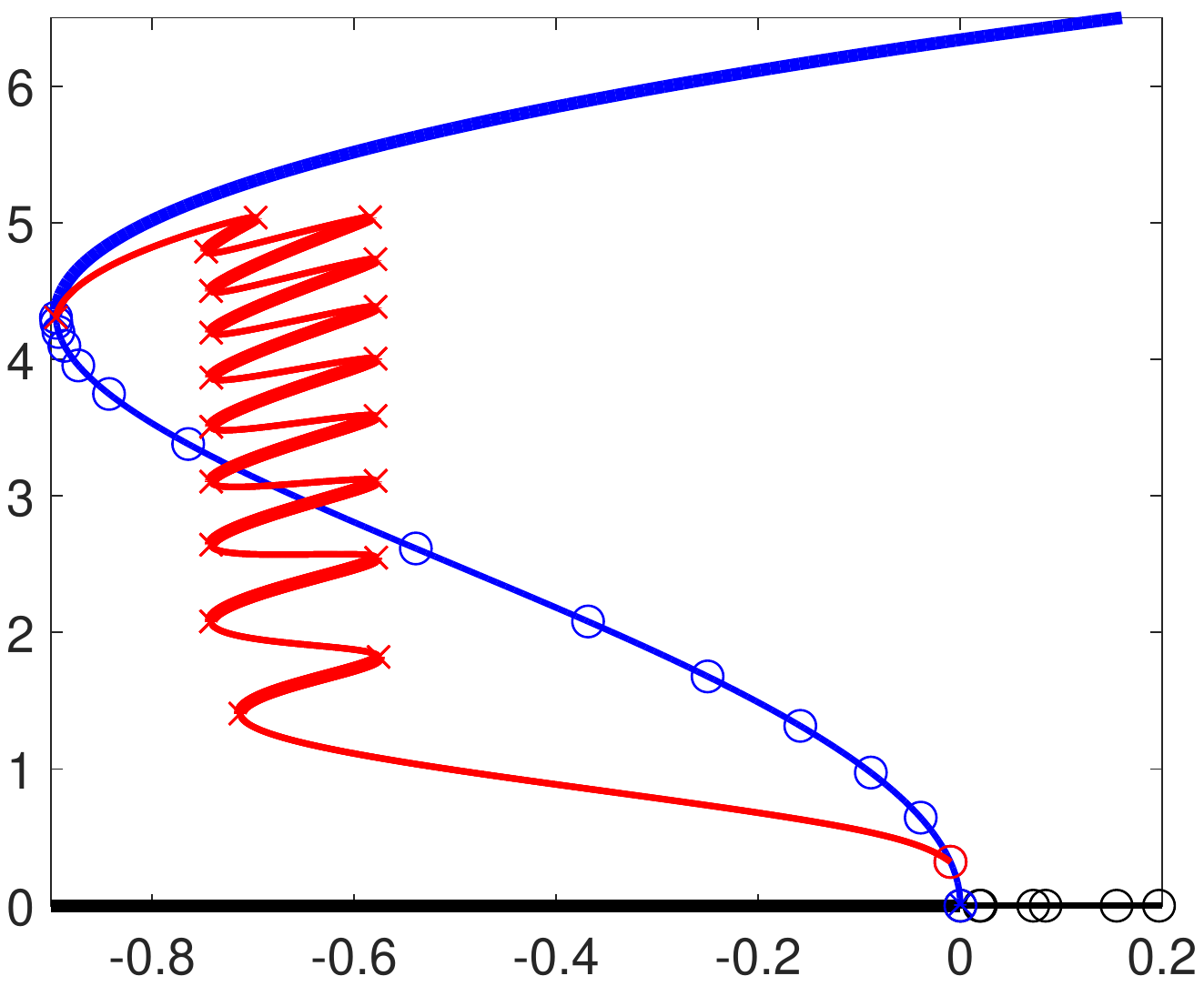}
         \put(-8,45){\rotatebox{90}{$\|u\|_{L^2}$}}
         \put(99,5){$\mu$}
         \end{overpic}
         \caption{$s=0.7$}\label{fig:sh-1Dbif-s07}
 \end{subfigure} \\[0.5cm]
 \begin{subfigure}[b]{0.45\textwidth}
         \centering
         \begin{overpic}[width=\textwidth]{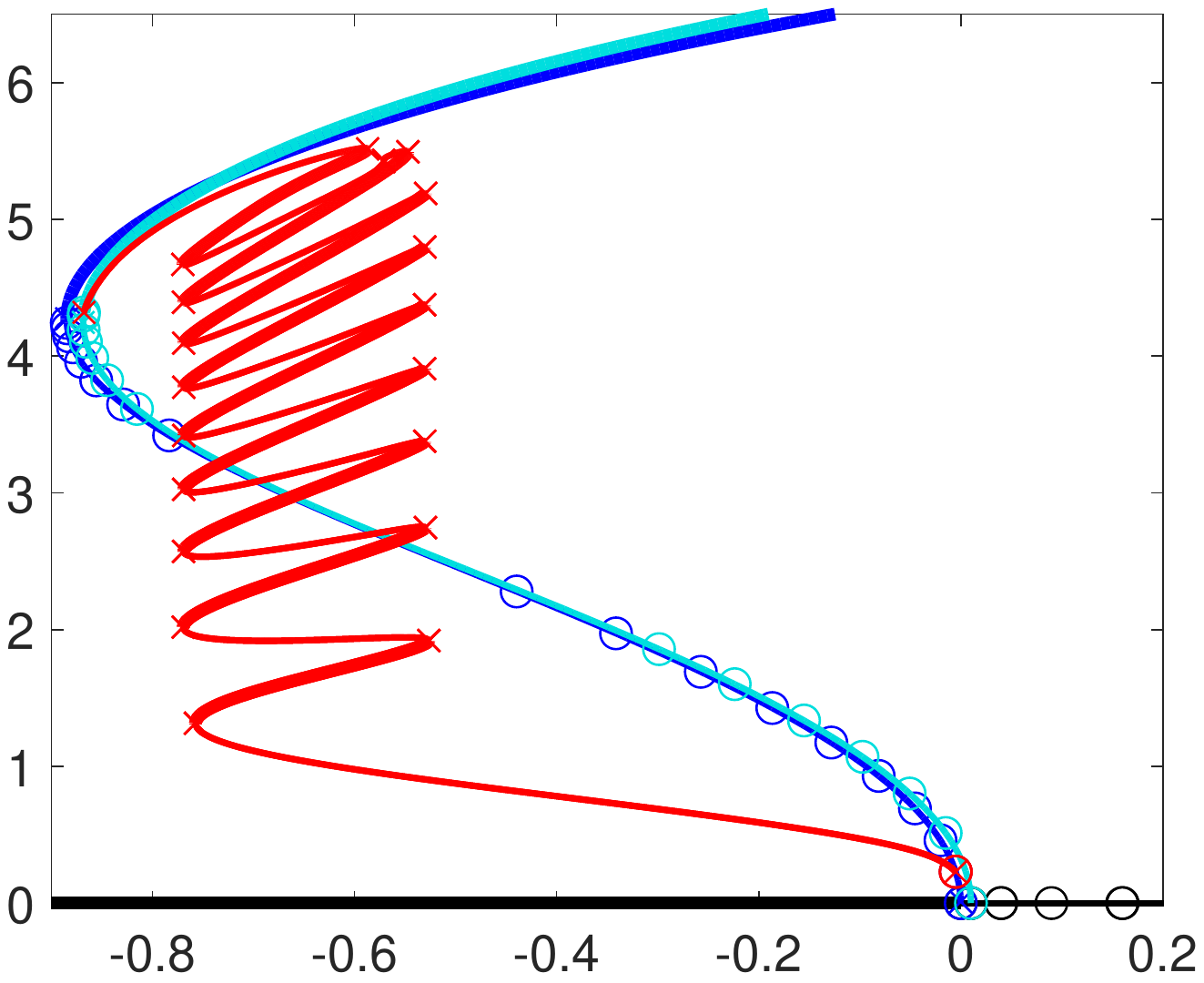}
         \put(-8,45){\rotatebox{90}{$\|u\|_{L^2}$}}
         \put(99,5){$\mu$}
         \end{overpic}
         \caption{$s=0.5$}\label{fig:sh-1Dbif-s05}
 \end{subfigure}
  \hspace{1cm}
  \begin{subfigure}[b]{0.45\textwidth}
         \centering
         \begin{overpic}[width=\textwidth]{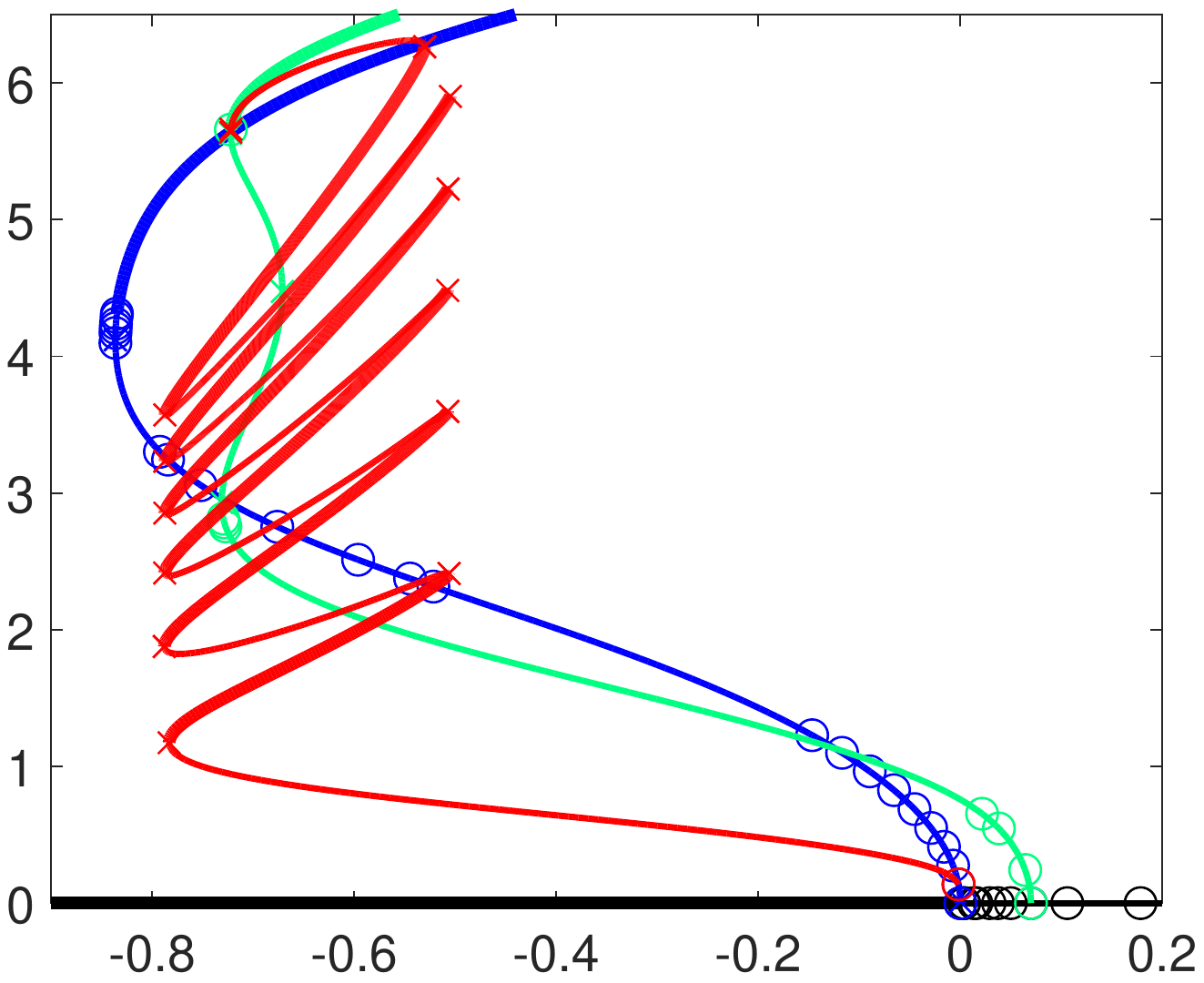}
         \put(-8,45){\rotatebox{90}{$\|u\|_{L^2}$}}
         \put(99,5){$\mu$}
         \end{overpic}
         \caption{$s=0.3$}\label{fig:sh-1Dbif-s03}
 \end{subfigure}
 \caption{Evolution of the bifurcation structure with respect to the bifurcation parameter $\mu$ of the fractional Swift--Hohenberg equation~\eqref{eq:SH-frac} with $\nu=2$ on a bounded domain $\Omega=(-5\pi, 5\pi)$ with homogeneous Dirichlet boundaries as the order of the fractional Laplacian decreases. We clearly see that the snaking pattern emerges for different fractional orders, yet it gets stretched (``bloated'') as the fractional order decreases. Thick and thin lines denote stable and unstable solutions, while circles and crosses indicate branch and fold points respectively.  (For interpretation of the references to color in this figure legend, the reader is referred to the web version of this article)}
 \label{fig:sh-bif-firstImpression}
\end{figure}

In order to get a general impression of the effects of fractional diffusion on the bifurcation structure of problem~\eqref{eq:SH-frac}, Figure~\ref{fig:sh-bif-firstImpression} shows the bifurcation diagram for four different fractional orders (namely $s = 0.9$, $s = 0.7$,  $s = 0.5$ and $s = 0.3$), which can be compared to the bifurcation structure of the standard Swift--Hohenberg equation presented in Section \ref{sec:SHe-theory}. The homogeneous branch is shown in black, the first branch bifurcating at $\mu=0$ is shown in blue and the first snaking branch in red. As we already said, bifurcation points relative to other wave numbers shift towards the primary bifurcation. Looking at the computed bifurcation structure, we notice two other important changes as $s$ decreases. Firstly, the width of the snaking branch is significantly increasing. Secondly, for $s=0.5$ and $s=0.3$ the snaking branch does not reconnect to the blue branch anymore, but to other branches with different wave numbers. Note that this behaviour also occurs for standard (Laplacian) Swift--Hohenberg equations varying the domain length \cite{bergeon2008eckhaus}.

\begin{figure}

\begin{multicols}{2}
\begin{subfigure}[t]{0.35\textwidth}
\centering
\hspace{-2cm}
\begin{overpic}[width=0.9\textwidth]{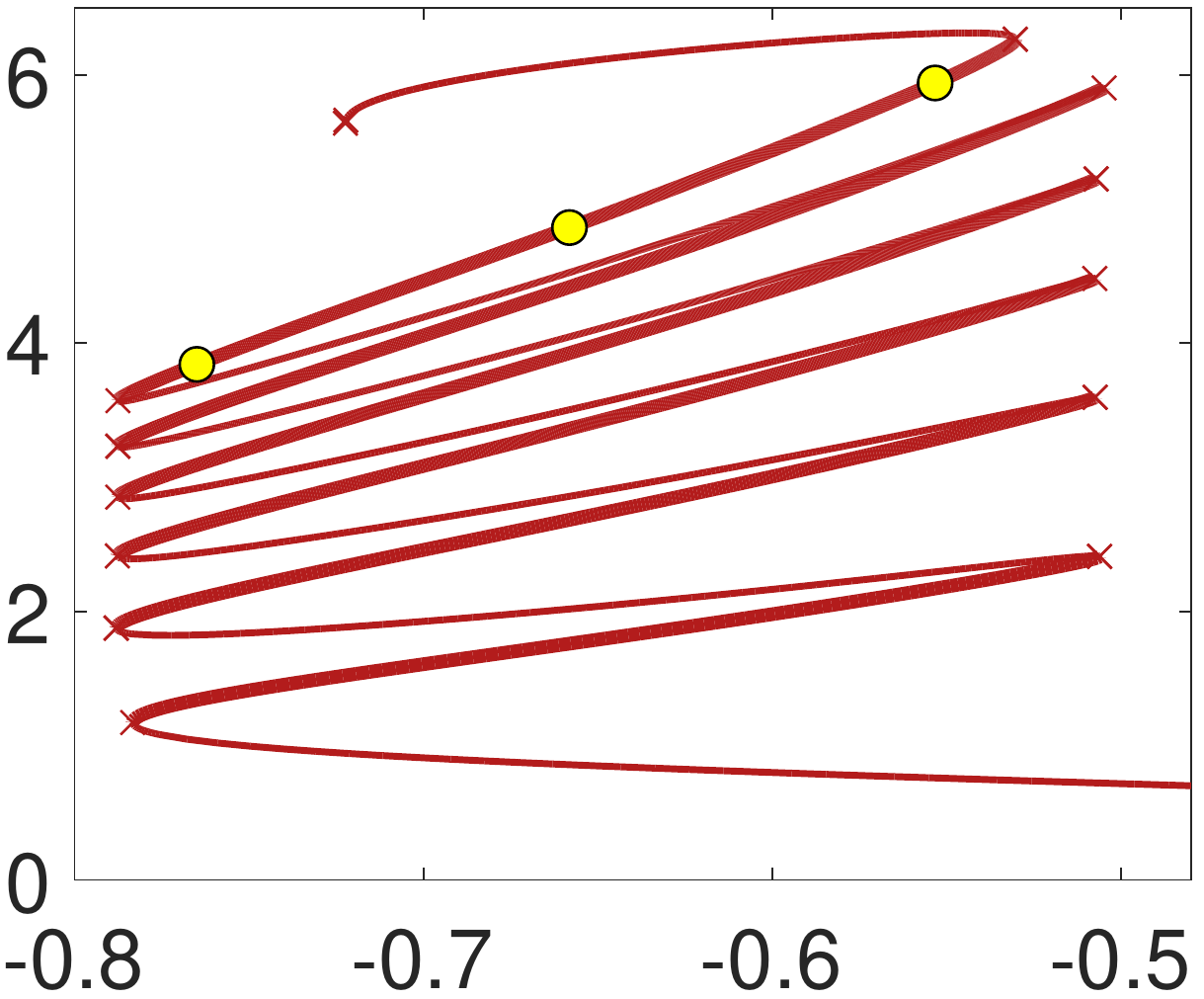}
\put(-10,45){\rotatebox{90}{$\|u\|_{L^2}$}}
\put(100,5){$\mu$}

\put(10,58){(a)}
\put(40,67){(b)}
\put(65,75){(c)}
\end{overpic}
\label{fig:sh-s03-tilted-labeled}
\end{subfigure}
\newpage
\vspace*{\fill}
\hspace{-2.5cm}
\begin{subfigure}[t]{0.2\textwidth}
 \centering
 \begin{overpic}[width=\textwidth]{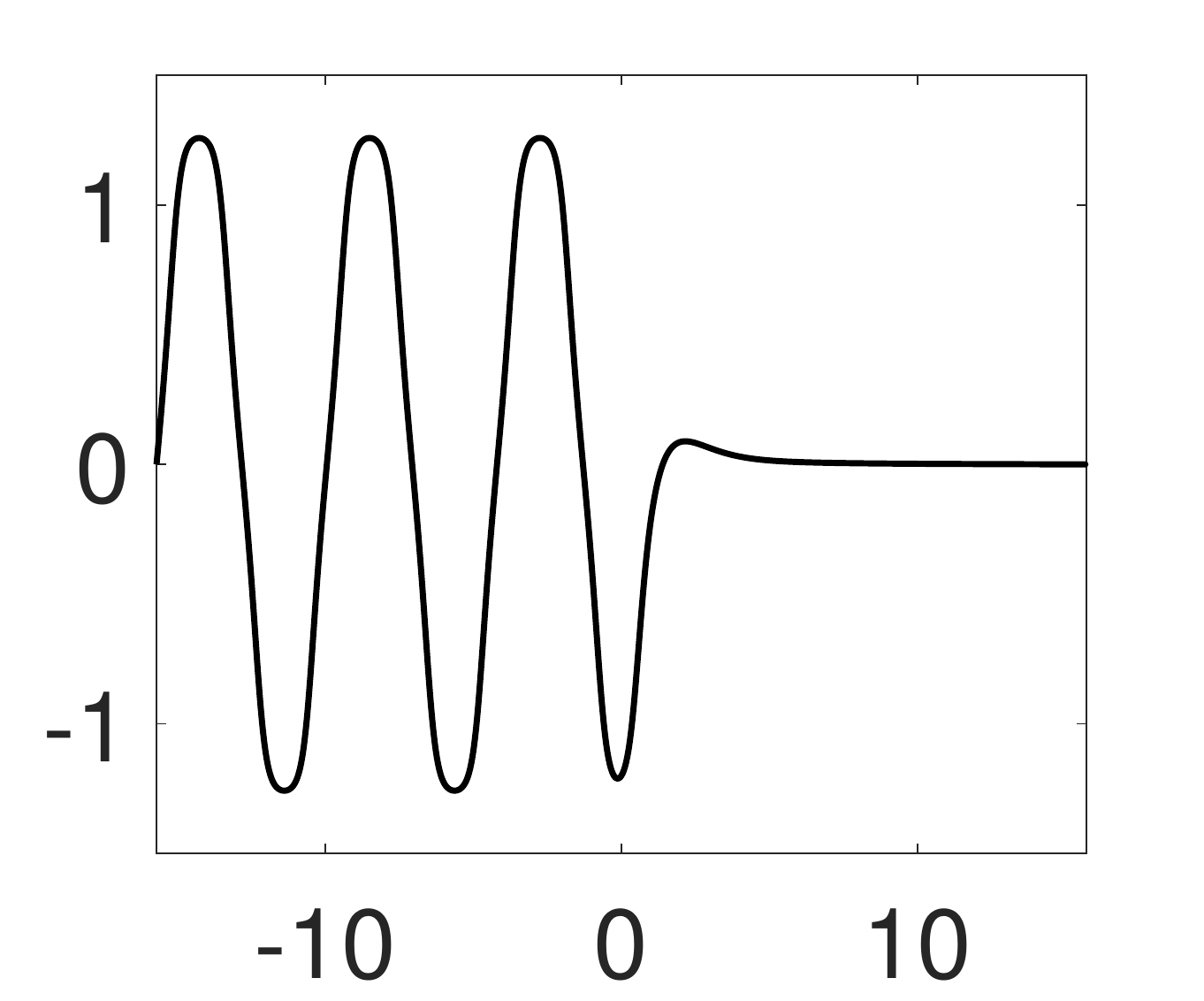}
 \put(-7,45){\rotatebox{90}{$u(x)$}}
 \put(95,5){$x$}
 \end{overpic}
 \caption{}\label{fig:sh-s03-pt5300}
\end{subfigure}
\begin{subfigure}[t]{0.2\textwidth}
 \centering
 \begin{overpic}[width=\textwidth]{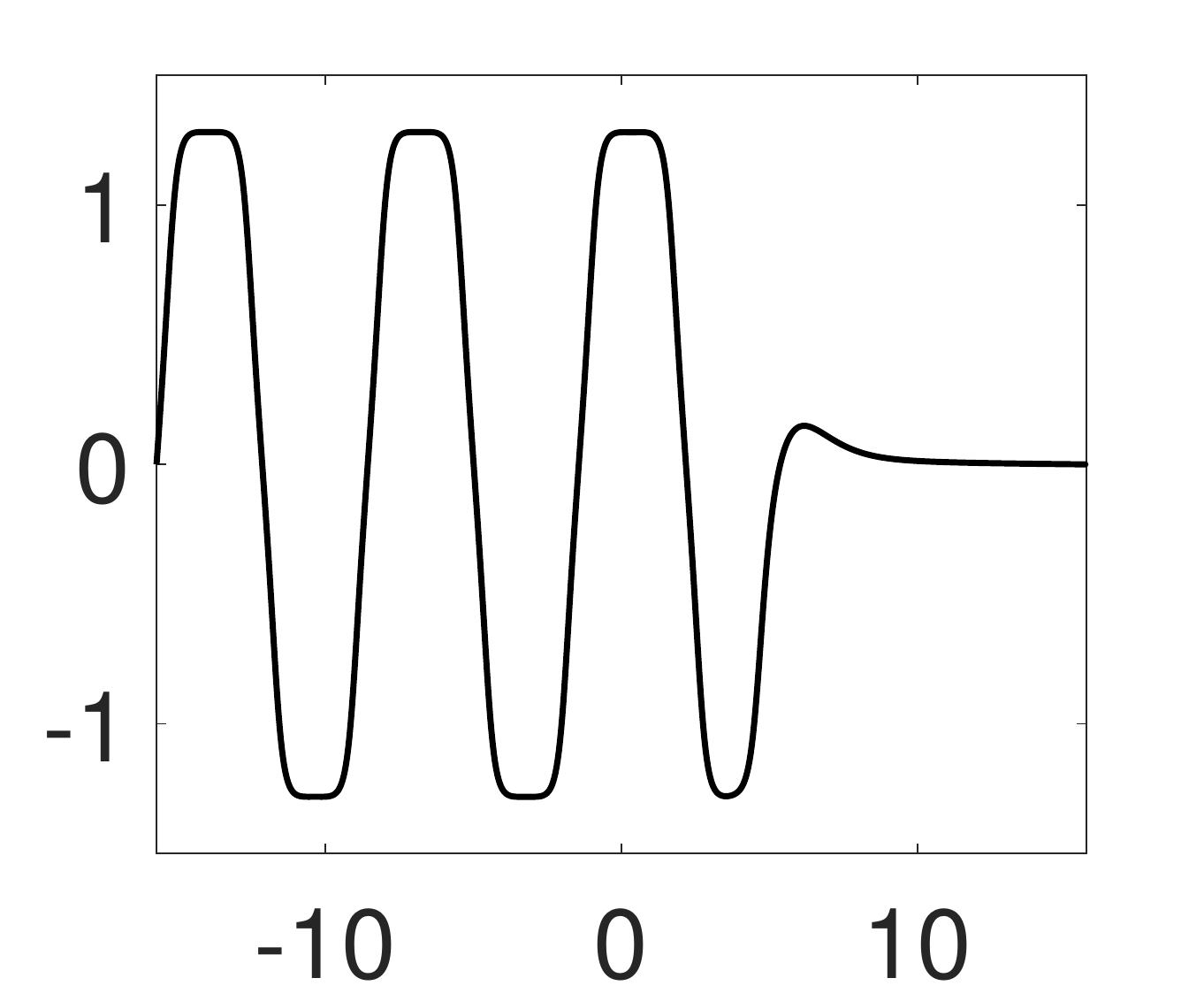}
 \put(95,5){$x$}
 \end{overpic}
 \caption{}\label{fig:sh-s03-pt5600}
\end{subfigure}
\begin{subfigure}[t]{0.2\textwidth}
 \centering
 \begin{overpic}[width=\textwidth]{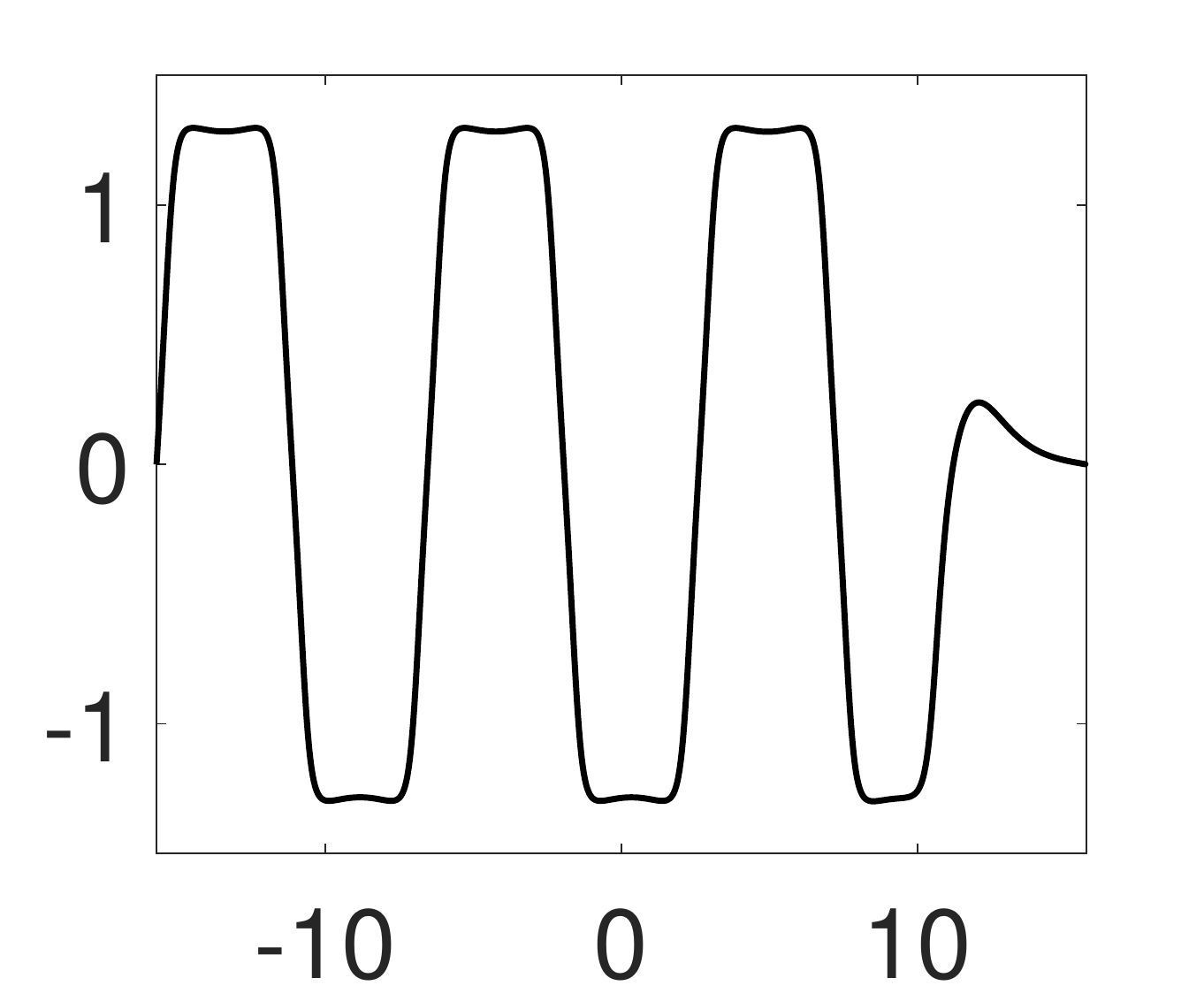}
 \put(95,5){$x$}
 \end{overpic}
 \caption{}\label{fig:sh-s03-pt5900}
\end{subfigure}
\vspace*{\fill}
\end{multicols}
\vspace{-0.5cm}
\caption{Solutions along the upper part of the snaking branch for $s=0.3$ in the fractional Swift--Hohenberg equation~\eqref{eq:SH-frac} on a bounded domain with homogeneous Dirichlet boundary conditions.}
\label{fig:sh-tilted-branch-s03}
\end{figure}

In addition to the increase in width, one can see in Figure~\ref{fig:sh-bif-firstImpression} that the back and forth oscillations forming the snaking tend to tilt on the left for smaller~$s$. This phenomenon is particularly visible for $s=0.3$. The solutions along the upper part of the snaking branch are illustrated in Figure~\ref{fig:sh-tilted-branch-s03}. We see that the oscillations significantly enlarge as we move up the snaking branch. They would then narrow again after the next turn in an ``accordion'' effect. The deformation of the snaking branch of the snaking branch is not a known effect for the fractional Laplacian but is known to occur for other nonlocal Swift--Hohenberg equations \cite{morgan2014SH}. Further, counting the number of turns on the snaking branches in Figure~\ref{fig:sh-bif-firstImpression}, we see that whereas for $s=0.9$ we have $8$ pair of turns, for $s=0.7$ and $s=0.5$ we have $9$ pairs and for $s=0.3$ only $6$. As shown in Figure~\ref{fig:sh-tilted-branch-s03}, the first snaking branch corresponds to front solutions. At each snaking turn, one oscillation is added to the front. Thus, the number of turns is directly related to the final number of oscillations in the solution. This is linked to the next observation which regards the reconnection of the snaking to a branch of non-homogeneous stationary solutions.

\begin{figure}
\begin{multicols}{2}
\centering
\vspace*{\fill}
\hspace{-2cm}
\begin{subfigure}[t]{0.25\textwidth}
 \centering
 \begin{overpic}[width=\textwidth]{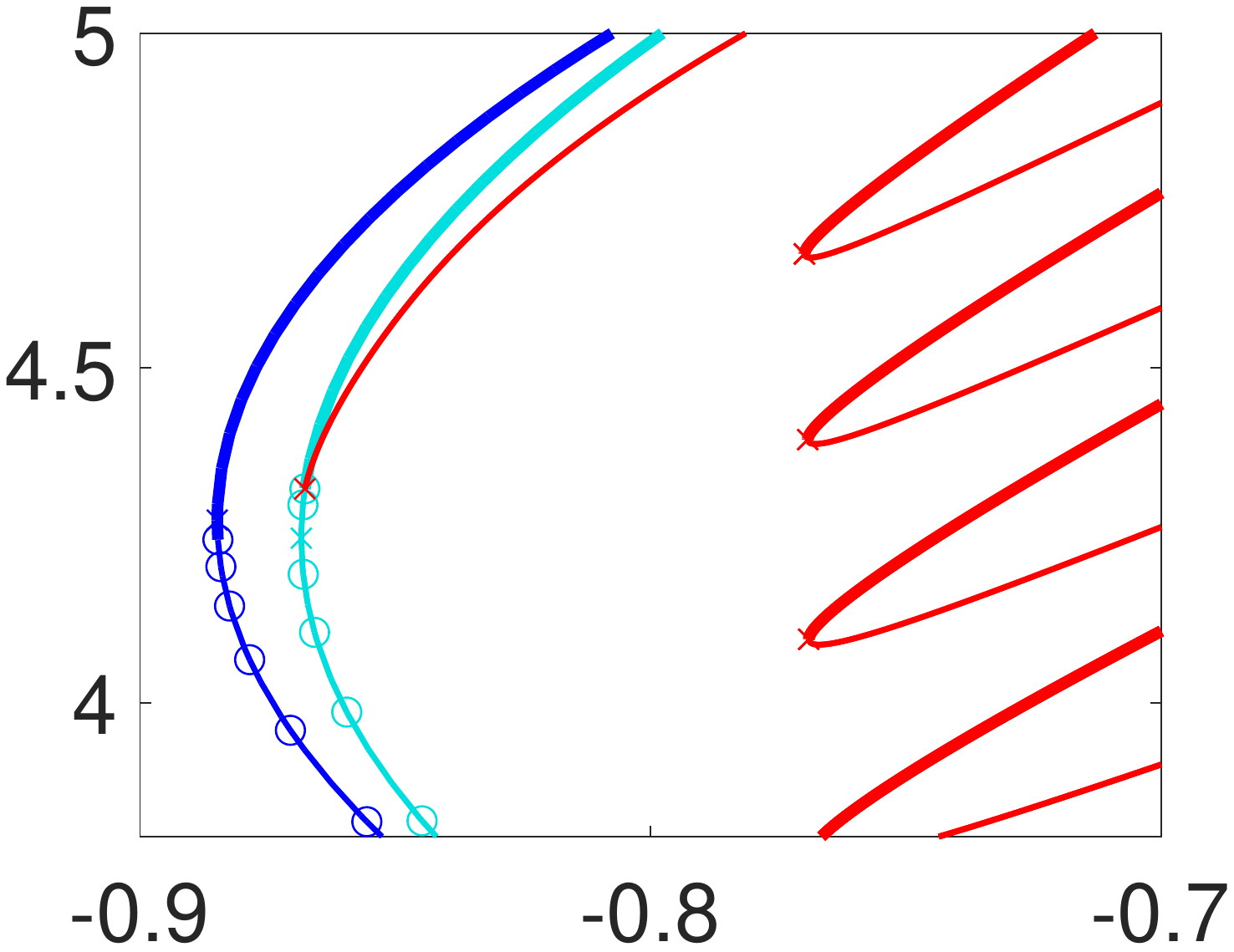}
 \put(-10,35){\rotatebox{90}{$\|u\|_{L^2}$}}
 \put(100,5){$\mu$}
 \end{overpic}
 \caption{}\label{subfig:sh-reconnection-zoom}
\end{subfigure}\\
\hspace{-2cm}
\begin{subfigure}[t]{0.25\textwidth}
 \centering
 \begin{overpic}[width=\textwidth]{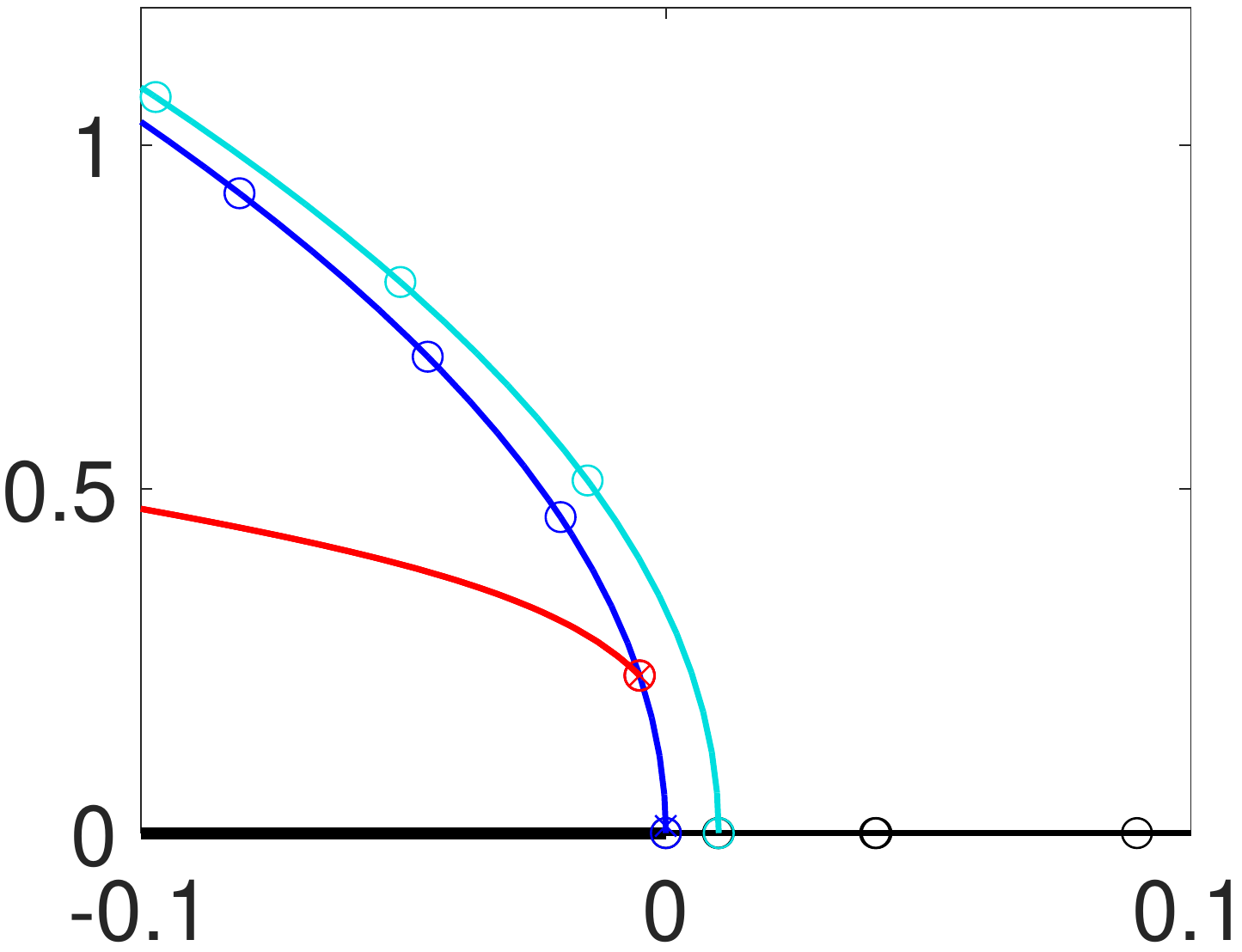}
 \put(-10,35){\rotatebox{90}{$\|u\|_{L^2}$}}
 \put(100,5){$\mu$}
 \end{overpic}
 \caption{}\label{subfig:sh-reconnection-zoom2}
\end{subfigure}
\vspace*{\fill}
\newpage
\vspace*{\fill}
\hspace{-3cm}
\begin{subfigure}[t]{0.45\textwidth}
 \centering
 \begin{overpic}[width=\textwidth]{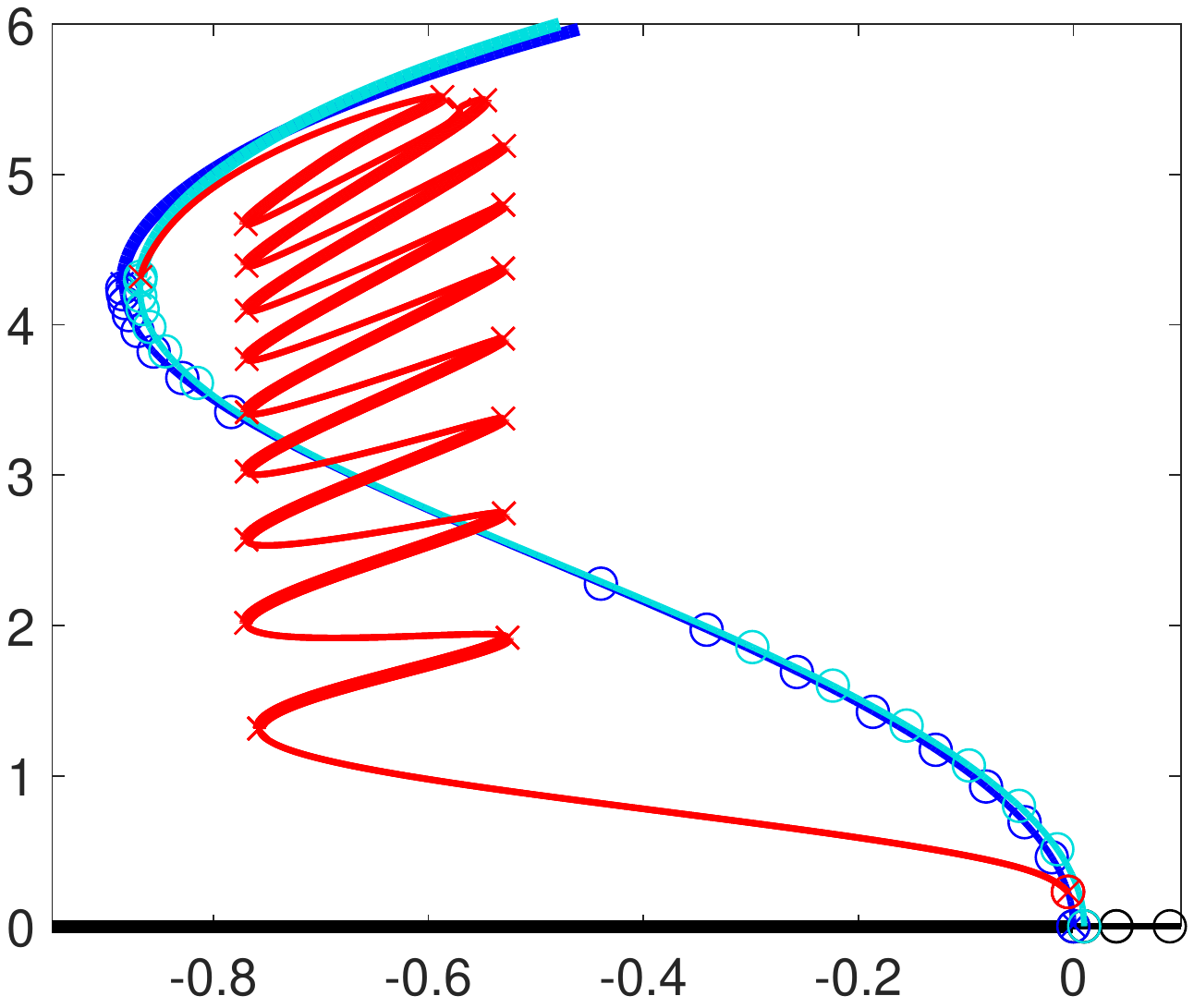}
 \put(4,58.5){\fbox{\rule{0.6in}{0pt}\rule[-2.4ex]{0pt}{6.8ex}}}
 \put(75,14){\fbox{\rule{0.6in}{0pt}\rule[-2.4ex]{0pt}{6.8ex}}}
 \put(-7,35){\rotatebox{90}{$\|u\|_{L^2}$}}
 \put(100,5){$\mu$}
 \end{overpic}
 \caption{}\label{subfig:sh-reconnection-bd}
\end{subfigure}\\
\vspace*{\fill}
\end{multicols}
\vspace{-0.5cm}
\caption{Magnified views of the reconnection (\subref{subfig:sh-reconnection-zoom}) and starting points (\subref{subfig:sh-reconnection-zoom2}) of the snaking branch of the fractional Swift--Hohenberg equation~\eqref{eq:SH-frac} on a bounded domain with homogeneous Dirichlet boundary conditions for fractional order $s=0.5$. The snaking branch in red originating from the first bifurcation point on the blue branch reconnects to the blue-green branch originating from the third bifurcation point on the homogeneous branch (black). The bifurcation diagram in Subfigure \subref{subfig:sh-reconnection-bd} is provided for reference purpose only.  (For interpretation of the references to color in this figure legend, the reader is referred to the web version of this article)}
\label{fig:sh-reconnection-plus-zoom}
\end{figure}

In fact, as $s$ decreases the snaking branch does not reconnect to the first periodic branch anymore but to periodic branches originating from subsequent bifurcations on the homogeneous branch. For instance, for $s=0.5$ and $s=0.3$ it connects to the branches originating from the the third and ninth bifurcation points, respectively, on the homogeneous (black) branch. Magnified views of the bifurcation diagram for $s=0.5$ close to the reconnection and the starting points are reported in Figure~\ref{fig:sh-reconnection-plus-zoom}.

\begin{figure}
\centering
 \begin{overpic}[width=0.4\textwidth]{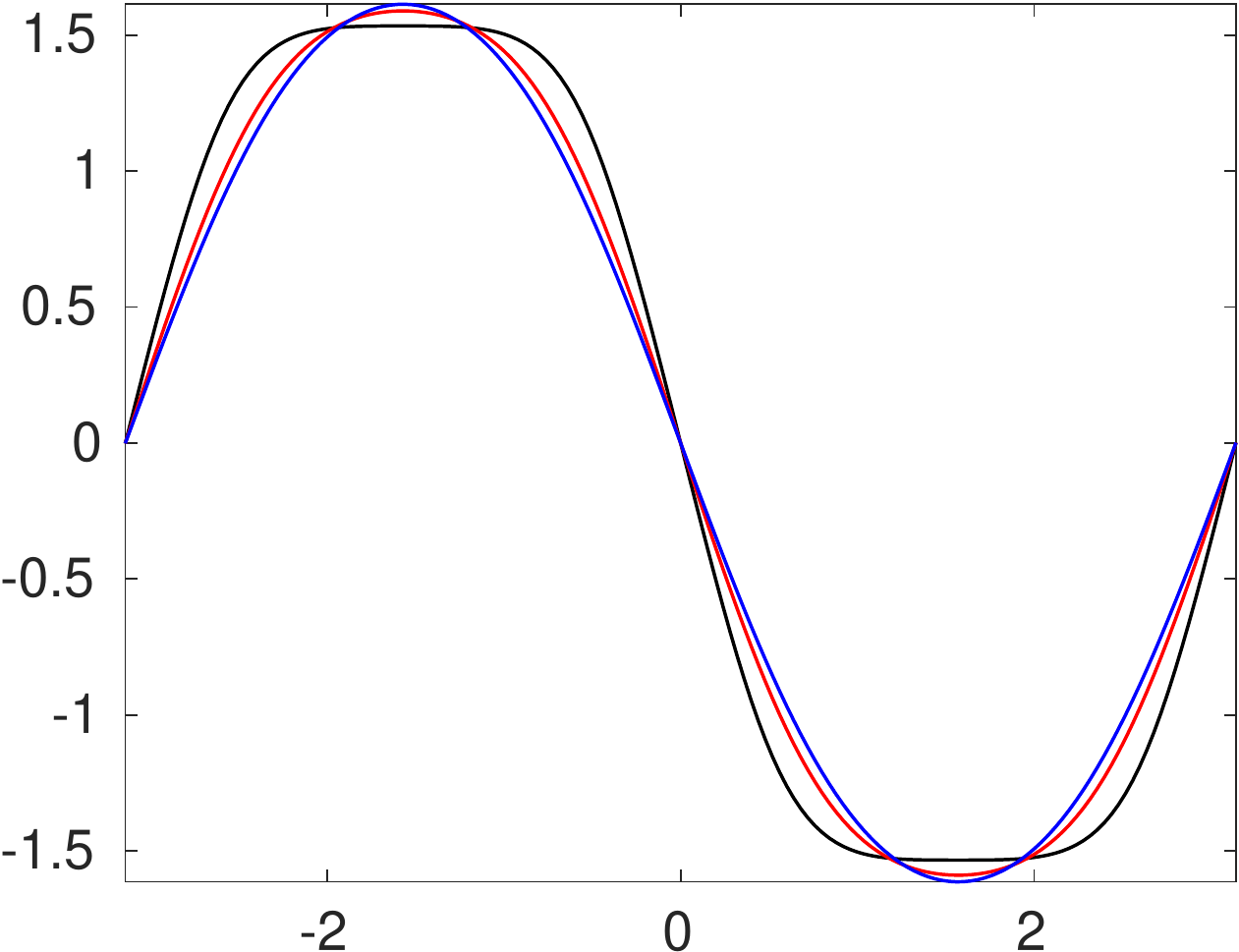}
 \put(-5,35){\rotatebox{90}{$u(x)$}}
 \put(100,5){$x$}
 \end{overpic}
\caption{Solution profiles to the fractional Swift--Hohenberg equation~\eqref{eq:SH-frac} on a bounded domain with homogeneous Dirichlet boundary conditions for $\mu=-0.4$ at different fractional orders: $s=0.9$ in blue, $s=0.5$ in red and $s=0.3$ in black. (For interpretation of the references to color in this figure legend, the reader is referred to the web version of this article)}
\label{fig:sh-sol-profiles-0907005}
\end{figure}

Finally, a last comment has to be made on the shape of the solutions. Looking at Figure~\ref{fig:sh-bif-firstImpression}, one can notice that, as the fractional order $s$ decreases, the non-homogeneous steady states (blue-branch) have slightly larger $L^2$-norm values. The solutions indeed tend to ``sharpen'' for smaller values of $s$ closer resembling a square-profile. This is illustrated in Figure~\ref{fig:sh-sol-profiles-0907005}, where the solution profiles on the first (blue) branch at $\mu=4$ are shown for different values of the fractional order $s$. We see that, whereas the solution at $s=0.9$ and $s=0.7$ do not differ much, at $s=0.5$ the profile significantly ``flattens'' towards the top of oscillation. Note that this pattern has also been observed in the previous section for the periodic solutions of the fractional Allen--Cahn equation.

\FloatBarrier

\section{The (fractional) Schnakenberg system}\label{sec:Schnaks-all}
Finally, we turn to the Schnakenberg system. The next subsection (\ref{schnak}) is devoted to the standard version of the Schnakenberg system. The standard steady state bifurcation diagrams and typical solutions are obtained using the continuation software \texttt{pde2path} and were already shown in \cite{uecker_pde2path_2014}. Subsection \ref{chapter:Schnakenberg-system-results} then presents the results for the fractional version of the system, exploiting the new version of \texttt{pde2path} adapted to treat fractional problems, see section \ref{sec:numerics}. \\

\subsection{The Schnakenberg system with standard Laplacian}\label{schnak}
The Schnakenberg system has been used to model the spatial distribution of a morphogen for tissue patterning, and it presents the activator-inhibitor structure crucial in Turing instability. The system reads
\begin{equation}
    \partial_t \begin{pmatrix} u_1 \\ u_2 \end{pmatrix}  =\begin{pmatrix} 1 & 0 \\ 0 & d \end{pmatrix}  \begin{pmatrix}  \Delta u_1 \\  \Delta u_2 \end{pmatrix} + F(u_1,u_2;\mu),
    \label{eq:Schnakenberg}
\end{equation}
where the reaction part $F$ is given in~\eqref{eq:F_m} and also depends on the parameter $\sigma$. The homogeneous state $\bar{u} = (\mu, 1/\mu)$ is a steady state solution of system $\eqref{eq:Schnakenberg}$, independent of the value of $\sigma$ in the reaction part. Further, following~\cite{Doelman2019}, we find that~$\bar{u}$ is stable for parameter values $\mu>\mu_c$ and that at $\mu_c$ non-homogeneous steady states bifurcate from the homogeneous branch with critical wavenumber $k_c$, where
\begin{equation*}
\mu_c = \sqrt{d(3-\sqrt{8})}, \qquad k_c = \sqrt{\frac{d-\mu_c^2}{2d}}=\sqrt{\sqrt{2}-1}.
\end{equation*}

\begin{figure}
    \centering
    \begin{overpic}[width=0.5\textwidth,trim=20 20 0 0,clip]{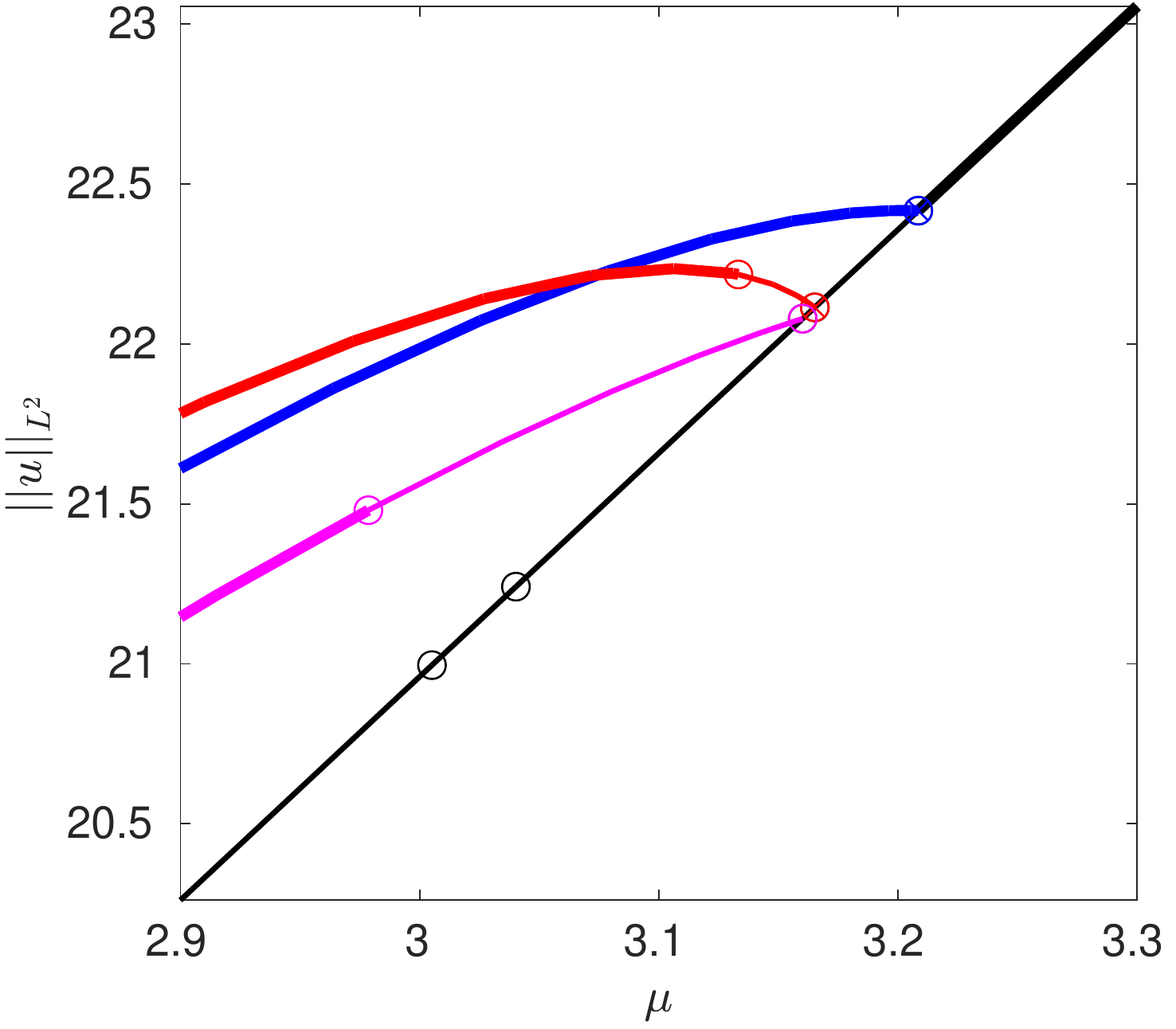}
    \put(-7,45){\rotatebox{90}{$\|u_1\|_{L^2}$}}
     \put(85,0){$\mu$}
    \end{overpic}
    \caption{Steady state bifurcation diagram of the Schnakenberg equation~\eqref{eq:Schnakenberg} with~$d=~60$ and~$\sigma=0$ in the reaction part~\eqref{eq:F_m}, on a bounded domain with length $L={10\pi}/{k_c}$, with homogeneous Neumann boundary conditions. For sake of clarity, only the first three branches bifurcating from the homogeneous states (black) are shown. Thicker lines denote stable solution while thinner lines correspond to unstable ones. Circles indicate branch points. (For interpretation of the references to color in this figure legend, the reader is referred to the web version of this article)}
    \label{fig:Schnak-bifurcations-theory}
 \end{figure}

On a finite domain, the first bifurcation is followed by a discrete set of further bifurcations in $\mu < \mu_c$ as further modes become unstable. The bifurcation diagram for the steady states of~\eqref{eq:Schnakenberg} with $\sigma=0$ is shown in Figure~\ref{fig:Schnak-bifurcations-theory}, considering the domain $\Omega=(-5\pi/k_c,5\pi/k_c)$ (chosen to accommodate the basic periodic patterns bifurcating from the homogeneous branch). The non-homogeneous (blue) branch bifurcates supercritically at $\mu_c$ from the homogeneous (black) branch. Two further branches are shown in red and magenta.

\begin{figure}
\begin{multicols}{2}
\hspace{0.5cm}
\begin{subfigure}[b]{0.5\textwidth}
    \centering
    \begin{overpic}[width=\textwidth,clip]{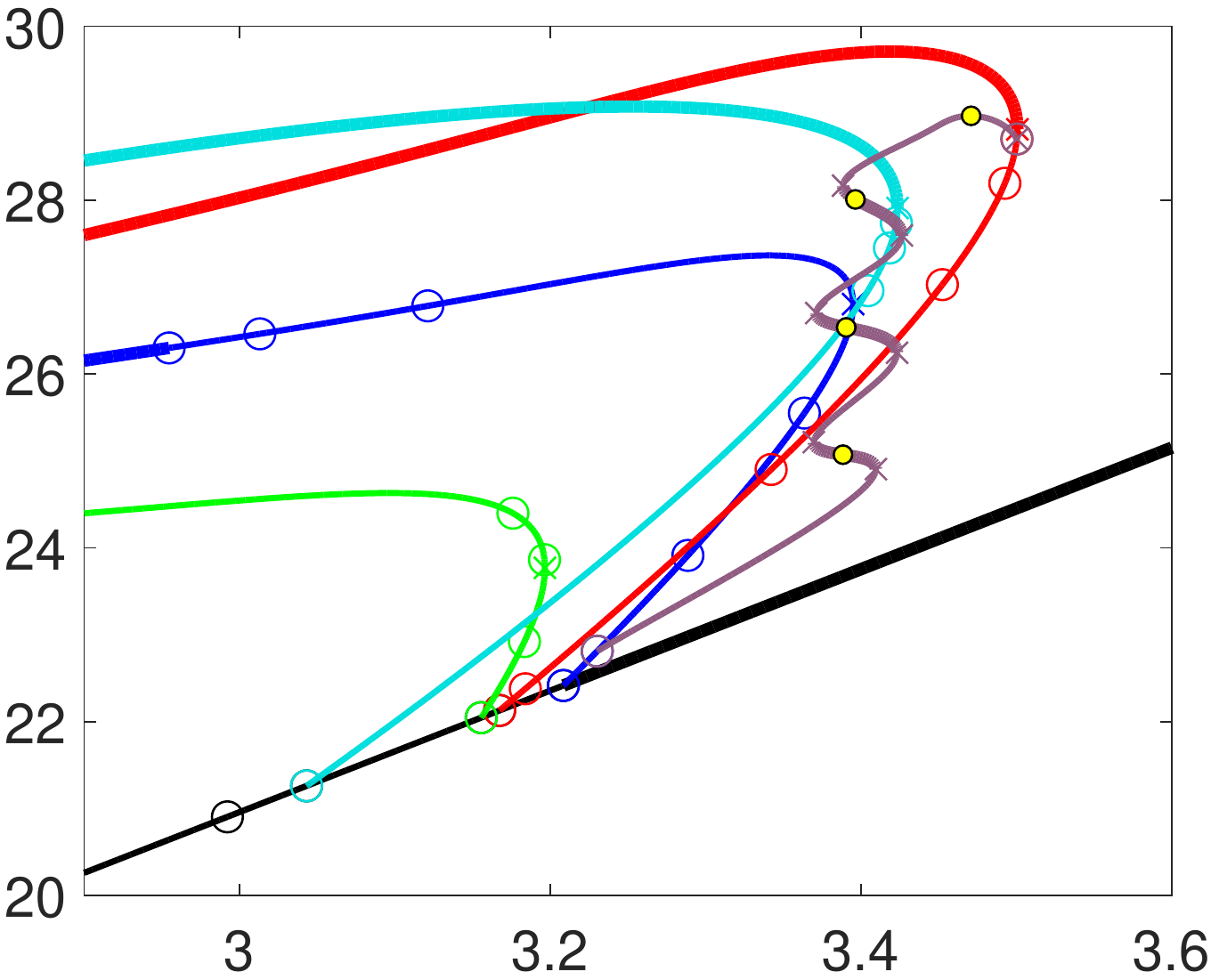}
    \put(85,-2){$\mu$}
    \put(-7,45){\rotatebox{90}{$\|u_1\|_{L^2}$}}
    \put(75,40){(a)}
    \put(59,50){(b)}
    \put(62,62){(c)}
    \put(71,71){(d)}
    \end{overpic}
\end{subfigure}
\newpage
\hspace{0.7cm}
\begin{subfigure}[b]{0.2\textwidth}
     \begin{overpic}[width=\textwidth]{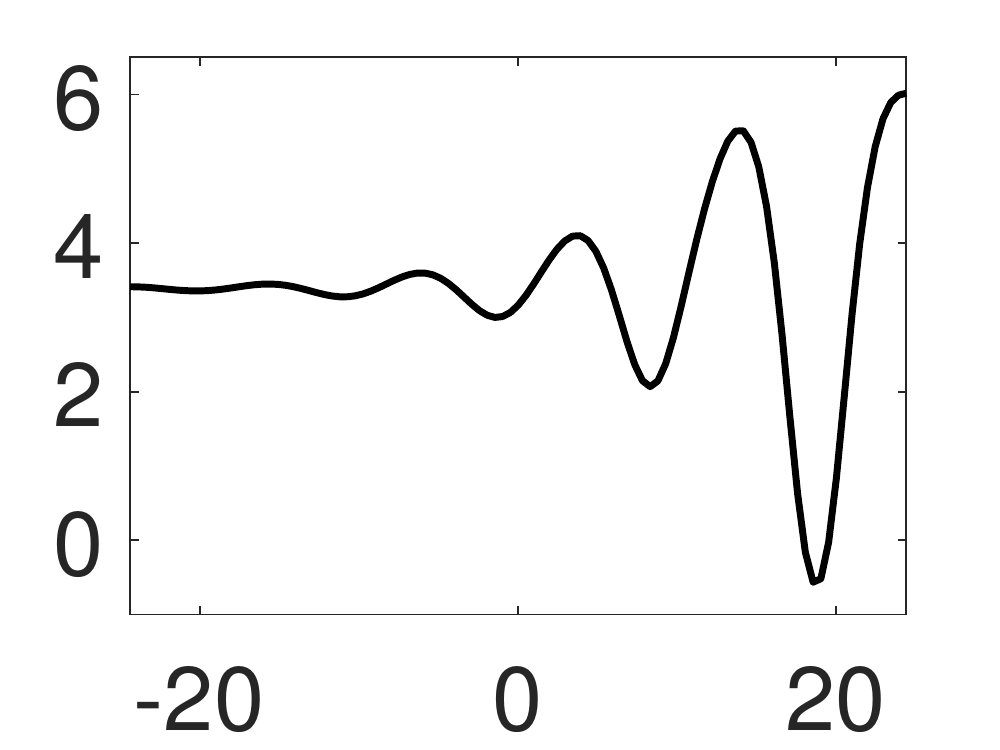}
     \put(-10,30){\rotatebox{90}{$u_1(x)$}}
     \put(100,5){$x$}
    \end{overpic}
     \caption{}
     \label{subfig-2:sn1-pt70}
\end{subfigure}
\hspace{0.2cm}
\begin{subfigure}[b]{0.2\textwidth}
     \begin{overpic}[width=\textwidth]{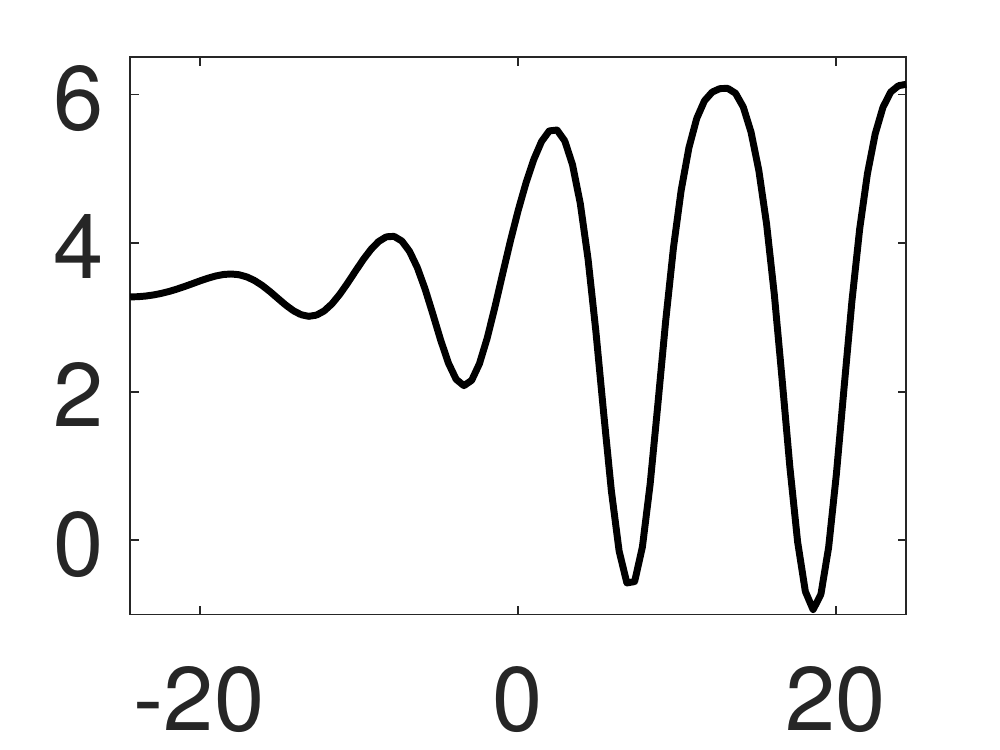}
     \put(-10,30){\rotatebox{90}{$u_1(x)$}}
     \put(100,5){$x$}
    \end{overpic}
     \caption{}
     \label{subfig-3:sn1-pt190}
\end{subfigure}\\

\hspace{0.7cm}
\begin{subfigure}[b]{0.2\textwidth}
     \begin{overpic}[width=\textwidth]{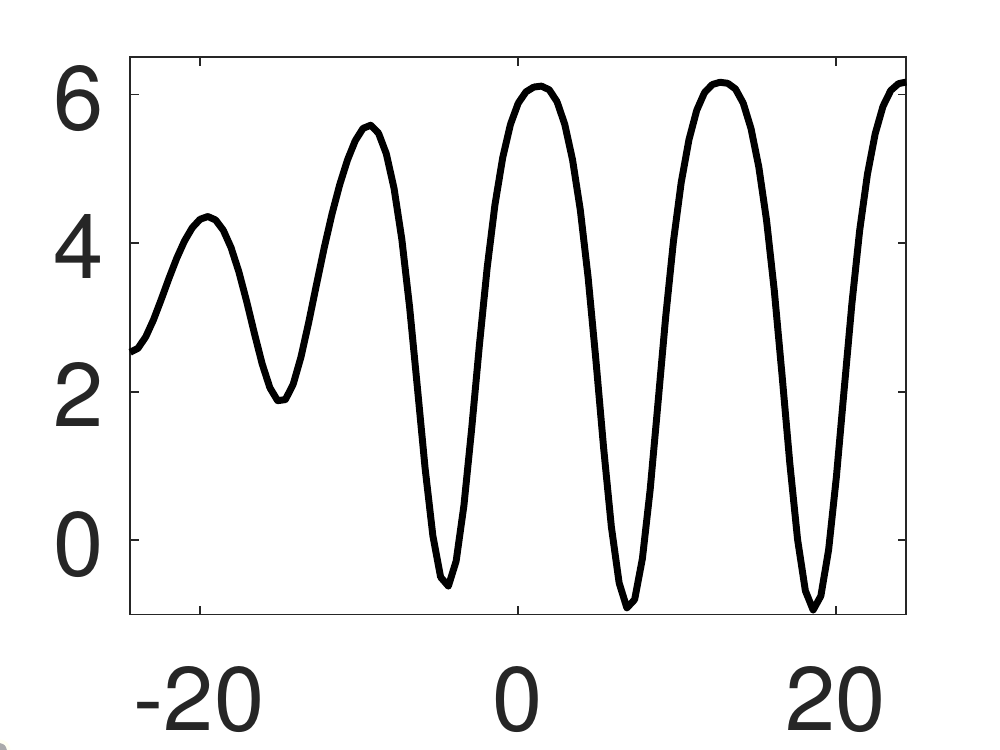}
     \put(-10,30){\rotatebox{90}{$u_1(x)$}}
     \put(100,5){$x$}
    \end{overpic}
     \caption{}
     \label{subfig-3:sn1-pt270}
\end{subfigure}
\hspace{0.2cm}
\begin{subfigure}[b]{0.2\textwidth}
     \begin{overpic}[width=\textwidth]{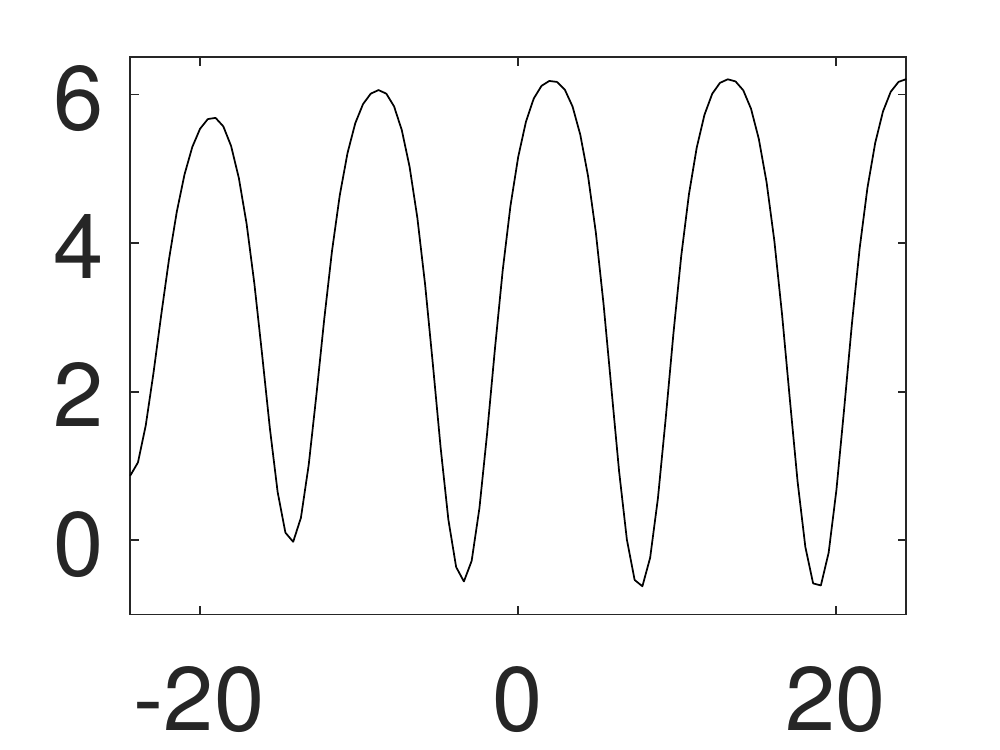}
     \put(-10,30){\rotatebox{90}{$u_1(x)$}}
     \put(100,5){$x$}
    \end{overpic}
     \caption{}
     \label{subfig-4:sn1-pt340}
\end{subfigure}
\end{multicols}
\vspace{-0.5cm}
\caption{Steady state bifurcation diagram of the Schnakenberg system~\eqref{eq:Schnakenberg} with $d=60$ and modified reaction terms~\eqref{eq:F_m} with~$\sigma=-0.6$, on a domain with length~$L=10\pi/k_c$ with homogeneous Neumann boundary conditions. Black: homogeneous states. Blue, red, green, light blue: branches of non-homogeneous stationary solutions presenting 5, 4.5, 5.5 and 4 bumps, respectively. Purple: snaking branch of front solutions shown in~Figures~\ref{subfig-2:sn1-pt70}~--~\ref{subfig-4:sn1-pt340}, connecting the blue branch to the red one (thick and thin lines denote stable and unstable solutions respectively, while circles and crosses indicate branch and fold points).  (For interpretation of the references to color in this figure legend, the reader is referred to the web version of this article)}
\label{fig:Schnak-sig06-theory}
\end{figure}
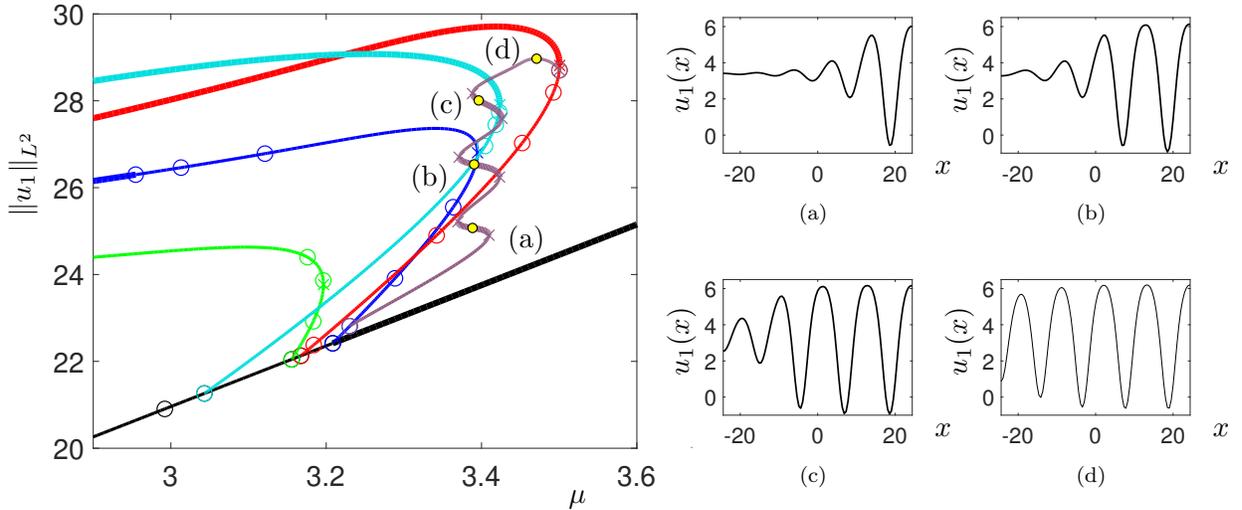

 Note that both ~$\mu_c$ and~$k_c$ are independent of $\sigma$. This is because $\sigma$ does not affect the linearization. However, it affects the nonlinear terms in~\eqref{eq:Schnakenberg} and hence the restabilization of modes at higher amplitudes \cite{Uecker2014}. Thus, we can use $\sigma$ to tune the primary bifurcation from super- to subcritical. This enables us to ``create'' snaking branches between periodic branches, as is illustrated in Figure~\ref{fig:Schnak-sig06-theory}. As for the Swift--Hohenberg equation, the snaking branch corresponds to front solutions illustrated in Figure~\ref{subfig-2:sn1-pt70}--\ref{subfig-4:sn1-pt340}.

\subsection{The fractional Schnakenberg system}
\label{chapter:Schnakenberg-system-results}

We finally study the fractional Schnakenberg system~\eqref{eq:Schnak-frac}. The problem is considered on a one dimensional domain of length $L$ with homogeneous Neumann boundary conditions. We study both classical ($\sigma=0$) and modified ($\sigma\neq 0$) reaction terms~\eqref{eq:F_m} as presented in the Introduction.

As for the Schnakenberg system with standard Laplacian~\eqref{eq:Schnakenberg}, the spatially homogeneous state $\bar{u}=(\mu, {1}/{\mu})$ is a stationary solution to~\eqref{eq:Schnak-frac}. Moreover, according to \eqref{eq:spectral_definition} the Neumann spectral fractional Laplacian has eigenpairs $(\phi_j, \lambda_j) = (\cos\left( {j\pi x}/{L} \right), -\left({j \pi}/{L} \right)^{2s})$ for~$j\geq 0$. The bifurcation points on the homogeneous branch then turn out to be located at
\begin{equation}\label{eq:bifurcations-schnak-frac}
\mu_j = k_j^s\sqrt{\frac{d(1-k_j^{2s})}{1+k_j^{2s}}}, \quad  k_j = \frac{j\pi}{L}, \quad j\geq 0.
\end{equation}
Further, the critical parameter $\mu_c$ at which the first instability occurs is not affected by the fractional diffusion, being
\begin{equation}\label{eq:mu_schnak}
\mu_c = \sqrt{d(3-\sqrt{8})}.
\end{equation}
The critical wavenumber, however, changes as a function of the fractional order $s$
\begin{equation}\label{eq:kc-schnak-frac}
k_c = \left(\frac{d-\mu_c^{2}}{2d}\right)^{\frac{1}{2s}}=\left(\sqrt{2}-1\right)^{\frac{1}{2s}}.
\end{equation}
Thus, we expect the wavelength of non-homogeneous stationary solutions from the first bifurcation to increase as $s$ decreases.
In the following, we fix the diffusion coefficient $d=60$. Thus, equation~\eqref{eq:mu_schnak} leads to the critical parameter $ \mu_c \approx 3.21$. Further, we choose the domain
\begin{equation}\label{eq:domain-schnak}
\Omega= \left(-\frac{m\pi}{k_c}, \frac{m\pi}{k_c}\right),\quad m\in\mathbb{N}_{>0},
\end{equation}
i.e. the domain is ``tuned'' to the primary bifurcation for each fractional order $s$. Precisely, according to equation~\eqref{eq:kc-schnak-frac}  we know that $k_c$ decreases with $s$, thus the size of $\Omega$, $L={2m\pi}/{k_c}$, is increased with $s$ in such a way that the first unstable mode always has $m$ periods in the domain. This choice guarantees that the primary bifurcation always occur at $\mu_c$ to modes with known wavenumber given by formula~\eqref{eq:kc-schnak-frac}.
Note that for the Allen--Cahn equation and for the Swift--Hohenberg equation, the domain size was always fixed since in these systems $k_c$ did not depend on $s$.
Furthermore, because of this choice and to be able to compare results for different fractional orders, it is important to use a measure for the bifurcation diagram that does not depend on the domain size. The normalized $L^2$-norm which was used in previous sections satisfies this criteria. In addition, we also use the normalized $L^8$-norm~\cite{Uecker2014}\\[-0.4cm]

$$
\quad \quad \|u\|_{L^8} = \left( \frac{1}{|\Omega|}\int_{\Omega} |u(x)|^8 \txtd x \right) ^{\frac{1}{8}},
$$

\noindent
since in some situations it has been shown to be more suitable for plotting than the standard normalized $L^2$-norm as it produces larger differences between states.

Finally, note that since the domain size increases as $s$ decreases, we must accordingly increase the number of discretization points to guarantee the accuracy of the computation (the values used in the following are listed in Table \ref{tab:Schnak-nph}). This induces a significant increase in computational cost. For this reason we studied the problem with fractional order not smaller than $s=0.5$.

\begin{table}
\centering
\begin{tabular}{cccc}
\multicolumn{4}{c}{$\sigma=0$}\\
\toprule
$s$ & $L$ & $n_p$ & $h$\\
\midrule
$0.9$ & $6.5\pi$ & $401$ & $0.05$\\
$0.7$ & $7.5\pi$ & $1501$ & $0.01$\\
$0.5$ & $9.7\pi$ & $2501$ & $0.009$\\
\bottomrule
\end{tabular}
\hspace{1cm}
\begin{tabular}{cccc}
\multicolumn{4}{c}{$\sigma=-0.6$}\\
\toprule
 $s$ & $L$ & $n_p$ & $h$\\
\midrule
$0.95$ & $15.9\pi$ & $1601$ & $0.0312$\\
$0.9$   & $16.4\pi$ & $1601$ & $0.0322$\\
$0.8$   & $17.2\pi$ & $1901$ & $0.028$\\
$0.73$ & $18.3\pi$ & $2601$ & $0.022$\\
\bottomrule
\end{tabular}
\caption{Values of the domain size $L$, the number of meshpoints $n_p$ and the corresponding meshsize~$h$ for different values of the fractional orders $s$, used for the fractional Schnakenberg system~\eqref{eq:Schnak-frac} with classical and modified reaction terms ($\sigma=0$ and $\sigma=-0.6$, respectively). }
\label{tab:Schnak-nph}
\end{table}

\subsubsection{Numerical results for the non-modified system}
In this section we consider the fractional Schnakenberg system~\eqref{eq:Schnak-frac} with classical reaction terms, i.e.~$\sigma~=~0$. We choose a relatively small domain, namely $m=2$ in~\eqref{eq:domain-schnak}, to limit the increase in computational cost as the fractional order $s$ gets smaller.

\begin{figure}
    \centering
    \begin{subfigure}[b]{0.45\textwidth}
        \centering
        \begin{overpic}[trim=0 0 0 0, clip, height=2.2in]{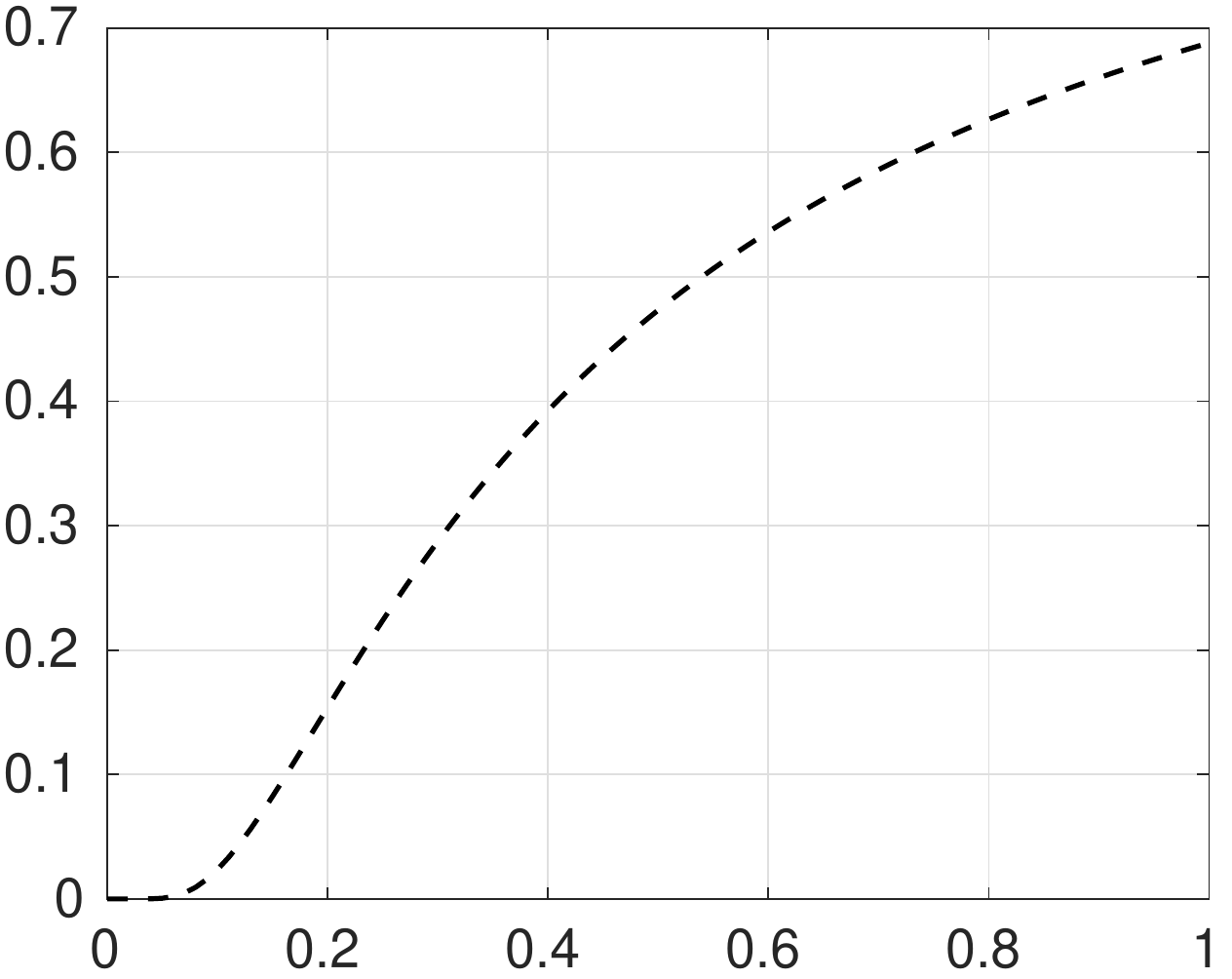}
        \put(-10,45){\rotatebox{90}{$k_c$}}
        \put(100,5){$s$}
        \end{overpic}
        \caption{}
        \label{fig:kc-schnak}
    \end{subfigure}
    \hspace{1cm}
    \begin{subfigure}[b]{0.45\textwidth}
        \centering
        \begin{overpic}[trim=0 1.15 0 0, clip, height=2.2in]{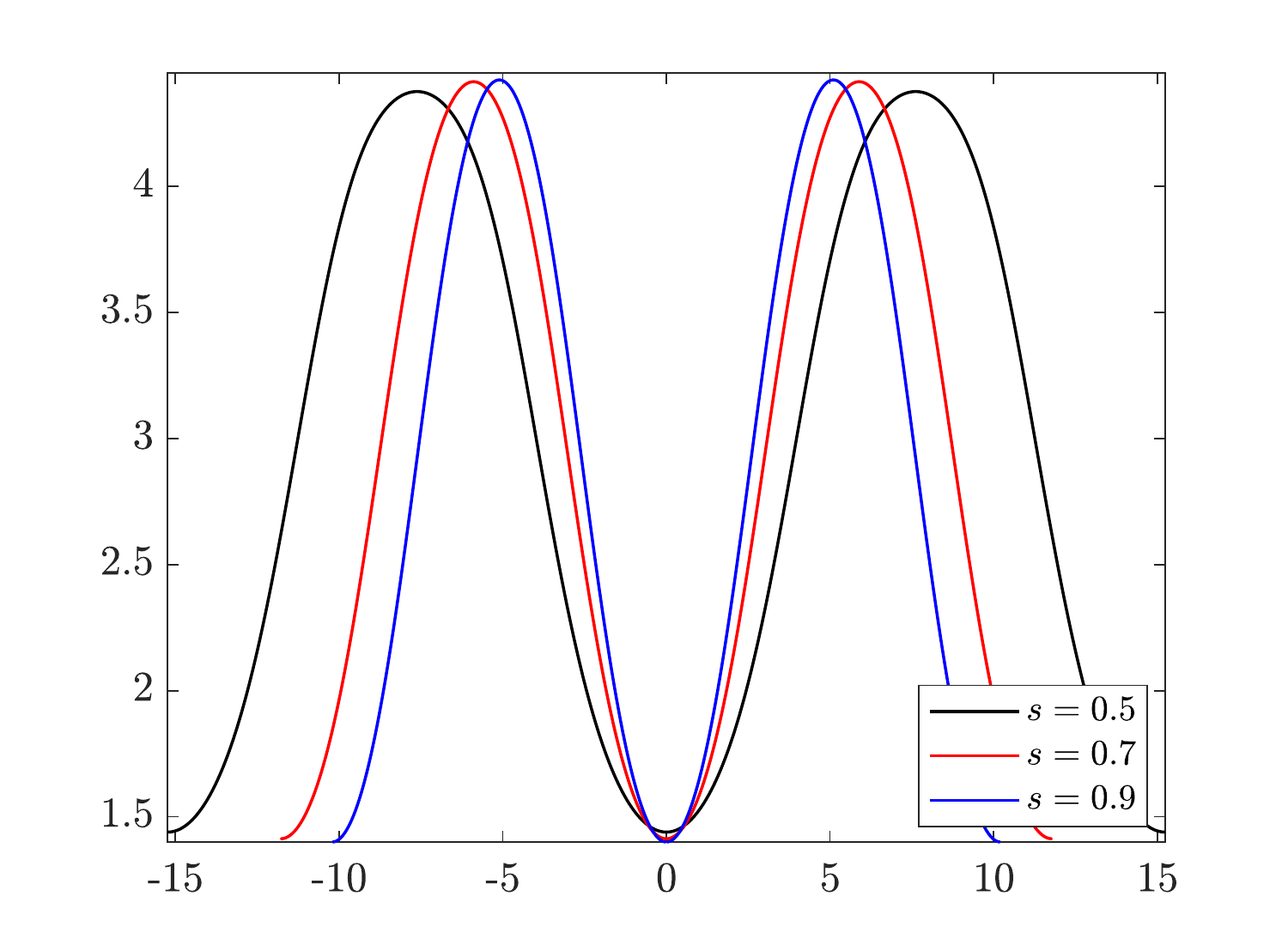}
        \put(-5,35){\rotatebox{90}{$u_1(x)$}}
         \put(100,5){$x$}
        \end{overpic}
        \caption{}
        \label{fig:waveelength-increase-schnak}
    \end{subfigure}
    \caption{Non-homogeneous stationary solutions of the Schnakenberg system~\eqref{eq:Schnak-frac} with classical reaction terms ($\sigma=0$) close to the first bifurcation from the homogeneous states. Left: decrease of the critical wavenumber $k_c$. Right: increase in the wavelength of the solutions for $\mu \approx \mu_c$ for different fractional orders $s$. Remember that to each fractional order corresponds a different domain size according to~\eqref{eq:domain-schnak}; in this figure we show only part of the solutions for a comparison.  (For interpretation of the references to color in this figure legend, the reader is referred to the web version of this article)}
    \label{fig:sols-1D1}
\end{figure}

Before showing the bifurcation diagram of system~\eqref{eq:Schnak-frac}, let us first illustrate the increase in the wavelength of solutions as $s$ decreases. The critical wavenumber, computed via equation~\eqref{eq:kc-schnak-frac}, for different values of $s$ is plotted in Figure~\ref{fig:kc-schnak}. Figure~\ref{fig:waveelength-increase-schnak} shows the corresponding change in wavelength for solutions close to the first bifurcation at $\mu \approx~3.2$. Further, one can see that the amplitude of the solution decreases with $s$. In fact, we notice that the $L^2$-norm of these solutions is conserved for all $s$ from $0.5$ and above (as mentioned, solutions below $s=0.5$ have not been computed).

\begin{figure}
\centering
\begin{subfigure}[b]{0.3\textwidth}
 \centering
 \begin{overpic}[width=\textwidth]{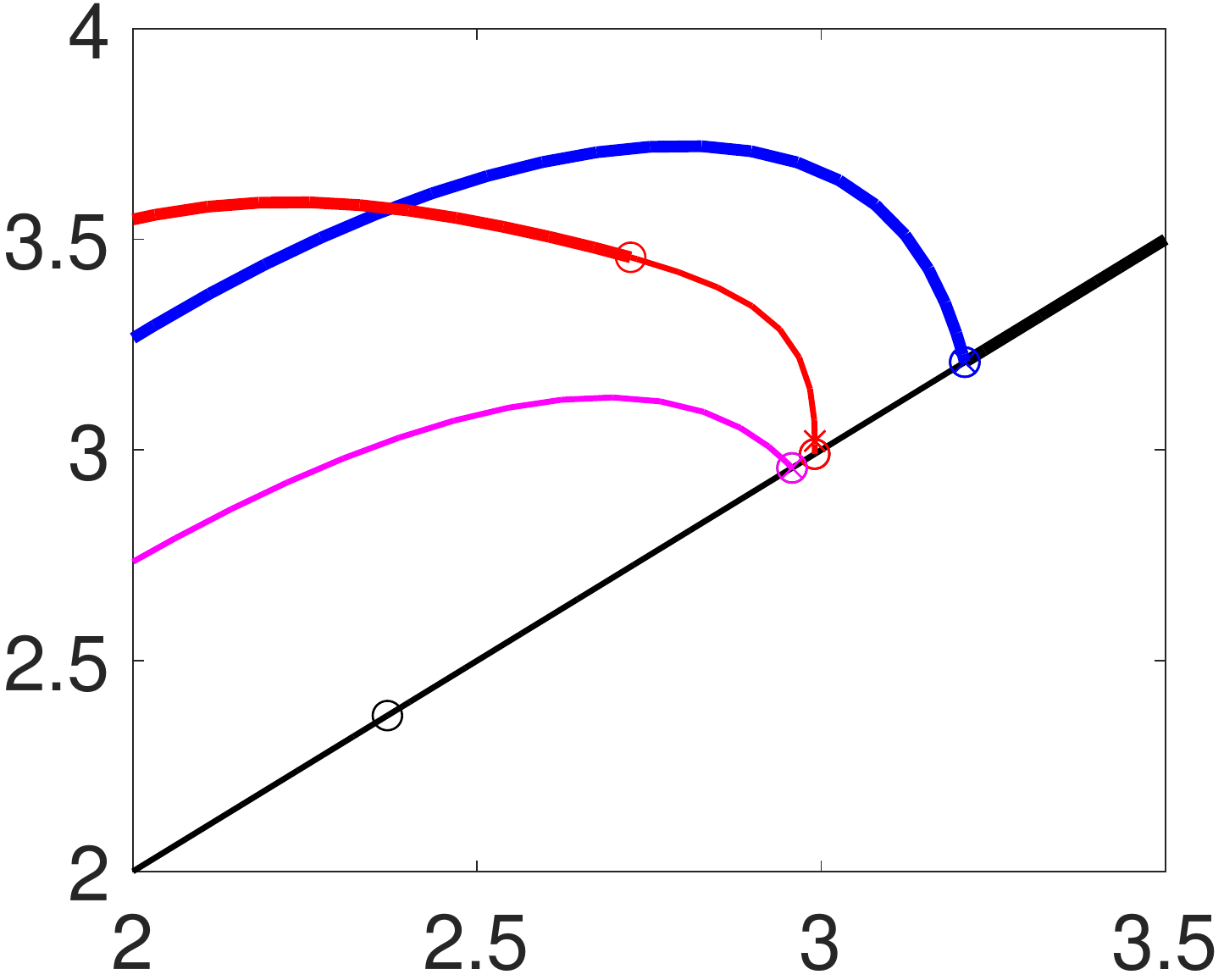}
 \put(-10,35){\rotatebox{90}{$\|u_1\|_{L^8}$}}
 \put(100,5){$\mu$}
 \end{overpic}
 \caption{$s=0.9$}\label{fig:bif-sig0-s09}
\end{subfigure}
\hspace{0.5cm}
\begin{subfigure}[b]{0.3\textwidth}
 \centering
 \begin{overpic}[width=\textwidth]{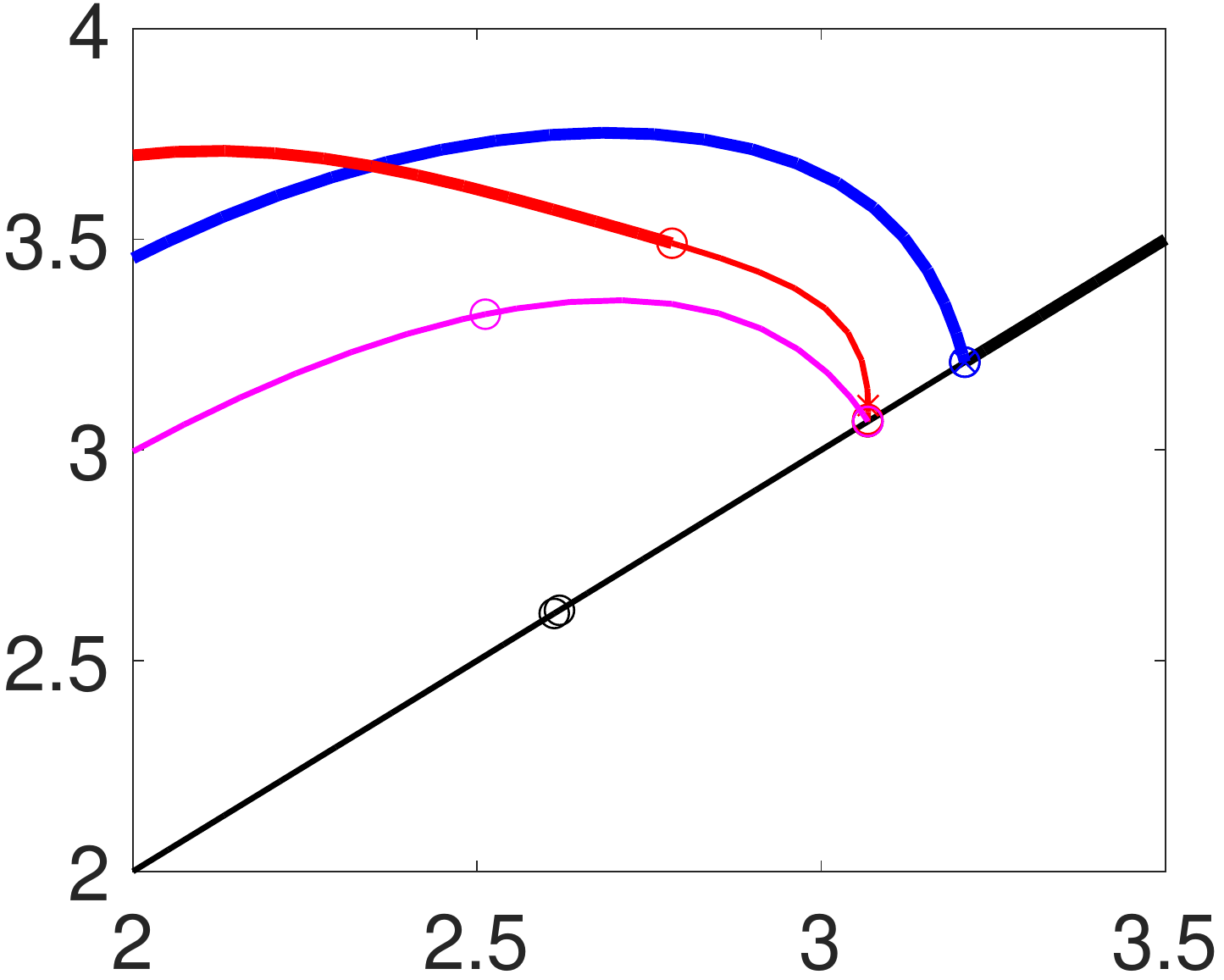}
\put(-10,35){\rotatebox{90}{$\|u_1\|_{L^8}$}}
 \put(100,5){$\mu$}
 \end{overpic}
 \caption{$s=0.7$}\label{fig:bif-sig0-s06}
\end{subfigure}
\hspace{0.5cm}
\begin{subfigure}[b]{0.3\textwidth}
 \centering
 \begin{overpic}[width=\textwidth]{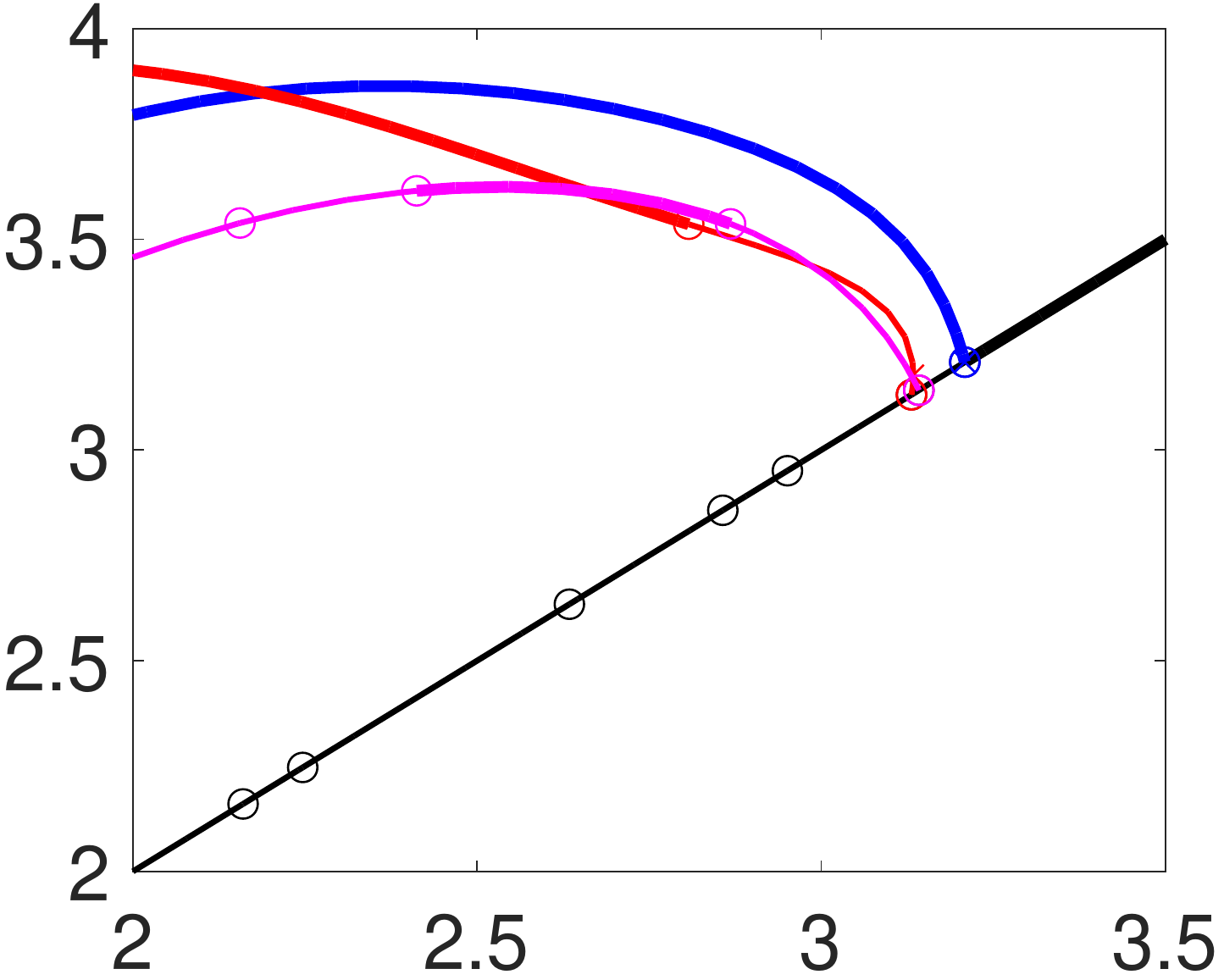}
\put(-10,35){\rotatebox{90}{$\|u_1\|_{L^8}$}}
 \put(100,5){$\mu$}
 \end{overpic}
 \caption{$s=0.5$}\label{fig:bif-sig0-s03}
\end{subfigure}
\caption{Evolution of the structure of the bifurcation diagram for the steady state Schnakenberg system~\eqref{eq:Schnak-frac} with classical reaction terms ($\sigma=0$) as the order of the fractional Laplacian decreases. Thick and thin lines denote stable and unstable solutions respectively, while circles indicate branch points. (For interpretation of the references to color in this figure legend, the reader is referred to the web version of this article)}
\label{fig:bif-sig0}
\end{figure}

\begin{figure}
\centering
\begin{overpic}[width=0.5\textwidth]{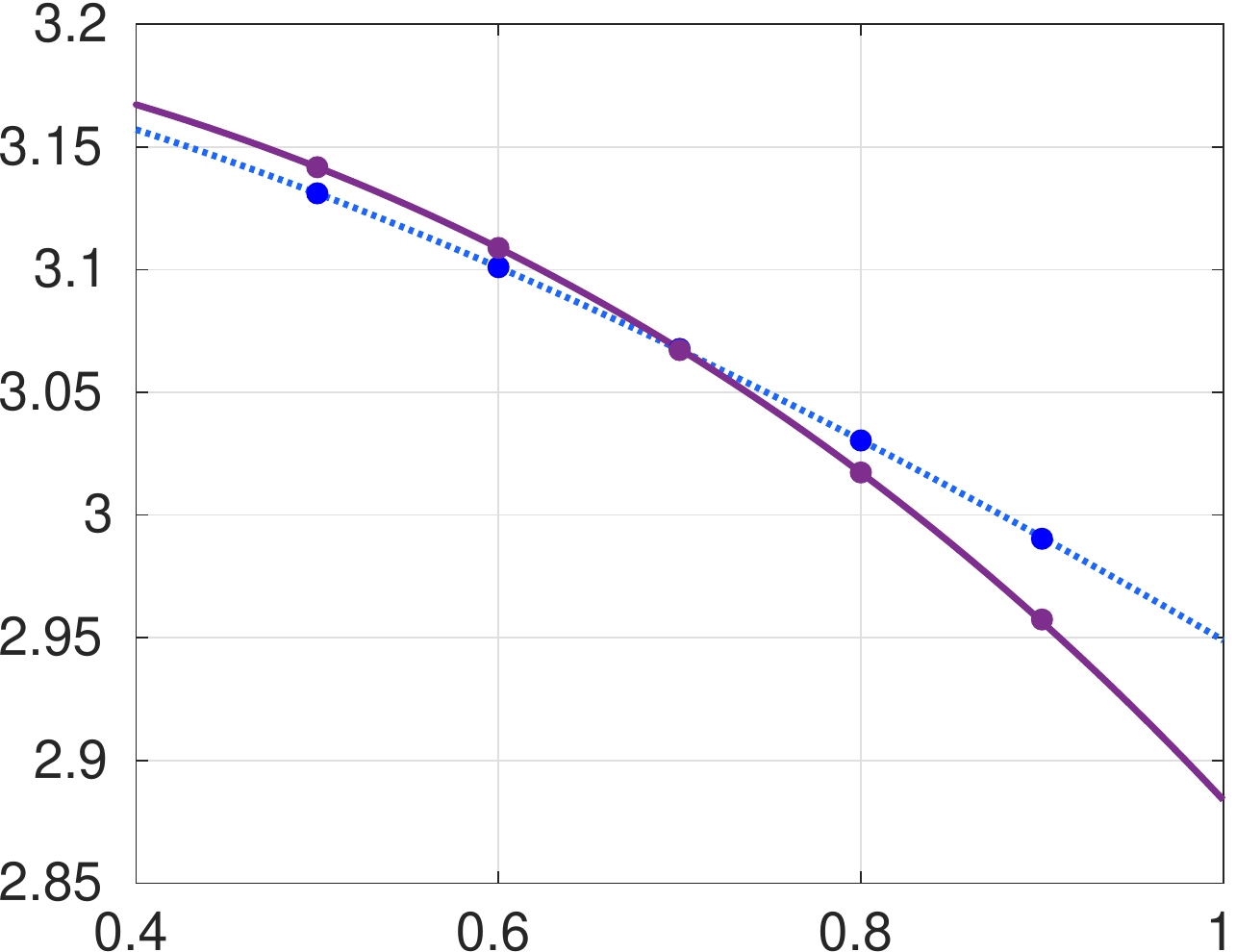}
\put(15,70){$\mu_3$}
\put(15,60){$\mu_2$}
 \put(100,5){$s$}
 \end{overpic}
\caption{Location of the second (blue dashed line) and third (purple continuous line) bifurcation points on the homogeneous branch of the Schnakenberg system~\eqref{eq:Schnak-frac} with classical reaction terms ($\sigma=0$) for a domain of size $L=4\pi/k_c$ as $s$ decreases. The dots correspond to numerical values and the lines are computed via equation~\eqref{eq:bifurcations-schnak-frac}. (For interpretation of the references to color in this figure legend, the reader is referred to the web version of this article)}
\label{fig:b1andb2goonawalk}
\end{figure}

In Figure~\ref{fig:bif-sig0} we report the bifurcation diagrams of system~\eqref{eq:Schnak-frac} for three different fractional orders (namely $s=0.9$, $s=0.7$ and $s=0.5$). We observe several changes as $s$ decreases. First, the bifurcation points on the homogeneous branch (in black) seem to move towards the first bifurcation at $\mu_c$. In addition, the second and third bifurcation points seem to exchange their position between $s=0.8$ and $s=0.7$. Correspondingly, the second (red) and third (magenta) bifurcating branches move closer and towards the first branch (blue) when $s$ becomes smaller. The numerical values of the bifurcation parameter $\mu$ for the location of the second and third bifurcation points, corresponding respectively to kernel vectors $(\cos\left( {3\pi x}/{L} \right), -\left({3 \pi}/{L} \right)^{2s})$ and $(\cos\left( {5\pi x}/{L} \right), -\left({5 \pi}/{L} \right)^{2s})$, are shown in Figure~\ref{fig:b1andb2goonawalk}. These values can easily be verified using formula~\eqref{eq:bifurcations-schnak-frac} for the domain size $L=4\pi/k_c$ and the corresponding curves are plotted as reference in Figure~\ref{fig:b1andb2goonawalk}. We observe that the numerical values agree with the analytical ones (with a maximum error of the order of $10^{-3}$). In addition, the bifurcation points come closer as the fractional order $s$ decreases and intersect close to $s=0.7$. Note that the bifurcation points eventually converge towards $\mu_c$ as $s\to0$.

\begin{figure}
 \centering
 \begin{subfigure}[b]{0.3\textwidth}
     \centering
     \begin{overpic}[width=0.8\textwidth,trim=20 20 0 0,clip]{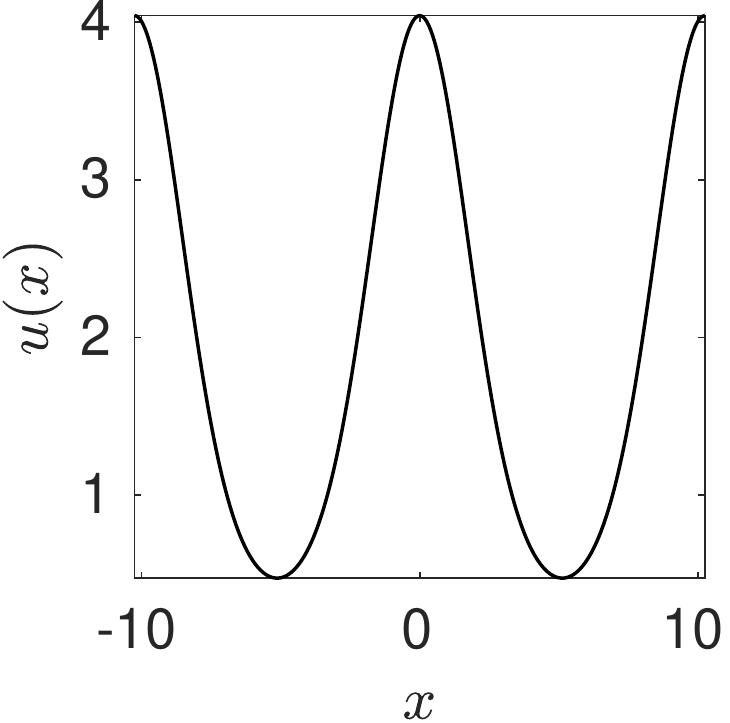}
     \put(-10,45){\rotatebox{90}{$u_1(x)$}}
     \put(100,5){$x$}
     \end{overpic}
     \caption{$s=0.9$}
     \label{fig:u-s09}
\end{subfigure}
\hspace{3cm}
\begin{subfigure}[b]{0.3\textwidth}
     \centering
    \begin{overpic}[width=0.8\textwidth,trim=20 21 0 0,clip]{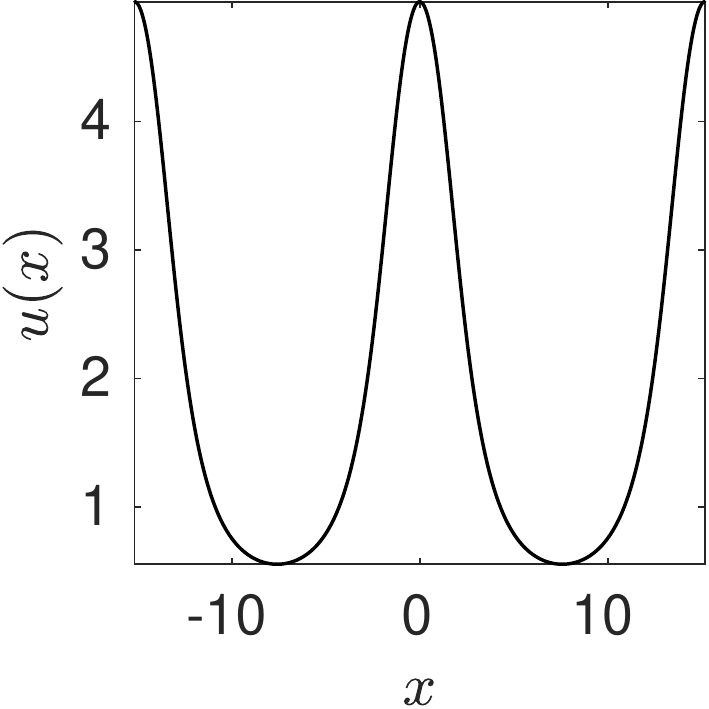}
    \put(-10,45){\rotatebox{90}{$u_1(x)$}}
     \put(100,5){$x$}
     \end{overpic}
     \caption{$s=0.5$}
     \label{fig:u-s03}
 \end{subfigure}
 \caption{Change in the solution profiles on the first bifurcating branch as the order of the fractional Laplacian decreases for the Schnakenberg system~\eqref{eq:Schnak-frac} with classical reaction terms ($\sigma=0$). }
 \label{fig:sol-bpt1-lam3}
\end{figure}

Second, we observe in Figure~\ref{fig:bif-sig0} that the $L^8$-norm of the branches of solutions is increasing. This is better understood by looking at the solution profiles. Figure~\ref{fig:sol-bpt1-lam3} shows the solution on the first bifurcating branch at $\mu = 3$ for $s=0.9$ and $s=0.5$. We see once more the wavelength increasing. In addition, for lower $s$ the solutions ``valleys'' tend to flatten. This is now a known effect, which we had already observed in the solutions of Allen--Cahn  and Swift--Hohenberg equations. In contrast to previous systems, however, the amplitude of the solution increases and the ``peak'' sharpens.

The last observation arising from Figure~\ref{fig:bif-sig0} regards the stability of solutions. As for the standard problem, with fractional diffusion the first branch bifurcating from the homogeneous states (blue) is entirely stable and the second one (red) becomes stable at some $\mu<\mu_c$. But, in contrast to the standard case, the third branch is unstable for $s=0.9$ and $s=0.7$ and a stable zone appears at $s=0.5$. This can be a consequence of the position exchange of the second and third bifurcation points on the homogeneous branch (see Figure~\ref{fig:b1andb2goonawalk}).

\subsubsection{Numerical results for the modified system}

We now turn to the fractional modified Schnakenberg system, i.e. equation~\eqref{eq:Schnak-frac} with $\sigma\neq 0$ in the reaction term. Recall from Section \ref{schnak} that in this case we expect to observe a snaking branch in the bifurcation diagram. Since in small domains snaking branches have very few turns, we choose here a domain larger than in the previous section in order to have more turns and facilitate our observations, namely $m=5$ in~\eqref{eq:domain-schnak}. Since we consider larger domains as the fractional order decreases, the computational costs increase, as explained in the previous section. Therefore we only consider fractional orders of the Laplacian greater than $s=0.7$.

\begin{figure}
\centering
\begin{overpic}[width=0.49\textwidth]{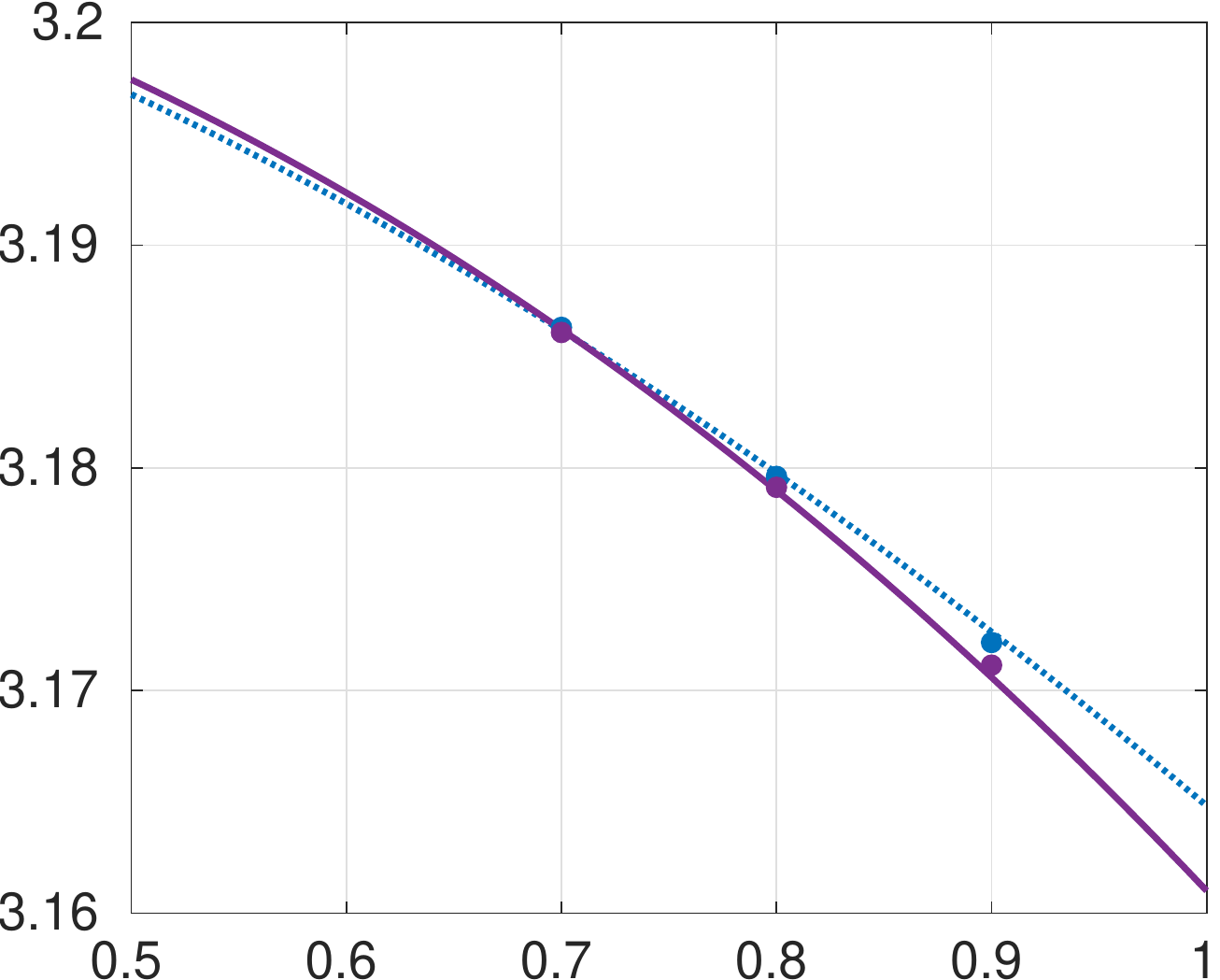}
\put(90,25){$\mu_2$}
\put(90,8){$\mu_3$}
\put(100,5){$s$}
\end{overpic}
\caption{Location of the second and third bifurcation points on the homogeneous branch of~\eqref{eq:Schnak-frac} for a domain of size $L=10\pi/k_c$ as $s$ decreases  for the fractional Schnakenberg system~\eqref{eq:Schnak-frac} with modified reaction terms ($\sigma=-0.6$).  (For interpretation of the references to color in this figure legend, the reader is referred to the web version of this article)}
\label{fig:b1b2-sig06}
\end{figure}

The numerical and analytical values for the location of the second and third bifurcation points (with kernel vectors $(\cos\left( {9\pi x}/{L} \right), -\left({9 \pi}/{L} \right)^{2s})$ and $(\cos\left( {11\pi x}/{L} \right), -\left({11 \pi}/{L} \right)^{2s})$) for fractional order~$s$ between $0.7$ and $0.9$ are shown in Figure~\ref{fig:b1b2-sig06}. The analytical values are computed using formula~\eqref{eq:bifurcations-schnak-frac} with $m=5$ and and agree with the numerical ones up to $3$ decimal digits. However this accuracy could not be maintained for fractional orders below $s=0.7$. Note that, from the analytical curves, we see that the bifurcation points cross each other close to $s=0.7$.

\begin{figure}
\centering
\begin{subfigure}[b]{0.45\textwidth}
\centering
\begin{overpic}[width=\textwidth]{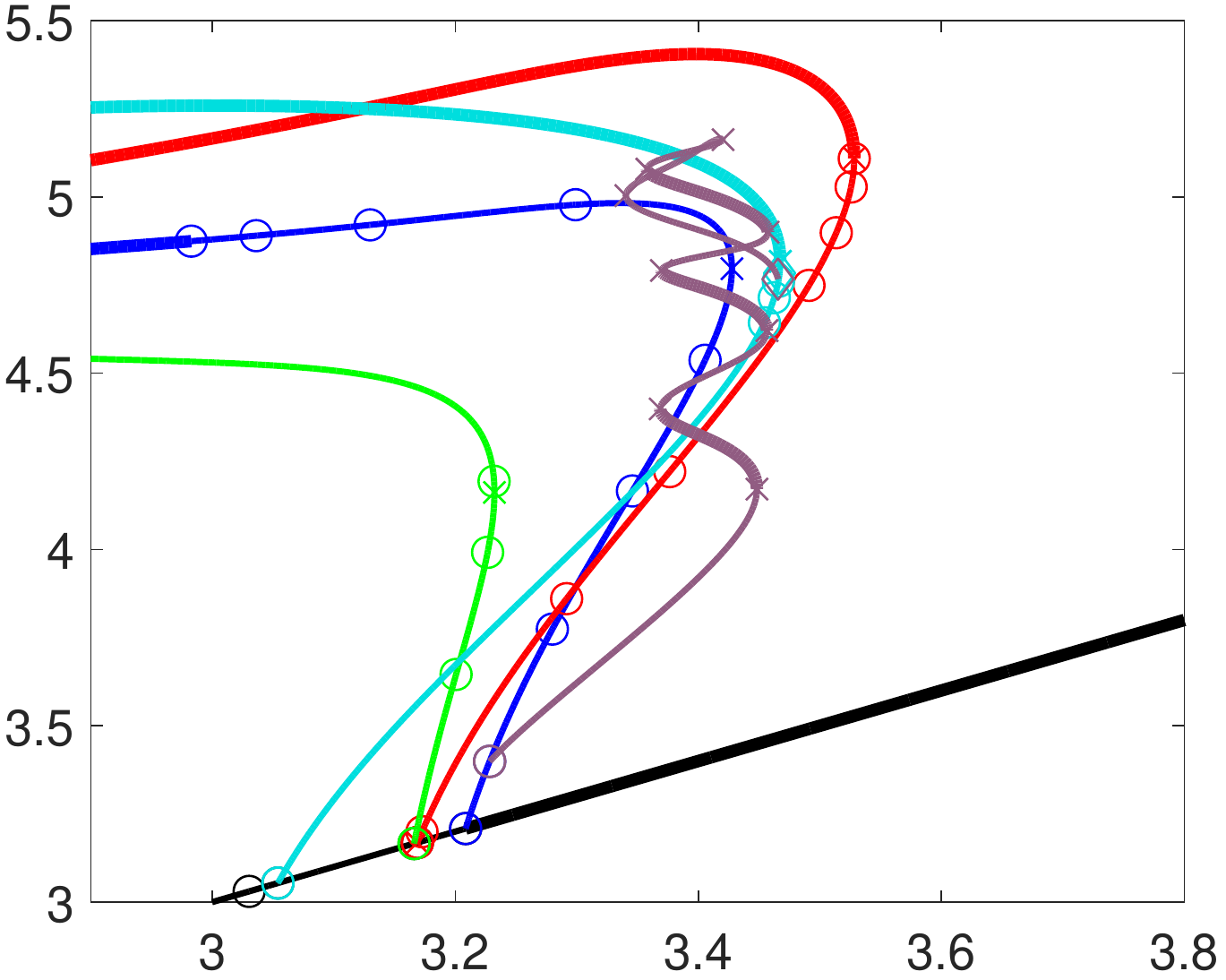}
\put(-5,40){\rotatebox{90}{$\|u_1\|_{L^8}$}}
\put(100,5){$\mu$}
\end{overpic}
\caption{$s=0.95$}
\label{fig:bif-sig06-s095}
\end{subfigure}
\hspace{1cm}
\begin{subfigure}[b]{0.45\textwidth}
\centering
\begin{overpic}[width=\textwidth]{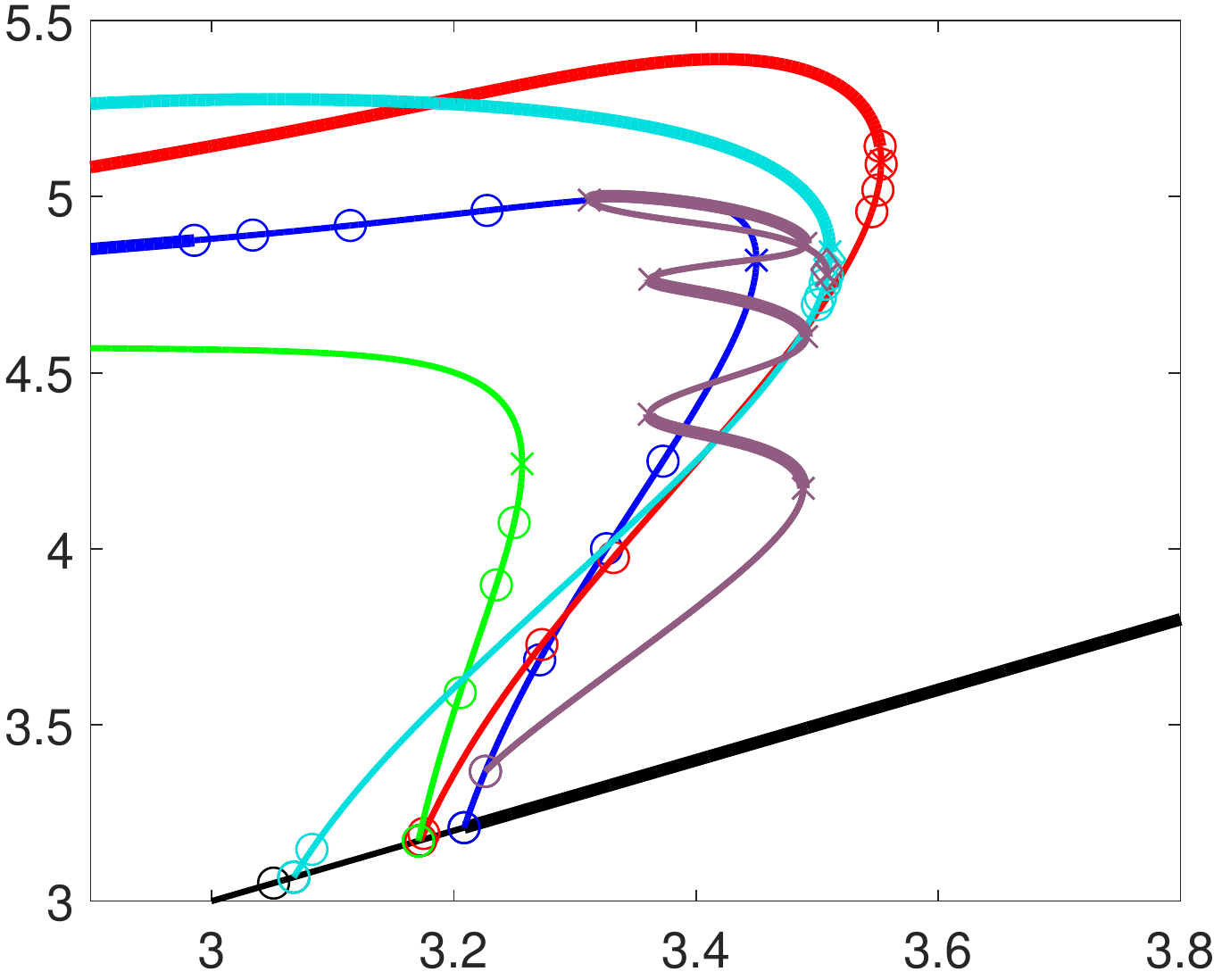}
\put(-5,40){\rotatebox{90}{$\|u_1\|_{L^8}$}}
\put(100,5){$\mu$}
\end{overpic}
\caption{$s=0.9$}
\label{fig:bif-sig06-s09}
\end{subfigure}
\hspace{1cm}
\begin{subfigure}[b]{0.45\textwidth}
\centering
\begin{overpic}[width=\textwidth]{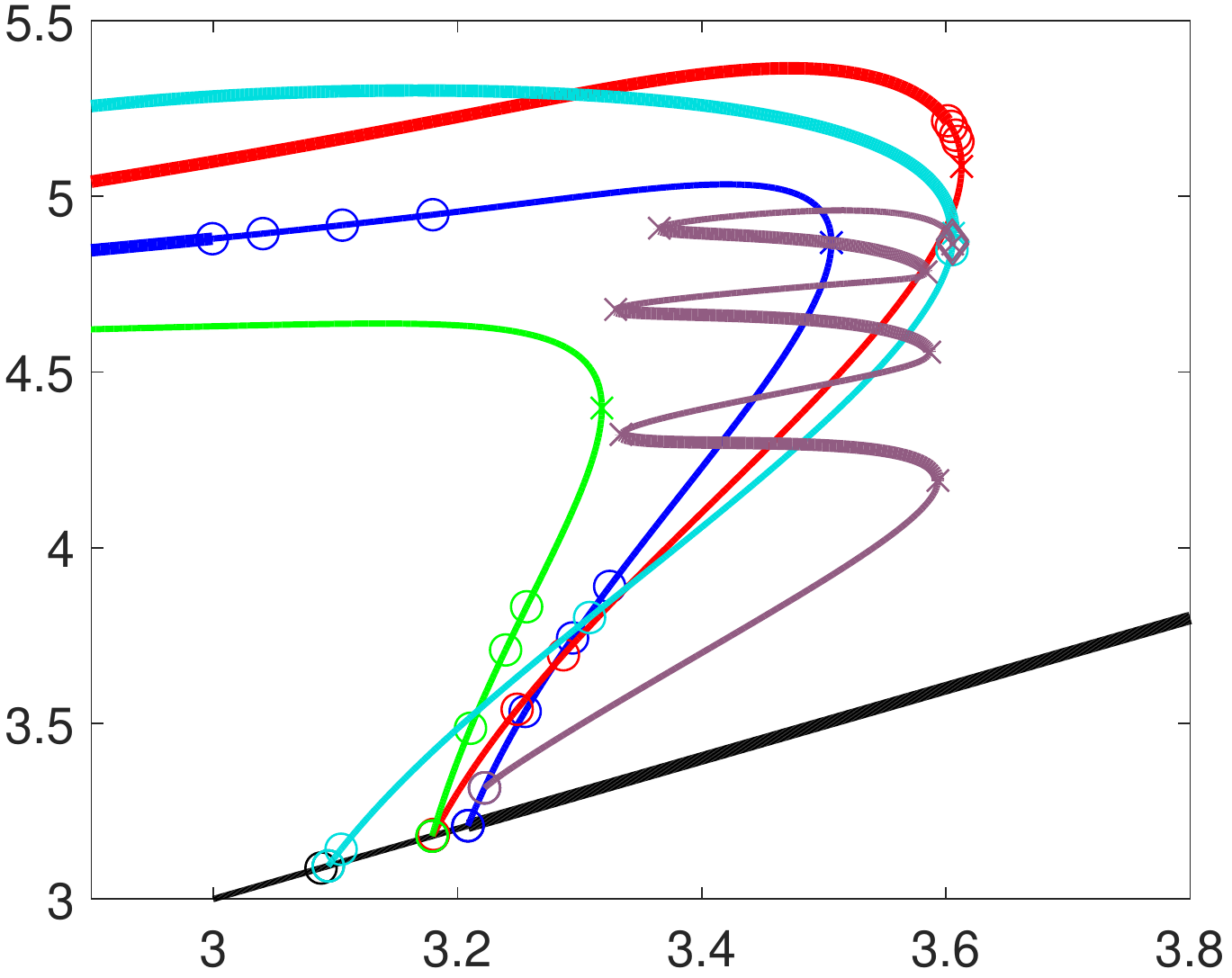}
\put(-5,40){\rotatebox{90}{$\|u_1\|_{L^8}$}}
\put(100,5){$\mu$}
\end{overpic}
\caption{$s=0.8$}
\label{fig:bif-sig06-s08}
\end{subfigure}
\hspace{1cm}
\begin{subfigure}[b]{0.45\textwidth}
\centering
\begin{overpic}[width=\textwidth]{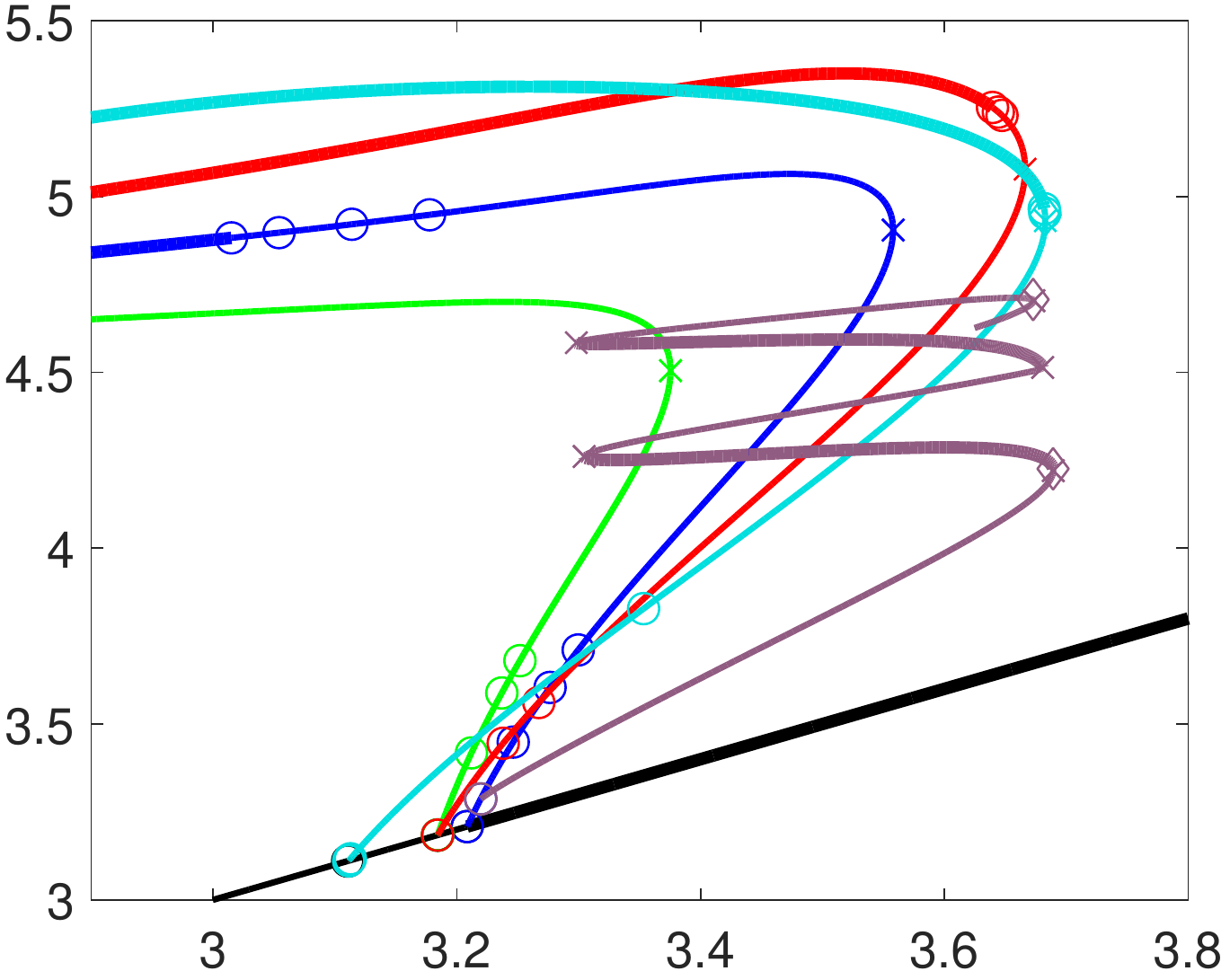}
\put(-5,40){\rotatebox{90}{$\|u_1\|_{L^8}$}}
\put(100,5){$\mu$}
\end{overpic}
\caption{$s=0.73$}
\label{fig:bif-sig06-s073}
\end{subfigure}
\caption{Evolution of the structure of the bifurcation diagram of the fractional Schnakenberg system~\eqref{eq:Schnak-frac} with modified reaction terms ($\sigma=-0.6$) as the order of the fractional Laplacian decreases. Thick and thin lines denote stable and unstable solutions, while circles, crosses and diamonds indicate branch, fold and Hopf points respectively.  (For interpretation of the references to color in this figure legend, the reader is referred to the web version of this article)}
\label{fig:bif-sig06}
\end{figure}
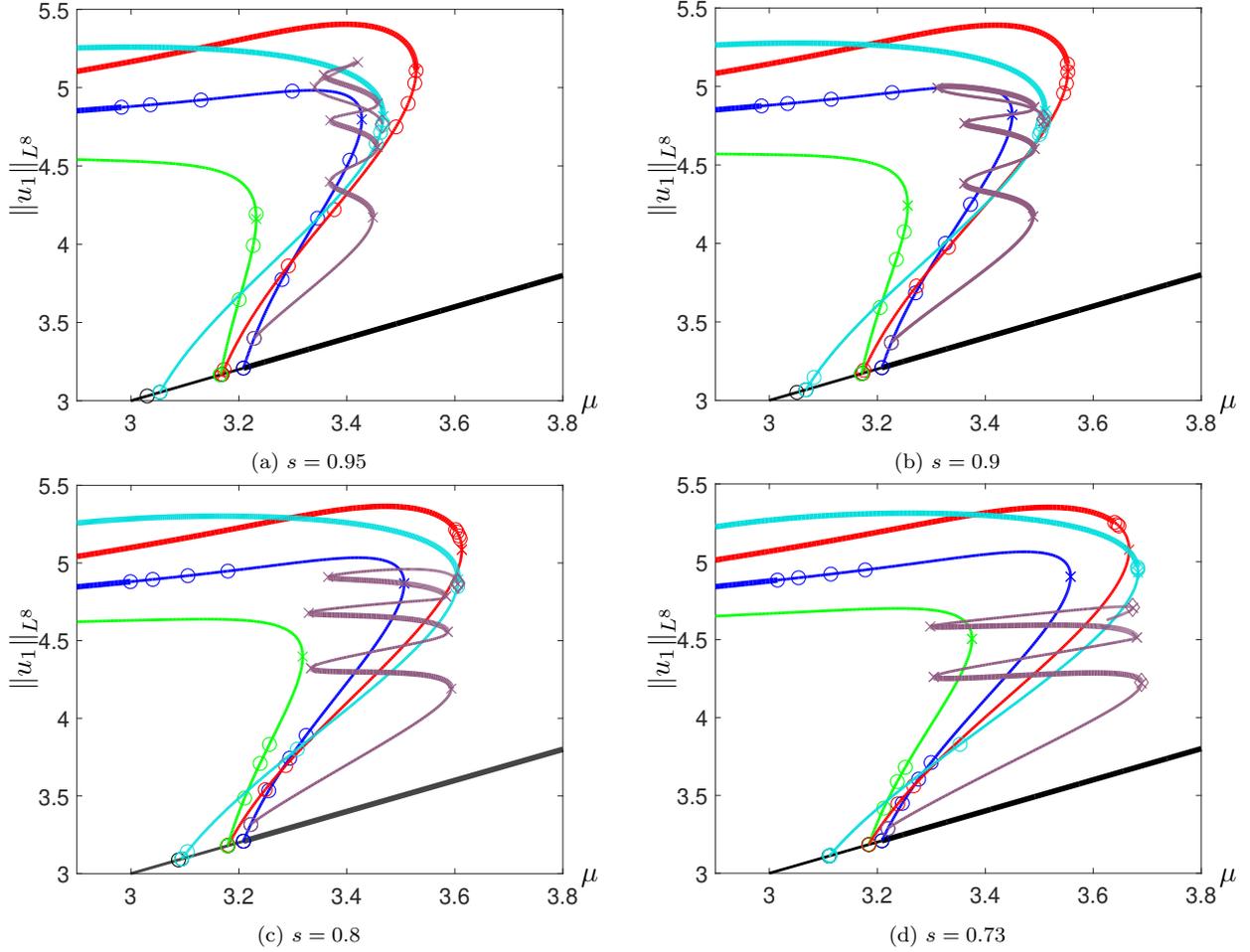

In order to get a first impression of the effects of the fractional Laplacian on the bifurcation structure of system~\eqref{eq:Schnak-frac} with modified reaction terms ($\sigma\neq 0$), Figure~\ref{fig:bif-sig06} shows the bifurcation diagram for four different fractional orders, namely $s=0.95$, $s=0.9$, $s=0.8$ and $s=0.73$. In the bifurcation diagrams, the blue branch, called $P5$, corresponds to non-homogeneous stationary solutions with five oscillations in the domain. The light blue branch, called $P4$, corresponds to periodic solutions with four oscillations in the domain. Analogously, the red ($P4.5$) and the green ($P5.5$) branches correspond to periodic solutions with 4.5 and 5.5 oscillations in the domain, respectively. Further, the snaking branch generated at the first bifurcation point on $P5$ is shown in purple. We observe three important changes in the bifurcation structure as $s$ decreases. First, as expected from the previous section, the bifurcation points are accumulating towards the critical point $\mu_c$ as $s$ gets smaller. Second, looking closely at the bifurcation diagrams, we see that the stable solutions on the blue branch shift towards the fold, that is towards higher $\mu$ values, as $s$ decreases. Recall that this effect was also observed in the fractional Allen--Cahn system (Section \ref{chapter:AC-results}).
Third, we observe that the snaking is widening as $s$ becomes smaller. In addition, at $s=0.73$, it neither reconnect to $P4$ nor ``climb'' to other branches. Moreover, the fractional order appears to have a significant influence on the behavior of the snaking branch. As mentioned above, we have seen in Figure~\ref{fig:bif-sig06}, that the number of turns and height of the snaking branch decreases progressively, and it does not ``climb'' anymore for $s=0.73$ (Figure~\ref{fig:bif-sig06-s073}). In addition, we can observe an increase in the width of the snaking as $s$ gets smaller.

\begin{figure}
\centering
\begin{multicols}{2}
\vspace*{\fill}
\begin{subfigure}{0.5\textwidth}
\centering
\begin{overpic}[width=\textwidth]{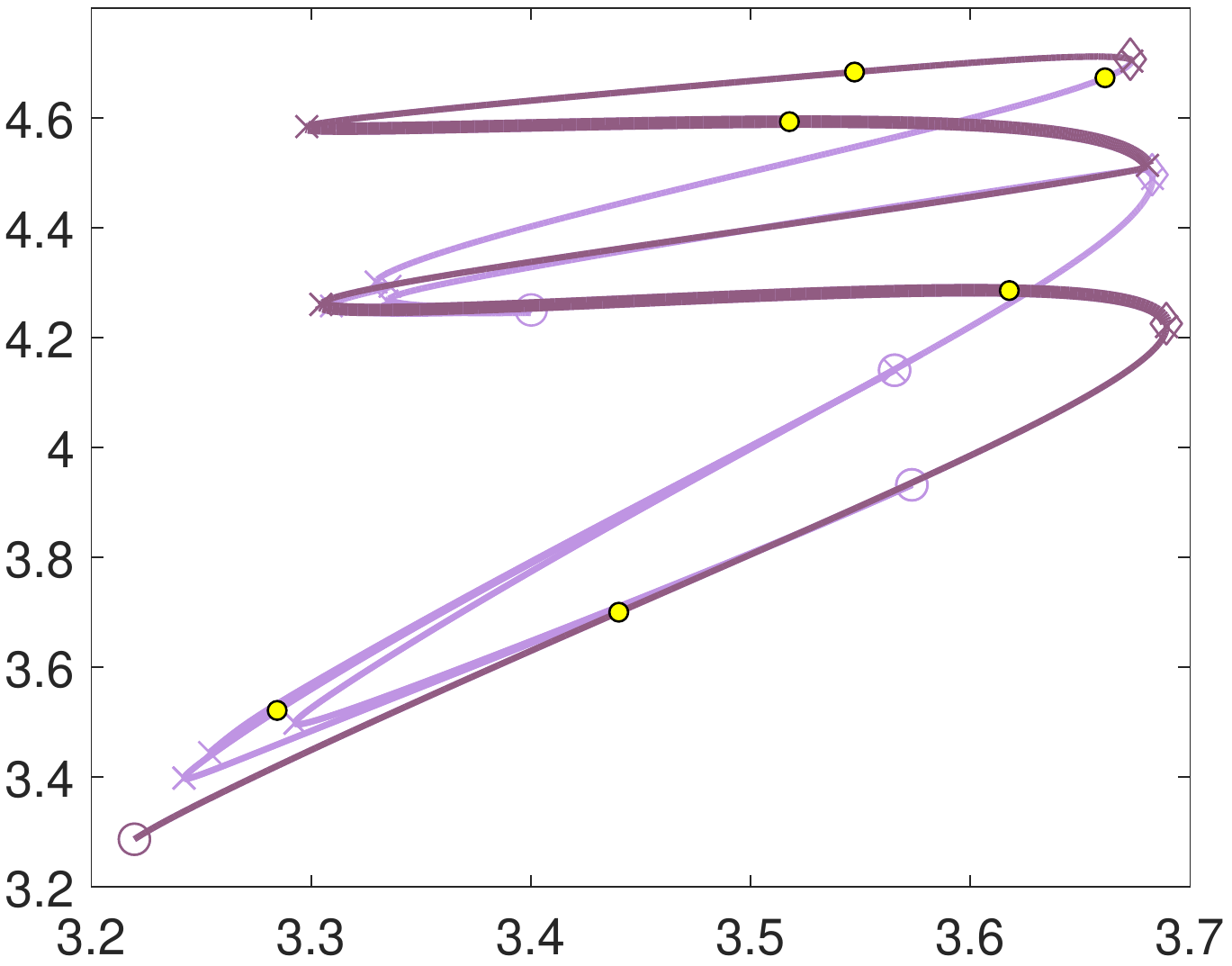}
\put(-5,40){\rotatebox{90}{$\|u_1\|_{L^8}$}}
\put(90,0){$\mu$}
\put(48,22){(a)}
\put(80,48){(b)}
\put(60,62){(c)}
\put(65,74){(d)}
\put(90,68){(e)}
\put(20,25){(f)}
\end{overpic}
\label{}
\end{subfigure}
 \vspace*{\fill}
\newpage
\begin{subfigure}{0.23\textwidth}
\centering
\begin{overpic}[width=\textwidth]{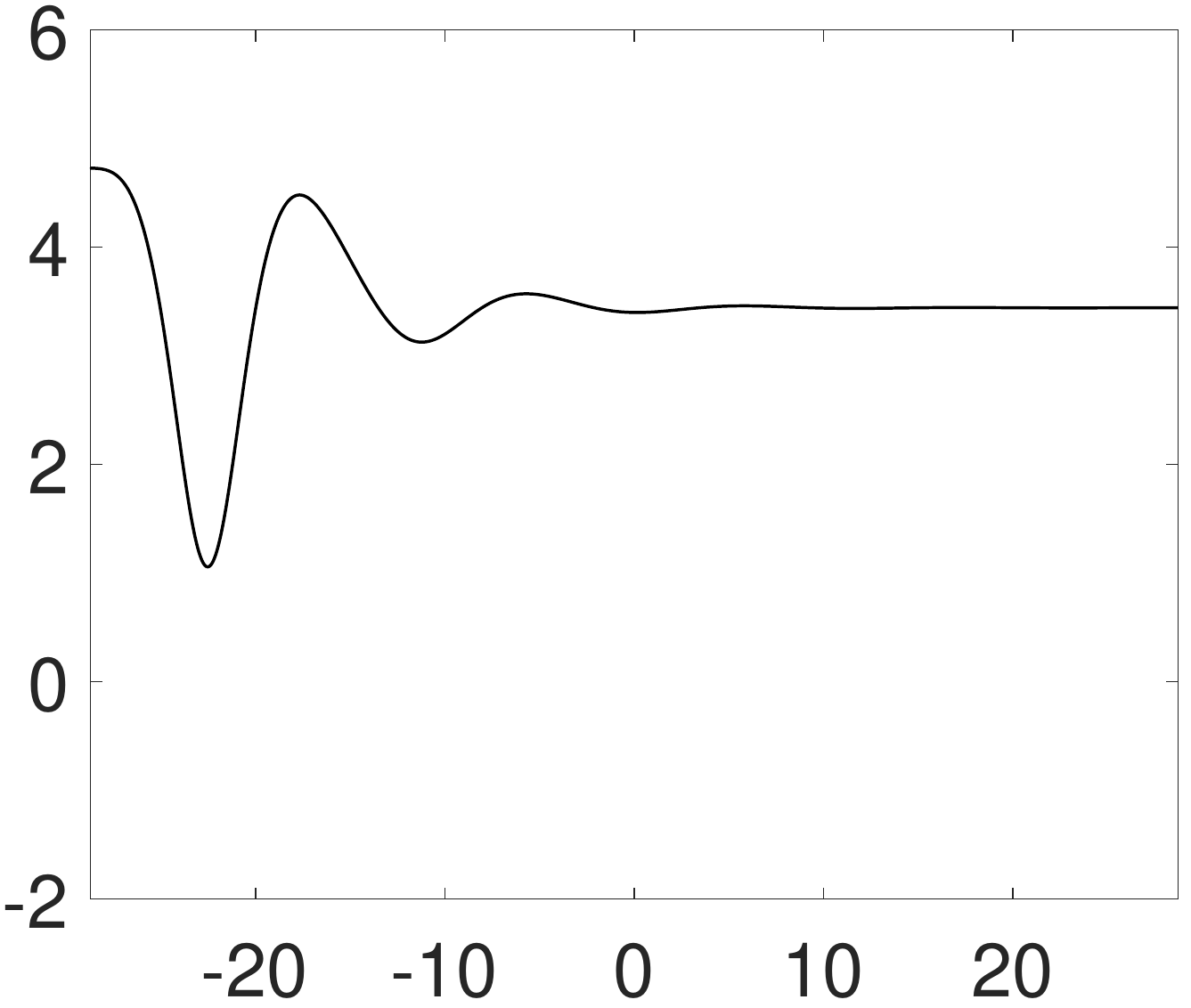}
\put(-5,40){\rotatebox{90}{$u_1$}}
\put(95,0){$x$}
\end{overpic}
\caption{}
\label{fig:schnak-sig06-sol-a}
\end{subfigure}
\begin{subfigure}{0.23\textwidth}
\centering
\begin{overpic}[width=\textwidth]{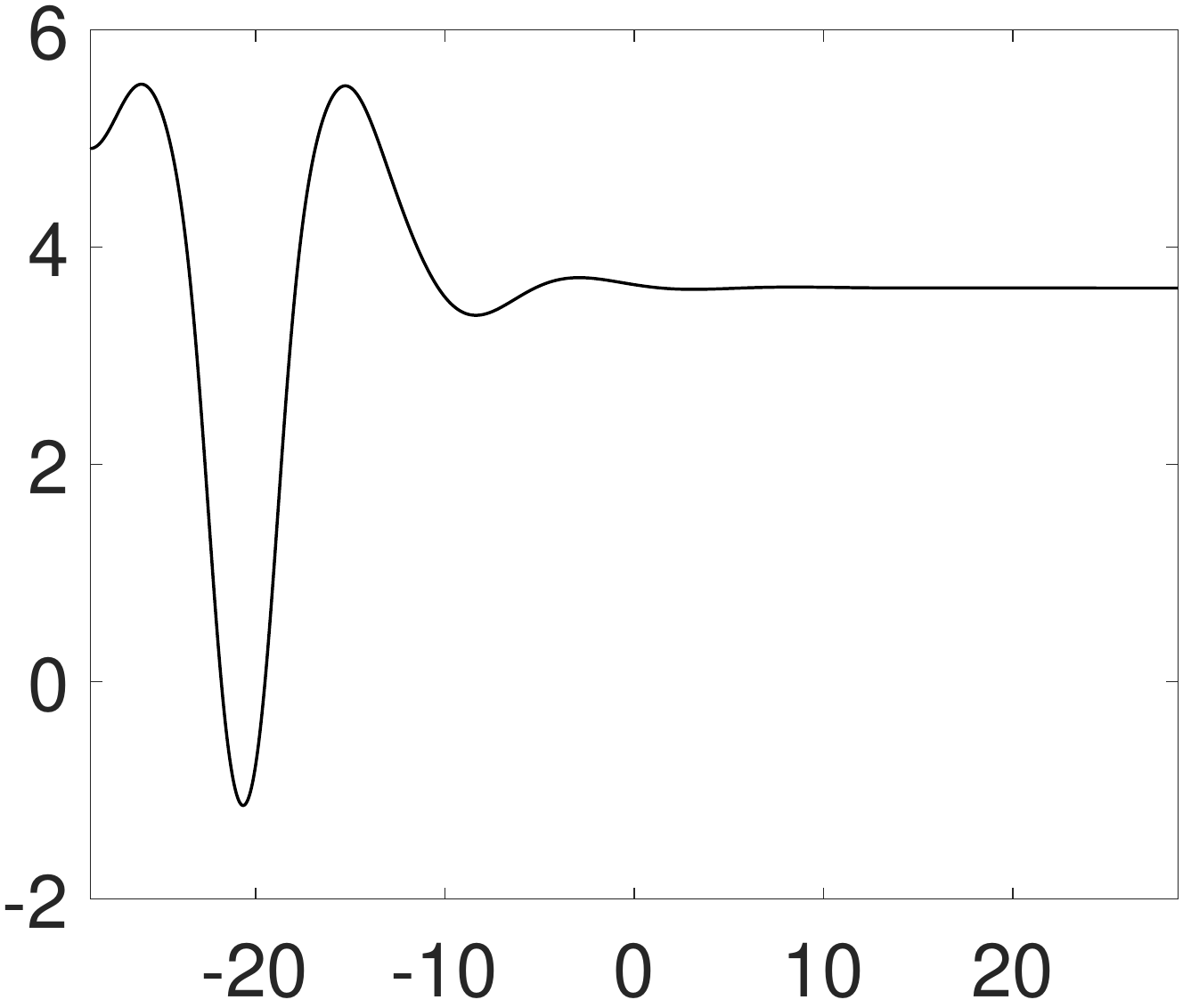}
\put(95,0){$x$}
\end{overpic}
\caption{}
\label{fig:schnak-sig06-sol-b}
\end{subfigure}
\hspace{-1cm}
\begin{subfigure}{0.23\textwidth}
\centering
\begin{overpic}[width=\textwidth]{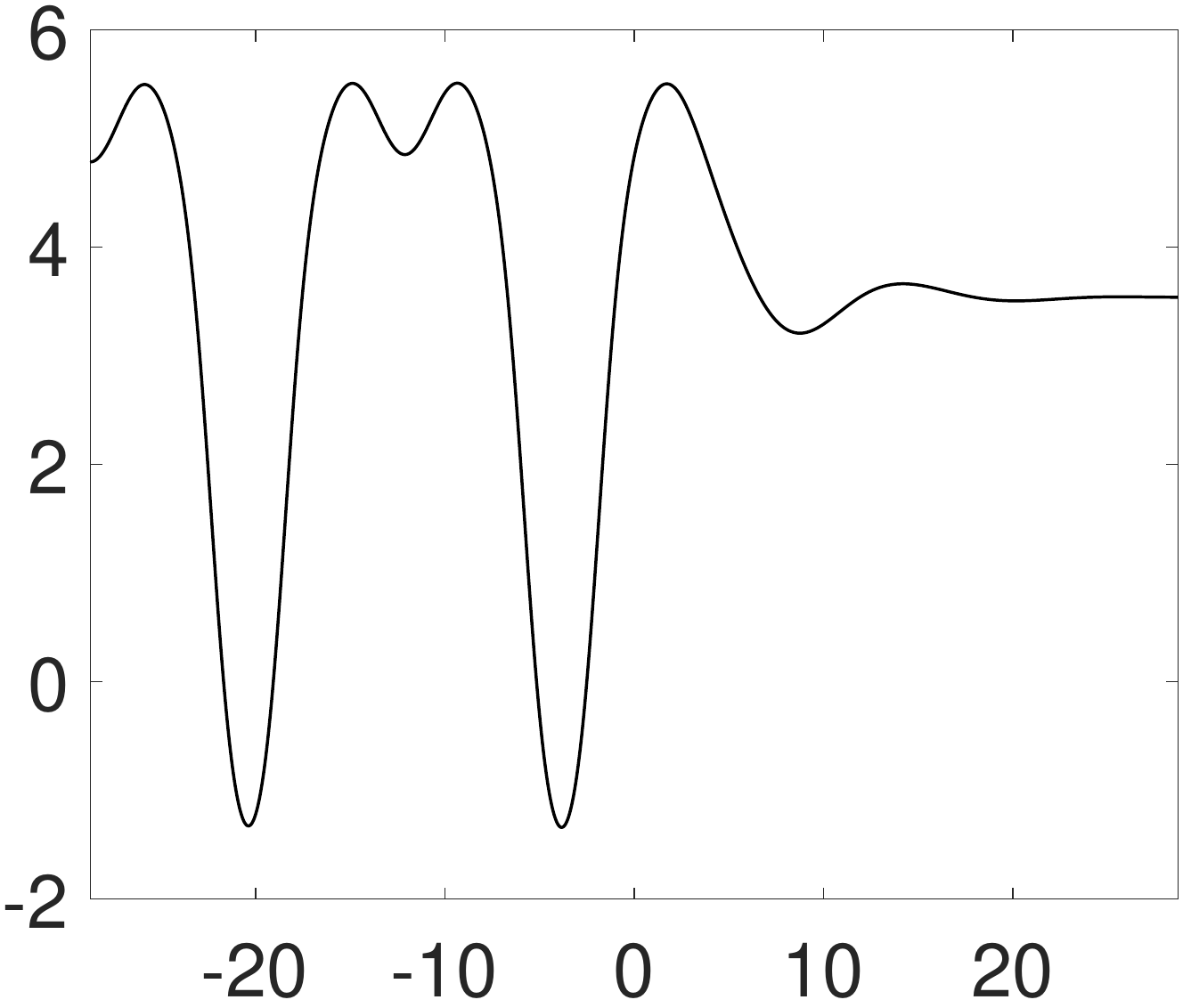}
\put(-5,40){\rotatebox{90}{$u_1$}}
\put(95,0){$x$}
\end{overpic}
\caption{}
\label{fig:schnak-sig06-sol-c}
\end{subfigure}
\begin{subfigure}{0.23\textwidth}
  \centering
\begin{overpic}[width=\textwidth]{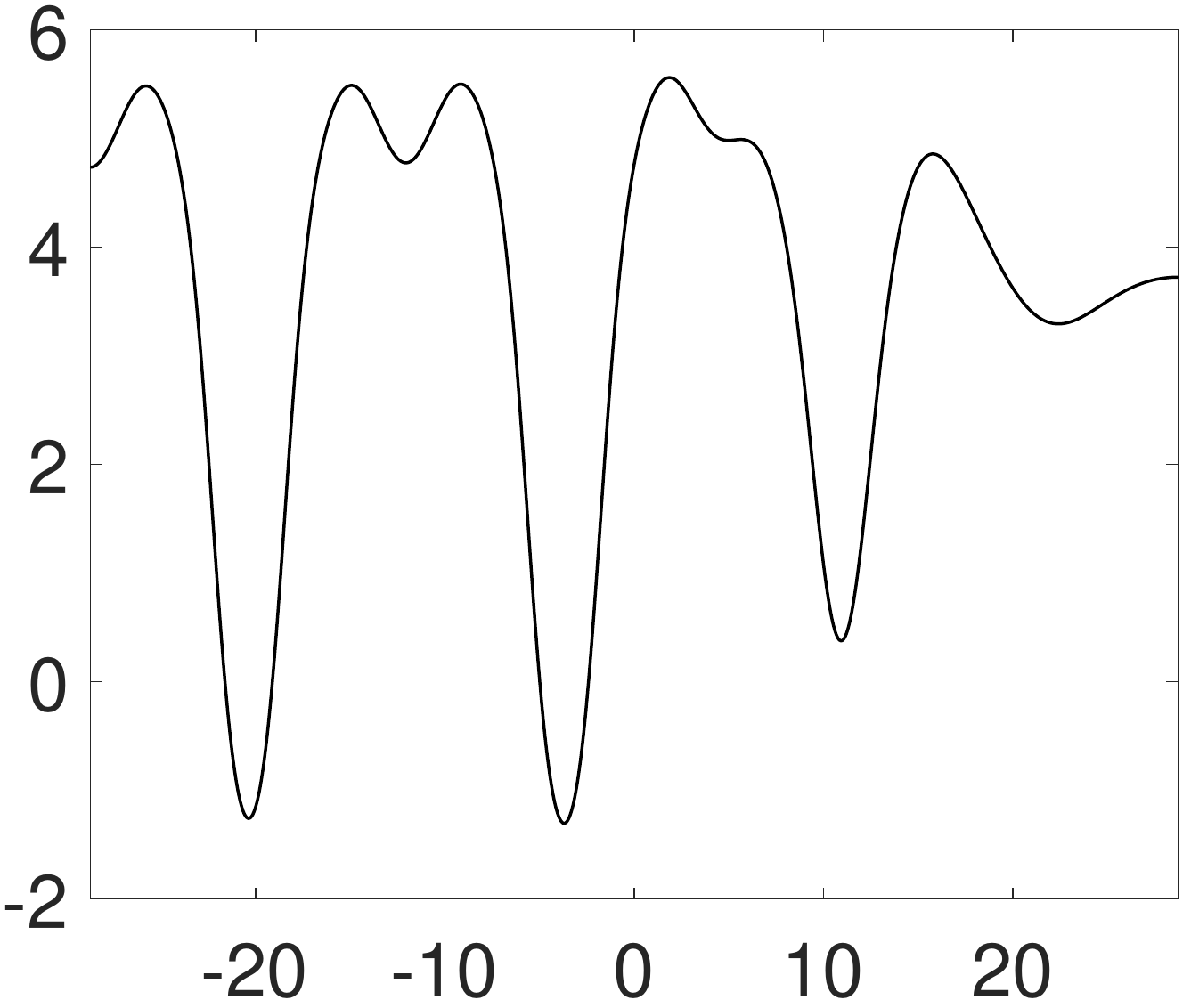}
\put(95,0){$x$}
\end{overpic}
  \caption{}
  \label{fig:schnak-sig06-sol-d}
\end{subfigure}
\hspace{-1cm}
\begin{subfigure}{0.23\textwidth}
  \centering
\begin{overpic}[width=\textwidth]{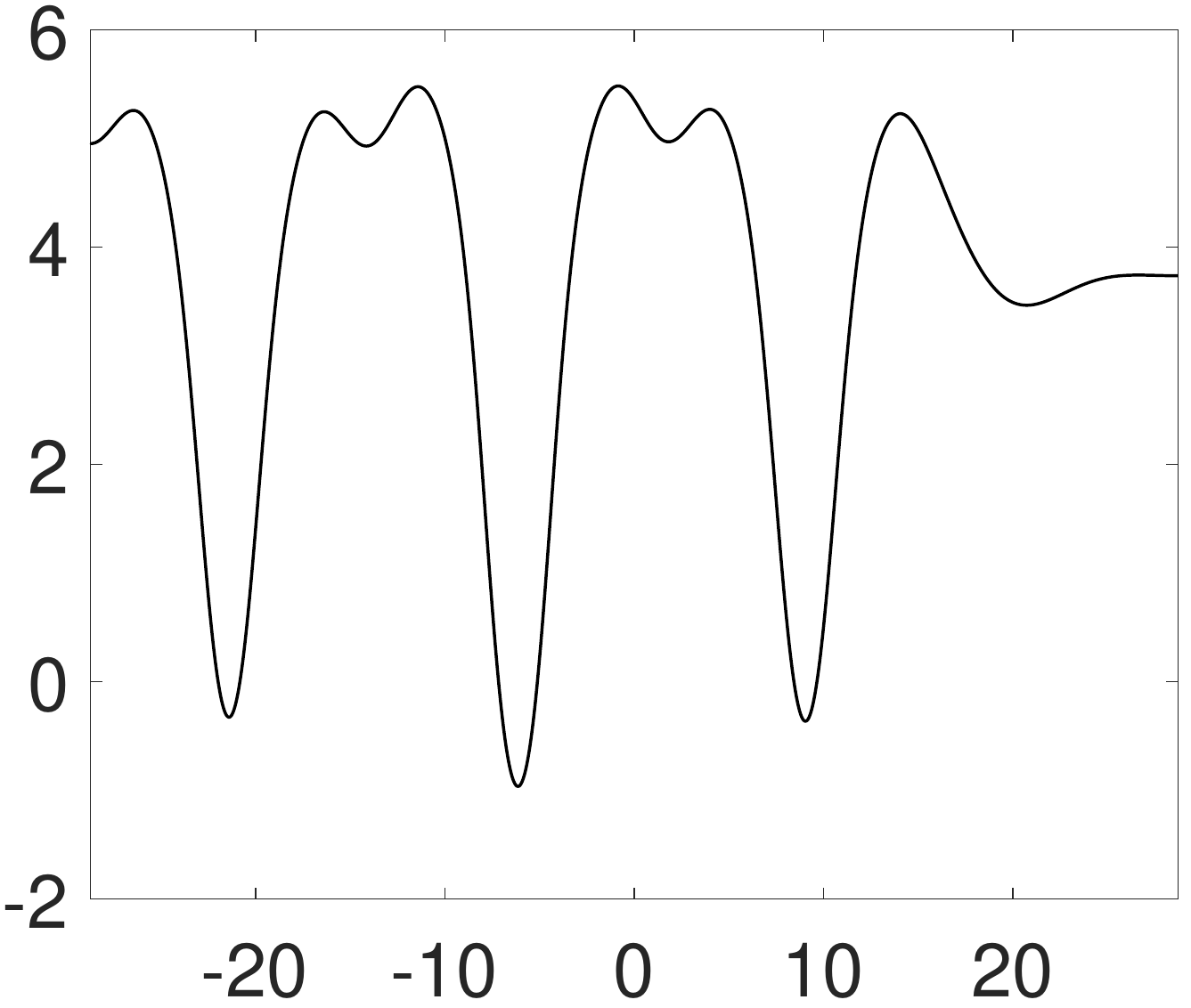}
\put(-5,40){\rotatebox{90}{$u_1$}}
\put(95,0){$x$}
\end{overpic}
  \caption{}
  \label{fig:schnak-sig06-sol-e}
\end{subfigure}
\begin{subfigure}{0.23\textwidth}
  \centering
\begin{overpic}[width=\textwidth]{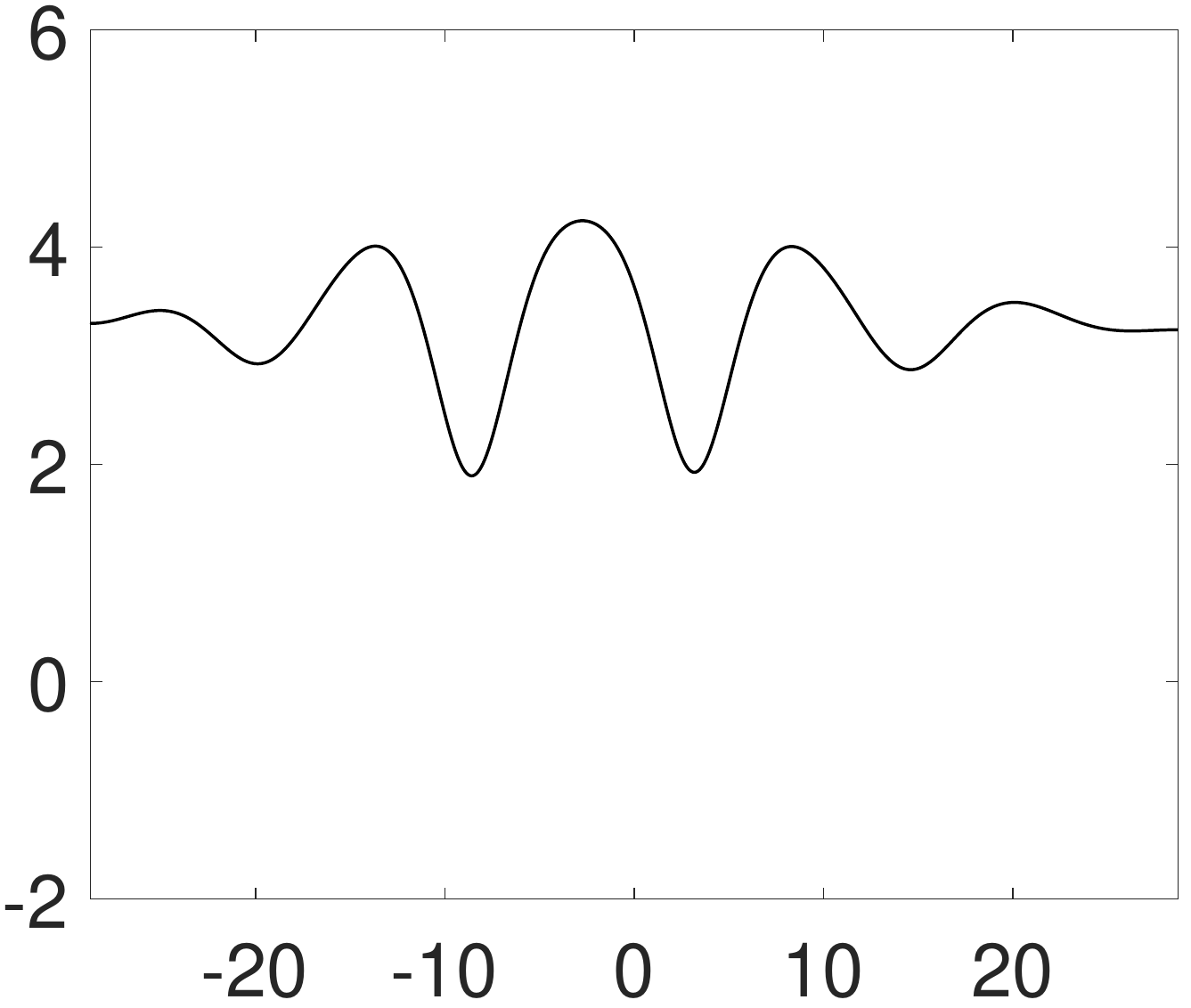}
\put(95,0){$x$}
\end{overpic}
  \caption{}
  \label{fig:schnak-sig06-sol-f}
\end{subfigure}
\end{multicols}
\vspace{-0.5cm}
\caption{Solutions along the first snaking branch of the fractional modified Schnakenberg system~\eqref{eq:Schnak-frac} ($\sigma=-0.6$) with fractional order $s=0.73$. (left) Snaking branch: purple denotes the first part of the snaking branch, while the second part is marked with light purple. (right) \subref{fig:schnak-sig06-sol-a}--\subref{fig:schnak-sig06-sol-c}: Solutions along the first part of the snaking branch. \subref{fig:schnak-sig06-sol-d}--\subref{fig:schnak-sig06-sol-f}: Solutions along the second part of the snaking branch.  (For interpretation of the references to color in this figure legend, the reader is referred to the web version of this article)}
\label{fig:solutions-along-sn2D}
\end{figure}

Furthermore, Figure~\ref{fig:solutions-along-sn2D} shows solution profiles along the snaking branch at $s=0.73$. As we had seen in Section~\ref{schnak}, they initially correspond to front solutions. Going up the snake, we expect the front to move ``forward'', that is gaining more and more oscillations until the domain is filled, to match the number of oscillations of the reconnecting branch. However, here, we see that the oscillations does not gain the fourth complete bump. In fact, at the upper fold point the front stops moving forward. In addition, unstable localized solutions appear after the fourth fold, which progressively present oscillations in the central part of the domain (Figure~\ref{fig:schnak-sig06-sol-f}). This change of behavior can be due to an interaction between the first snaking branch of front solutions with the snaking branch of homoclinic solutions originating (in the standard Schnakenberg system) from the second bifurcation point on $P5$~\cite{Uecker2014}.

One last observation has to be made looking at Figures \ref{fig:bif-sig06}--\ref{fig:solutions-along-sn2D}. Along the snaking branches the software detects Hopf bifurcations (indicated by a diamond marker in the figures) which are not present with the standard Laplacian. However, further investigations of the presence of time-periodic orbits and the branch switching at these Hopf points will be matter of future works.

\section{Conclusion and Outlook} \label{sec:conclusions}

In this paper we have successfully extended for the first time numerical continuation methods to nonlinear space-fractional PDEs. Specifically, in the context of \texttt{pde2path}, we have interfaced FEM with the fractional Laplacian via a discretization of the Balakrishnan formula.
This enabled us to investigate some of the effects of super-diffusion on the steady state bifurcation structure and solution profiles of the Allen--Cahn equation~\eqref{eq:AC-frac}, the Swift--Hohenberg equation~\eqref{eq:SH-frac} and the Schnakenberg system~\eqref{eq:Schnak-frac}. We highlight here the common effects that were observed, which contribute to a better understanding of the effects of the fractional Laplacian on the steady solutions of generic reaction--diffusion systems on bounded domains.
\begin{itemize}[leftmargin=0.3cm]
\item[-] Firstly, in all systems we have observed that the bifurcation points on the branch of homogeneous solutions accumulate at a precise value of the bifurcation parameter (at which a change of stability of the homogeneous solution occurs) as the fractional order $s$ tends to zero. This can be seen from equations  \eqref{eq:AC-bifurcation-points-s}, \eqref{eq:subsequent-bif} and \eqref{eq:bifurcations-schnak-frac}, and it is illustrated in Figures \ref{fig:AC-bifurcation-pts-1D}, \ref{fig:sh-bif1and2-shifting}, \ref{fig:b1andb2goonawalk} and \ref{fig:b1b2-sig06}. In particular, for the fractional Allen--Cahn all the bifurcation points vary and tend to the value $\mu=1$ as $s\to 0^+$. For the Swift--Hohenberg equation and the Schnackenberg system the first bifurcation point (where the change in stability occurs) remains fixed, while the others move towards $\mu=0$ and $\mu=\mu_c$ (formula \eqref{eq:mu_schnak}), respectively.
\item[-] Secondly, we have seen that, for a fixed value of the bifurcation parameter, the non-homogeneous stationary solutions to the three equations tend to have flatter peaks or valleys as the fractional order of the Laplacian decreases. This effect can be memorized as ``sharpening the teeth'' of the spatial solution profiles. We suppose that this effect comes from the fact that diffusion acts more ballistically as $s$ decreases, producing a larger spread and larger tails.
\item[-] Finally, in both systems featuring a snaking branch of localized solutions, that is the Swift--Hohenberg equation and the Schnakenberg system, we have observed that decreasing the fractional order of the Laplacian leads to a significant widening of the snaking, see Figures \ref{fig:sh-bif-firstImpression} and \ref{fig:bif-sig06}. This can be memorized as ``bloated snaking''. Furthermore, the re-connection properties of the snaking branch change as the fractional order is decreased.
\end{itemize}

Hence, our results pave the way for further analytical investigations of nonlinear reaction--diffusion systems involving fractional operators.\\
Our investigation focused on one-dimensional fractional problems, where matrices associated with the discretized problem can be kept sufficiently small. Thus, these problems can still be treated numerically on today's standard desktop computers. To discretize the fractional operator in higher dimensions poses additional numerical challenges. In this regard, the choice of interfacing the \texttt{pde2path} package with the fractional Laplacian via a discretization of the Balakrishnan formula makes the solver flexible and, in particular, suitable for higher dimensional problems involving adaptive meshes. To obtain accurate results the Balakrishnan formula involves the solution of a large number of linear systems of equations of significant size. Therefore, it might be interesting to check the performance of iterative linear system solvers in order to make the computation more efficient. In particular, a method based on Krylov subspaces has recently been implemented~\cite{dohr2019fem}, which could be helpful in future work. Of course, it would be interesting to use other discretizations of the spectral fractional Laplacian on bounded domains based on the integral formulation of the operator via the heat-semigroup formalism as proposed in \cite{cusimano2018discretizations} or based on the very weak FEM \cite{antil2018fractional}. Furthermore, the choice of parameters $k,\,n_+,\, n_-,$ appearing in the FEM discretization of the fractional Laplacian as well as the detailed study of the numerical accuracy could be addressed in future work. Since we have benchmarked the method already against analytical results, which gave excellent agreement with the numerics, we expect that the numerical analysis is also going to yield favorable results. \\
Beyond perspectives of computational and numerical optimizations, the preliminary results obtained within this work bring up fascinating issues, which will be matter of future work. For instance, the deformations and changes of the snaking branch in the Schnakenberg system for fractional order smaller than $0.7$ is an open issue. Another point of interest is the nature of secondary bifurcation points, the bifurcating branches and their fate as the fractional order $s$ goes to zero, as well as the presence of time-periodic orbits and the branch switching at the Hopf points as the fractional order becomes smaller. Furthermore, a natural extension at this point is the numerical investigation of these models in higher-dimensional domains. Regarding other models, it could be interesting to investigate if the fractional Laplacian can lead to the appearance of unexpected patterns (which are not shown with standard diffusion). In particular, looking towards the applications, we could also focus on models that arise in the literature, for which the application of the continuation software could lead to a better undertanding of the studied phenomena. Finally, implementing other discretizations for different definitions of the fractional operator \cite{bonito2018numerical, hofreither2019unified}, we will be able to contrast and compare them by looking at the bifurcation diagrams and stationary solutions.

In conclusion, the paper constitutes the meeting point between theoretical results, bifurcation theory, numerical analysis and scientific computing. Extending advanced continuation methods, such as implemented in \texttt{pde2path}, allowed us to treat a much larger class of systems. In particular, we gained a new tool to investigate bifurcations and pattern formation in fractional PDEs. Thus, we established a novel interaction between continuation techniques and fractional PDEs. This point of view will shed new light on the field as it allows for global exploration of nonlinear systems involving the fractional Laplacian. We can now systematically apply our approach to investigate a large class of problems in many fields (biology, epidemiology, population dynamics). Furthermore, the paper can be viewed as a key contribution towards a dynamical system approach to fractional reaction--diffusion systems.

\bigskip\noindent
\textbf{Acknowledgments:}
CK and CS have been supported by a Lichtenberg Professorship of the VolkswagenStiftung. CK also acknowledges partial support of the EU within the TiPES project funded by the European Unions Horizon 2020 research and innovation programme under grant agreement No. 820970. CS has received funding from the European Union's Horizon 2020 research and innovation programme under the Marie Sk\l odowska--Curie grant agreement No. 754462. Partial support by the Italian National Group of Mathematical Physics (GNFM--INdAM) is also gratefully acknowledged by CS. NE acknowledges the support of Deutschlandstipendium. The authors acknowledge the financial support of the University of Graz.

\bibliographystyle{plainurl}
\bibliography{bibliography}

\appendix
\section{The MatLab code}\label{app:code}
We provide here the Matlab code which implements the FEM discretization of the fractional Laplacian presented in Algorithm \ref{alg:FractionalLaplacian} in Section \ref{sec:numerics}. We do not explain in detail the basic setup of \texttt{pde2path}; see \cite{Uecker2014, uecker_pde2path_2014,rademacher2018oopde} for a complete overview of the continuation software for beginner users.

In order to numerically study fractional reaction--diffusion equation of the form~\eqref{eq:reaction-diffusion-ss-frac} in \texttt{pde2path}, we must adapt the standard setup providing a new function \texttt{FractionalLaplacian} which will be called by the function \texttt{oosetfemops}. We recall that the routine \texttt{oosetfemops} generates the matrices~\eqref{eq:matrixtot} appearing in the algebraic system~\eqref{eq:algebraic_system_for_G}. In the fractional case, this file differs from the standard setting only for the line by which the FEM stiffness matrix $K$, stored in the Matlab variable \texttt{p.mat.K}, is replaced by $-MK_s$, the matrix approximation of the fractional Laplacian multiplied with the mass matrix. Else the usage of \texttt{pde2path} stays the same for the basic continuation calls. A fully exacutable code reproducing Figures \ref{fig:AC-bifurcations-1D}, \ref{fig:AC-solutions-mu2-1D}, \ref{fig:sh-bif-firstImpression} and \ref{fig:bif-sig06} can be found at \cite{VideoFolder}. 

The new function \texttt{FractionalLaplacian} implementing Algorithm \ref{alg:FractionalLaplacian} is presented in Listing~\ref{lst:FractionalLaplacian}. The required input arguments are the standard stiffness and mass matrices \texttt{K} and \texttt{M}, which must encode the proper boundary conditions, the number of mesh points \texttt{np}, the (uniform) mesh size \texttt{h} and the fractional order \texttt{s}. The output of the function is the matrix approximation of the fractional Laplacian \texttt{Ks}.

Few remarks need to be made concerning the Matlab implementation. The first comment regards the choice of parameters $\kappa$, $n_+$ and $n_-$ in the quadrature approximation~\eqref{eq:quadrature_approx}, denoted by \texttt{k, Np, Nm} in Listing \ref{lst:FractionalLaplacian}. Given a spacial discretization with mesh size $h$, the value of $\kappa$ must be chosen such that the sinc quadrature error balances with the finite element error. Following \cite[Corollary 2]{dohr2019fem} where a similar discretization method was used, we choose
\begin{equation*}
\kappa = \frac{1}{|log(h)|},
\end{equation*}
corresponding to Line \ref{line:kvalue} in Listing \ref{lst:FractionalLaplacian}. Then, again following \cite{dohr2019fem}, the values $n_+$ and $n_-$ are chosen proportional to $1/\kappa^2$:
\begin{equation*}
n_{+} \coloneqq \left\lceil\frac{\pi^2}{4(1-s)\kappa^2}\right\rceil\quad \quad \text{and} \quad \quad n_{-} \coloneqq \left\lceil\frac{\pi^2}{4s\kappa^2}\right\rceil,
\end{equation*}
and they are assigned in Lines \ref{line:Npvalue}, \ref{line:Nmvalue} in Listing \ref{lst:FractionalLaplacian}.

The second comment regards the efficiency and correctness of the computation. For each column of $K_s$ we need to solve $n_{+}+n_{-}+1$ systems of equations (see equation~\eqref{eq:linsys}) of the type
 \begin{equation*}
(\txte^{\kappa l}M + K) y = K \hat{e}_i,
\end{equation*}
where $-n_{-} \leq l \leq n_{+}$. Since all columns as well as all systems are independent, ideally, we would like to parallelize both loops in Lines \ref{line:parfor} and \ref{line:innerloop}. However, nested parallelism is not allowed in Matlab. Therefore, we only parallelize the outer loop (using \texttt{parfor} in Line \ref{line:parfor}).

Finally, when $\txte^{\kappa l}$ becomes sufficiently small, we have that $(\txte^{\kappa l}M + K)\approx K$. Since $K$ is singular, there exist infinitely many solutions to systems of the form $Kv=z$. In the Matlab implementation, the typical \textit{backslash} operator would produce \textit{one} solution but we would not know \textit{which} solution it produces. In order to have a more robust routine, we use \texttt{lsqminnorm} which computes the minimum norm least-squares solution to the system using the complete orthogonal decomposition of $(\txte^{\kappa l}M + K)$ (see Line \ref{line:lswminnorm} in Listing \ref{lst:FractionalLaplacian}). 

\begin{Listing}
\lstset{language=Matlab,
    breaklines=true,
    escapechar=|,
    keywordstyle=\color{blue},
    identifierstyle=\color{black},
    stringstyle=\color{mylilas},
    commentstyle=\color{mygreen},
    showstringspaces=false,
    numbers=left,
    numberstyle={\color{black}},
    numbersep=9pt,
    emph=[1]{pi,log,zeros,exp,sin},emphstyle=[1]\color{black},
    emph=[2]{parfor},emphstyle=[2]\color{blue},
}
\begin{lstlisting}
function Ks = FractionalLaplacian(K,M,np,h,s)
% Ks discretization of the fractional operator
% K,M standard stiffness and mass matrices
% np number of mesh points
% h meshsize
% s fractional order

% parameters
k = -1/log(h) |\label{line:kvalue}|
Np = ceil(pi^2/(4*(1-s)*k^2));|\label{line:Npvalue}|
Nm = ceil(pi^2/(4*s*k^2));|\label{line:Nmvalue}|
coeff = -k*sin(s*pi)/pi;

% fractional Laplacian
Ks = zeros(np,np);
parfor i = 1: np |\label{line:parfor}|
    t = zeros(np,1);
    t(i,1) = 1; % test vector
    col = zeros(np,1);
    % compute the Balakrishnan sum
    z = K*t;
    for l =[-Nm:Np] |\label{line:innerloop}|
        alpha = exp(k*l);
        A = alpha*M + K;
        v = lsqminnorm(A,z); |\label{line:lswminnorm}|
        col = col + exp(s*k*l)*v;
    end
    col = coeff*col;
    % set the i-th colomn of Ks
    Ks(:,i) = col;
end
\end{lstlisting}
\caption{\texttt{FractionalLaplacian.m}}
\label{lst:FractionalLaplacian}
\end{Listing}
\end{document}